\documentclass[11pt]{article}
\usepackage{amsmath, amsthm, amsfonts, amssymb, color}
\usepackage{mathrsfs, bbm}

\newtheorem{thm}{Theorem}[section]
\newtheorem{cor}[thm]{Corollary}
\newtheorem{lem}[thm]{Lemma}
\newtheorem{remark}[thm]{Remark}
\newtheorem{prp}[thm]{Proposition}

\theoremstyle{definition}

\definecolor{wco}{rgb}{0.5,0.2,0.3}

\numberwithin{equation}{section} \theoremstyle{remark}

\def\1{{\mathbbm 1}}
\def\non{\nonumber \\}
\def\wt{\widetilde}
\def\eps{\varepsilon}

\def\R{\mathbb R}

\def\E{\mathbb E}

\def\L{\mathcal L}
\def\S{\mathcal S}
\def\J{\mathcal J}

\def\M{\mathcal M}
\def\wh{\widehat}

\def\<{\langle}
\def\>{\rangle}

\def\LL{\Lambda}

\def\pf{\noindent{\bf Proof.} }

 \def\beq{\begin{equation}}

 \def\P{\mathbb P} 
  \def\ee{\varepsilon}

\begin{document}
\bibliographystyle{plain}

\title{\Large \bf Perturbation  by Non-Local Operators}

\author{
{\bf Zhen-Qing Chen}\thanks{Research partially supported
by NSF Grant DMS-1206276, and NNSFC Grant 11128101.}
  \quad and
\quad {\bf Jie-Ming Wang}\thanks{Corresponding author.
Research partially supported  by NNSFC Grant 11401025.}
}
\date{(September 28, 2016)}

\maketitle

\begin{abstract}
Suppose that $d\ge 1$ and $0<\beta<\alpha<2$. We establish
the existence and uniqueness of the fundamental solution
$q^b(t, x, y)$ to 
a class of (typically nonsymmetric)
non-local operators $\L^b=\Delta^{\alpha/2}+\S^b$,
where
 $$
 \S^bf(x):={\cal A}(d, -\beta)
\int_{\R^d} \left( f(x+z)-f(x)- \nabla f(x) \cdot
z\1_{\{|z|\leq 1\}} \right) \frac{b(x, z)}{|z|^{d+\beta}}dz
$$
and $b(x, z)$ is a bounded measurable function on $\R^d\times \R^d$
with $b(x, z)=b(x, -z)$ for $x, z\in \R^d$. Here ${\cal A}(d, -\beta)$
is a normalizing constant so that $\S^b=\Delta^{\beta/2}$
when $b(x, z)\equiv 1$. We show that if
$b(x, z) \geq -\frac{{\cal A}(d, -\alpha)}
{{\cal A}(d, -\beta)}\, |z|^{\beta -\alpha}$, then
$q^b(t, x, y)$ is a strictly positive continuous function
and it uniquely determines a conservative Feller process $X^b$,
which has strong Feller property.  The Feller process $X^b$ is the unique
solution to the martingale problem of $(\L^b, {\cal S} (\R^d))$,
where ${\cal S}(\R^d)$ denotes the space of tempered functions on $\R^d$.
Furthermore, sharp two-sided estimates on $q^b(t, x, y)$ are derived.
In stark contrast with the gradient perturbations,
these estimates exhibit different behaviors for different types of $b(x, z)$.
The model considered in this paper
contains the following as a special case.
Let $Y$ and $Z$ be (rotationally) symmetric
$\alpha$-stable process and symmetric $\beta$-stable processes
on $\R^d$, respectively,
that are independent to each other.
Solution to
stochastic differential equations $dX_t=dY_t + c(X_{t-})dZ_t$
has infinitesimal generator $\L^b$ with $b(x, z)=| c(x)|^\beta$.
\end{abstract}



\section{Introduction}

Let  $d\geq 1$ be an integer and $0< \beta< \alpha <2$.
For integer $k\geq 1$, denote by $C^k_b(\R^d)$ (resp. $C^k_c (\R^d)$)
the space of continuous functions on $\R^d$ that have bounded
continuous partial derivatives up to order $k$ (resp. the space
of continuous functions on $\R^d$ with compact support
that have continuous  partial derivatives up to order $k$).
Recall that a stochastic process $Y=(Y_t, \P_x, x\in \R^d)$
on $\R^d$ is called a (rotationally)
symmetric $\alpha$-stable process on $\R^d$ if it is
a L\'evy process having
$$ \E_x \left[ e^{i\xi \cdot (Y_t-Y_0)}\right] =e^{-t |\xi|^\alpha}
\qquad \hbox{for every } x , \xi \in \R^d.
$$
Let $\wh f(\xi):=\int_{\R^d} e^{i\xi \cdot x } f(x) dx$
 denote the Fourier transform of a function $f$ on $\R^d$.
The fractional Laplacian $\Delta^{\alpha/2}$ on $\R^d$ is defined as
\begin{equation}\label{e:1.2n}
 \Delta^{\alpha/2} f(x)= \int_{\R^d} \left( f(x+z)-f(x)- \nabla f(x)
 \cdot  z\1_{\{|z|\leq 1\}}  \right) \frac{{\cal A}(d, -\alpha)}{|z|^{d+\alpha}}dz
\end{equation}
for   $f\in C_b^2(\R^d)$. Here  ${\cal A}(d, -\alpha)=\Gamma((d+\alpha)/2)/(2^{-\alpha}\pi^{d/2}|\Gamma(-\alpha/2)|),$
which is the normalizing constant so that
$\widehat {\Delta^{\alpha/2} f} (\xi)
= -|\xi|^\alpha \widehat f (\xi) $.
Hence $\Delta^{\alpha/2}$ is the infinitesimal generator
for the symmetric $\alpha$-stable process on $\R^d$.

Throughout this paper, $b(x, z)$ is
a  real-valued bounded function
on $\R^d\times \R^d$ satisfying
\begin{equation}\label{e:b}
b(x, z)=b(x, -z) \qquad \hbox{for every } x, z\in \R^d.
\end{equation}
This paper is concerned with the existence, uniqueness
and   
{\it sharp two-sided estimates} 
on the ``fundamental solution" of
the following non-local operator  on $\R^d$,
$$
\L^b f(x)=\Delta^{\alpha/2} f(x) + \S^b f(x),
\quad f\in C^2_b(\R^d),
$$
where
 \begin{equation}\label{e:1.1}
\S^bf(x):={\cal A}(d, -\beta)
\int_{\R^d} \left( f(x+z)-f(x)- \nabla f(x) \cdot
z \1_{\{|z|\leq 1\}}  \right) \frac{b(x, z)}{|z|^{d+\beta}}dz.
\end{equation}
We point out that
since
 $b(x, z)$ satisfies condition \eqref{e:b},
the truncation $|z|\leq 1$ in \eqref{e:1.1} can be replaced
by $|z|\leq \lambda$
for any $\lambda >0$; that is, for every $\lambda>0$,
\begin{equation}\label{e:1.4a}
\S^bf(x) ={\cal A}(d, -\beta)
\int_{\R^d} \left( f(x+z)-f(x)-\langle\nabla f(x),
z \> \1_{\{|z|\leq \lambda\}} \right) \frac{b(x, z)}{|z|^{d+\beta}}dz.
\end{equation}
In fact, under condition \eqref{e:b},
\begin{eqnarray}\label{e:1.5a}
\S^bf(x)&=& {\cal A}(d, -\beta) \ {\rm p.v.}
\int_{\R^d} \left( f(x+z)-f(x)\right)
\frac{b(x, z)}{|z|^{d+\beta}}dz \nonumber \\
&:=& \mathcal{A}(d,-\beta) \,
\lim_{\varepsilon\to 0}\int_{\{z\in \R^d:
|z|>\varepsilon\}}\left(f(x+z)-f(x)\right)\frac{b(x,z)}{|z|^{d+\beta}}\,dz.
\end{eqnarray}
 
Condition \eqref{e:b}
allows us to reduce general  bounded measurable function $b$ on $\R^d\times
\R^d$ to the situation
where $\| b\|_\infty$
is sufficiently small through a scaling argument
(see \eqref{e:4.2} and Lemma \ref{L:3.5b}).
The operator $\L^b$ is in general non-symmetric.
Clearly,  $\L^b=\Delta^{\alpha/2}$ when $b \equiv 0$ and
$ \L^b=\Delta^{\alpha/2}+\Delta^{\beta /2}$ when $b\equiv 1$.

We are led to the study of this non-local operator $\L^b$ by
 the consideration of the following
stochastic differential equation (SDE) on $\R^d$:
\begin{equation}\label{e:sde}
dX_t =dY_t + c(X_{t-}) dZ_t,
\end{equation}
where $Y$ is a  symmetric $\alpha$-stable process
on $\R^d$ and $Z$ is an independent symmetric $\beta$-stable
process with $0<\beta<\alpha$.
Such SDE arises naturally in applications when there are more than one
sources of random noises.
When $c$ is a bounded Lipschitz function
on $\R^d$,  it is easy to show
using Picard's iteration method
that for every $x\in \R^d$,
SDE \eqref{e:sde} has a  unique strong solution with
$X_0=x$. We denote the law of such a solution by  $\P_x$.
The collection of the solutions $(X_t, \P_x, x\in \R^d)$
 forms a strong Markov process $X$ on $\R^d$.
 Using Ito's formula, one concludes that the
infinitesimal generator of $X$ is  $\L^b$ with
 $b(x, z)=|c(x)|^\beta$ and so in this case
 $X$ solves the martingale problem for
 $(\L^b, C^2_b(\R^d))$.
The following questions arise naturally:
does  the Markov process $X$
have a transition density function? If so, what
is its sharp two-sided estimates?
Is there a solution to the martingale problem
for $\Delta^{\alpha/2}+|c(x)|^\beta \Delta^{\beta/2}$
when $c$ is not Lipschitz
continuous?
We will address these questions for the more general operator
$\L^b$ in this paper.

Heat kernel analysis is an important subject in analysis and in probability theory,
as heat kernel encodes all the information about the corresponding infinitesimal 
generator and  the corresponding Markov processes. 
Since explicit formula can only be derived in some very special and limited  cases,
the main focus of the heat kernel analysis is on its sharp estimates. While it is relatively 
easy to get some crude bounds, obtaining sharp two-sided bounds on the heat kernel
is typically quite  challenging. It requires deep understanding of
the corresponding generator. 
 Heat kernel estimates for discontinuous Markov processes
have been
under intense study recently.
 Most results obtained so far are
 mainly for symmetric Markov processes. See \cite{C} for a recent survey.
 It is well known that the study of non-symmetric operators requires   different
 approaches and techniques  than that for symmetric operators. 
Results of this paper can also be viewed as an attempt in
establishing heat kernel estimates for non-symmetric discontinuous
Markov processes. For example,   Corollary \ref{C:1.4} and 
Theorem \ref{T:1.4}  can be viewed
as the non-symmetric analogy, though in a restricted setting,
of the two-sided heat kernel estimates for symmetric stable-like
processes and mixed stable-like processes established in \cite{CK}
and \cite{CK2}, respectively. See Remark \ref{R:1.4} below for more information
on heat kernel analysis.  

\medskip

For $a\geq 0$, denote by   ${p}_a(t,x, y)$ the
fundamental function   of   $\Delta^{\alpha/2}+a\Delta^{\beta/2}$
(or equivalently, the transition density function of the L\'evy process
$Y_t+a^{1/\beta} Z_t$).  Clearly,  $  p_a(t, x, y)$
is a function of $t$ and $x-y$, so sometimes we
also  write it as $p_a (t, x-y)$.
It is known (see \eqref{e:h} of  Section 2 for details)
 that on $(0, \infty)\times \R^d \times \R^d$,
\begin{eqnarray}
p_0 (t, x, y)  &\asymp& t^{-d/\alpha} \wedge \frac{t}{|x-y|^{d+\alpha}} ,
\label{e:1.5} \\
  p_a(t, x, y) &\asymp& \left(t^{-d/\alpha}\wedge (at)^{-d/\beta}\right)
 \wedge \left( \frac{t}{|x-y|^{d+\alpha}} + \frac{at}{|x-y|^{d+\beta}}\right). \label{e:1.6}
\end{eqnarray}
Here for two non-negative functions $f$ and $g$,
the notation $f\asymp g$ means that there is a constant $c\geq 1$ so that $c^{-1} f\leq g\leq cf$ on
their common domain of definitions. For real numbers $a, c\in \R$,
we use $a\vee c$ and $a\wedge c$ to denote $\max\{a, c\}$ and $\min \{a, c\}$, respectively.
We point out that the comparison constants in \eqref{e:1.6} is independent
of $a>0$; see \eqref{e:h} in Section 2.
Using the observation that $a\wedge b \asymp \frac{ab}{a+b}$,
one concludes from \eqref{e:1.5} that
\begin{equation} \label{e:1,8}
p_0(t, x, y) \asymp \frac{t}{(t^{1/\alpha} + |x-y|)^{d+\alpha}}
\qquad \hbox{on } (0, \infty) \times \R^d \times \R^d.
\end{equation}
Note that $(at)^{-d/\beta} \geq t^{-d/\alpha}$
whenever $0<t\leq a^{-\alpha/(\alpha-\beta)}$. Thus for every $k>0$,
\begin{equation}\label{e:1.7a}
  p_a(t, x, y) \asymp t^{-d/\alpha}
\wedge \left( \frac{t}{|x-y|^{d+\alpha}} + \frac{at}{|x-y|^{d+\beta}}\right)
\qquad \hbox{on } (0, k a^{-\alpha/(\alpha-\beta)}]\times \R^d \times \R^d,
\end{equation}
with the comparison constants depending only on $d$, $\alpha$, $\beta$
and $k$.

\medskip

Since $\L^b=\Delta^{\alpha/2}+\S^b$ is a lower order perturbation of $\Delta^{\alpha/2}$ by $\S^b$, heuristically the fundamental solution
(or kernel) $q^b(t, x, y)$ of $\L^b$ should satisfy the following
Duhamel's formula:
\begin{equation}\label{e:1.7n}
q^b(t, x, y)=p_0(t, x, y)+\int_0^t \int_{\R^d} q^b (t-s, x, z)
\S^b_z p_0 (s, z, y)dz ds
\end{equation}
for $t>0$ and $x, y\in \R^d$.
Here the notation $S^b_z p_0 (s, z, y)$ means the non-local operator $\S^b$
is applied to the function $z\mapsto p_0(s, z, y)$. Similar notation will also be used for other operators, for example, $\Delta^{\alpha/2}_z$.
Applying \eqref{e:1.7n} recursively, it is reasonable to
conjecture that $\sum_{n=0}^\infty q^b_n(t, x, y)$, if convergent,
 is a solution to \eqref{e:1.7n}, where $q^b_0(t, x, y):=p_0(t, x, y)$
 and
\begin{equation}\label{e:qn}
q^b_n (t, x, y):=\int_0^t \int_{\R^d} q^b_{n-1} (t-s, x, z)
\S^b_z p_0(s, z, y) dz ds \quad \hbox{for } n\geq 1.
\end{equation}

For each bounded function $b(x,z)$ on $\R^d\times\R^d$ and
$\lambda >0$,  define
\begin{equation}\label{b1}
m_{b, \lambda} = {\rm essinf}_{x, z\in \R^d, |z|>\lambda} b(x, z)
\quad \hbox{and} \quad
 M_{b, \lambda} = {\rm esssup}_{x, z\in \R^d, |z|>\lambda} |b(x, z)| .
\end{equation}
The followings are the main results of this paper.

\begin{thm}\label{T:1.1}
For every bounded   function $b$ on
$\R^d\times \R^d$ satisfying condition \eqref{e:b},
there is a unique  continuous function $q^b(t,x,y)$
 on $(0, \infty)\times \R^d\times \R^d$
  that   satisfies \eqref{e:1.7n}  on
 $(0, \ee]\times \R^d \times \R^d$
 with $|q^b(t, x, y)| \leq c   p_1(t, x, y)$ on
 $(0, \ee]\times \R^d \times \R^d$ for some $\ee, c >0$, and that
\begin{equation}\label{e:1.4}
\int_{\R^d} q^b(t, x, y) q^b(s, y, z) dy =q^b(t+s, x, z)
\quad \hbox{for every }  t, s >0 \hbox{ and } x, z\in \R^d.
\end{equation}
Moreover, the following holds.
\begin{description}
\item{\rm (i)}
There is a constant $A_0=A_0(d, \alpha, \beta)>0$
  so that
$q^b(t, x, y)=\sum_{n=0}^\infty q^b_n(t, x, y)$
on \hfill \break
$(0,  ( A_0/\|b\|_\infty)^{\alpha/(\alpha -\beta)}]\times \R^d \times \R^d$, where
$q^b_n(t, x, y)$ is defined by \eqref{e:qn}.

\item{\rm (ii)}
$q^b(t,x,y)$ satisfies the Duhamel's formula \eqref{e:1.7n} for all $t>0$ and $x,y\in\R^d$. Moreover, $\S^b_x q^b(t, x, y)$ exists pointwise in the sense of \eqref{e:1.5a} and
\begin{equation}\label{e:Du}
q^b(t, x, y)=p_0(t, x, y)+\int_0^t \int_{\R^d} p_0 (t-s, x, z)
\S^b_z q^b (s, z, y)dz ds
\end{equation}
for $t>0$ and $x, y\in \R^d$.

\item{\rm (iii)}
 For each $t>0$ and $x\in \R^d$,
$\int_{\R^d} q^b(t, x, y) dy =1$.

\item{\rm (iv)} For every $f\in C^2_b (\R^d)$,
$$
T^b_t f(x)-f(x)=\int_0^t T^b_s \L^b f(x) ds,
$$
where $T^b_t f(x)=\int_{\R^d} q^b(t, x, y) f(y) dy$.

\item{\rm (v)} Let   $A>0$ and $\lambda >0$.
There is a  positive constant $C =C (d,\alpha,\beta,  A, \lambda )
\geq 1$
so that
for any $b$ satisfying \eqref{e:b} with $\| b \|_\infty \leq   A$,
\begin{equation}\label{e:1.7}
|q^b(t,x,y)| \leq C e^{Ct} {p}_{M_{b, \lambda}}
(t,x,y)  \quad \hbox{on }
(0, \infty) \times \R^d\times \R^d.
\end{equation}

\end{description}
\end{thm}

\medskip
 We remark that estimate \eqref{e:1.7} allows one to get sharper bound on $|q^b(t, x, y)|$ by selecting optimal $\lambda>0$.
When $Z_t$ is the deterministic process $t$ and $c$ is an $\R^d$-valued
bounded Lipschitz function on $\R^d$, the solution of \eqref{e:sde}
is a symmetric $\alpha$-stable process with drift.
Its infinitesimal generator is $\Delta^{\alpha/2}+c (x) \nabla$.
Existence of integral kernel to $\Delta^{\alpha/2}+c (x) \nabla$
and its estimates have been studied recently in \cite{BJ}
(in fact, $c$ there can be an $\R^d$-valued function in certain Kato class).

Unlike the  gradient perturbation for $\Delta^{\alpha/2}$, in general
the kernel $q^b(t, x, y)$ in Theorem \ref{T:1.1} can
take negative values.
For example, this is the case when $b\equiv -1$, that is,
when $\L^b=\Delta^{\alpha/2}-\Delta^{\beta /2}$,
according to the next theorem.
Observe that
$$ \L^bf (x)=\int_{\R^d} \left( f(x+z)-f(x)-\langle\nabla f(x), z \>
\1_{\{|z|\leq 1 \}} \right) j^b(x, z)dz,
$$
where
\begin{equation}\label{e:jb}
j^b(x, z)= \frac{{\cal A}(d, -\alpha)}{|z|^{d+\alpha}}
\left( 1+ \frac{{\cal A}(d, -\beta)}{{\cal A}(d, -\alpha)}\,
b(x, z)\, |z|^{\alpha-\beta} \right) .
\end{equation}
The next result gives a necessary and sufficient condition for the
kernel $q^b(t, x, y)$ in Theorem \ref{T:1.1} to be non-negative
when $b(x, z)$ is continuous in $x$ for a.e. $z$.

\medskip

\begin{thm}\label{T:1.2}  Let $b$ be a bounded   function on
$\R^d\times \R^d$  that satisfies \eqref{e:b} and that
\begin{equation}\label{e:5.1}
x\mapsto b(x, z) \hbox{ is continuous for a.e. } z\in \R^d.
\end{equation}
Then $q^b(t, x, y)\geq 0$ on $(0, \infty)\times \R^d \times \R^d$
if and only if for each $x\in \R^d$, $j^b(x, z)\geq 0$ for
a.e. $z\in \R^d$; that is,
if and only if
\begin{equation}\label{e:1.15}
b(x, z)\geq -\frac{{\cal A} (d, -\alpha)}{{\cal A} (d, -\beta)} \,
|z|^{\beta -\alpha} \quad
\hbox{for  a.e. } z\in \R^d.
\end{equation}
In particular, if $b(x, z)=b(x)$ is a function of $x$ only, then
$q^b(t, x, y)\geq 0$ on
 $(0, \infty)\times \R^d \times \R^d$
if and only if $b(x)\geq 0$ on $\R^d$.
\end{thm}

Next theorem drops the assumption \eqref{e:5.1},
gives lower bound estimates and refines upper bound estimates
 on $q^b (t, x, y)$
for $b(x, z)$ satisfying condition \eqref{e:1.15}
and makes connections to the martingale problem
for $\L^b$.
To state it, we need first to recall some definitions.
\bigskip

Let ${\mathbb D}([0, \infty), \R^d)$ be the space of right continuous
$\R^d$-valued functions having left limits on $[0, \infty)$, equipped
with Skorokhod topology. Denote by $X_t$ the
projection coordinate map on ${\mathbb D}([0, \infty), \R^d)$.
Let ${\cal C}$ be a subspace of $C^2_b(\R^d)$.
A probability measure $Q$ on the Skorokhod space
${\mathbb D}([0, \infty), \R^d)$ is said to
to be a solution to the martingale problem
for $(\L^b, {\cal C})$ with initial value $x\in \R^d$
if   $Q(X_0=x)=1$ and for every
$f\in {\cal C}$,
$$ M^f_t:=f(X_t) -f (X_0)-\int_0^t \L^b f(X_s)ds
$$
is a $Q$-martingale.
The martingale problem $(\L^b, {\cal C})$ with initial value $x\in \R^d$
is said to be well-posed if it has a unique solution.

Let $C_\infty (\R^d)$ be the space of continuous functions on $\R^d$
that vanish at infinity,  equipped with supremum norm. Set
$$
 C^2_\infty (\R^d)=\left\{ f\in C_\infty (\R^d): \,
\hbox{the first and second derivatives of } f \hbox{  are all in }
 C_\infty (\R^d) \right\}.
$$
A Markov process on $\R^d$ is called a Feller process if
its transition semigroup is a strongly continuous semigroup
in  $C_\infty (\R^d).$
Feller processes is a class of nice strong Markov processes,
called Hunt processes (see \cite{Chung}). Let $\overline{p}_0(t,x,y)$ be the fundamental solution of the truncated  operator
$$
\overline{\Delta}^{\alpha/2} f(x)
=\int_{|z|\leq 1} \left( f(x+z)-f(x)- \nabla f(x)
 \cdot  z\1_{\{|z|\leq 1\}}  \right) \frac{{\cal A}(d, -\alpha)}{|z|^{d+\alpha}}dz;
 $$
 or, equivalently, $\overline p_0 (t, x, y)$ is the transition
 density function for the finite range $\alpha$-stable (L\'evy) process
 with L\'evy measure ${\cal A}(d, -\alpha) |z|^{-(d+\alpha)}
 \1_{\{|z|\leq 1\}}$.
It is established in \cite{CKK} that
$\overline p_0(t, x, y)$ is jointly continuous and
enjoys the following two sided estimates:
\begin{equation}\label{e:trun1}
 \overline p_0(t, x, y) \asymp t^{-d/\alpha} \wedge
\frac{t}{|x-y|^{d+\alpha}}
\end{equation}
for $ t\in (0, 1]$ and  $|x-y|\leq 1$,
  and there are constants $c_k=c_k(d, \alpha)>0$, $k=1, 2, 3, 4$
  so that
\begin{equation}\label{e:trun2}
c_1 \left( \frac{t}{|x-y|}\right)^{c_2 |x-y|}
\leq \overline p_0(t, x, y) \leq
c_3 \left( \frac{t}{|x-y|}\right)^{c_4 |x-y|}
\end{equation}
for $ t\in (0, 1]$ and  $|x-y|> 1$.

\medskip
Define $b^+(x,z)=\max\{b(x,z), 0\}.$

\begin{thm}\label{T:1.3}
For every $A>0$ and $\lambda >0$,
there are positive constants
$C_k=C_k(d,\alpha,\beta, A)$, $k= 1,2$,
and  $C_3=C_3(d,\alpha,\beta, A, \lambda)$
such that
for any bounded $b$ satisfying \eqref{e:b} and \eqref{e:1.15}
with $\| b\|_\infty\leq A$,
\begin{equation}\label{e:1.14'}
C_1\overline{p}_0(t,C_2 x, C_2 y)\leq q^b(t,x,y)
\leq C_3 {p}_{M_{b^+, \lambda}}
(t,x,y)  \quad \hbox{for }
t\in (0,1] \hbox{ and }  x,y\in\R^d .
\end{equation}
Moreover, for every $\eps >0$, there is a positive constant
$C_4=C_4(d,\alpha,\beta, A, \lambda, \eps)$
such that for any   $b$ on $\R^d\times \R^d$
satisfying \eqref{e:b}  with $\| b\|_\infty\leq A$
so that
\begin{equation}\label{e:1.21a}
j^b(x, z)\geq  \eps\, |z|^{-(d+\alpha)} \quad
\hbox{for a.e. } x, z \in \R^d
\end{equation}
we have
\begin{equation}\label{e:1.14}
C_4  p_{m_{b^+, \lambda}}
(t,x,y)\leq q^b(t,x,y)
\leq C_3 {p}_{M_{b^+, \lambda}} (t,x,y)
 \quad \hbox{for } t\in (0,1] \hbox{ and }  x,y\in\R^d.
\end{equation}
The kernel $q^b(t, x, y)$ uniquely determines a Feller process $X^b
=(X^b_t, t\geq 0, \P_x, x\in \R^d)$
on the canonical Skorokhod space
${\mathbb D}([0, \infty), \R^d)$
such that
$$ \E_x \left[ f(X^b_t)\right] =\int_{\R^d} q^b (t, x, y) f(y) dy
$$
for every bounded continuous function $f$ on $\R^d$.
The Feller process $X^b$ is conservative
 and has a L\'evy system $(J^b(x, y)dy, t)$,
where $J^b(x, y)=j^b(x, y-x)$.
\begin{equation}\label{e:Jb}
 J^b(x, y)=j^b(x, y-x)= \frac{{\cal A}(d, -\alpha)}{|x-y|^{d+\alpha}}
+ \frac{{\cal A}(d, -\beta) \, b(x, y-x)}{|x-y|^{d+\beta}}.
\end{equation}
Moreover, for each $x\in \R^d$, $(X^b, \P_x)$ is the unique
solution to the martingale problem $(\L^b, {\cal S}(\R^d ))$
with initial value $x$. Here ${\cal S} (\R^d)$ denotes the
space of tempered functions on $\R^d$.
\end{thm}

\bigskip
Here we say $(J^b(x, y)dy, t)$ is a L\'evy system for $X^b$
if  for any non-negative
measurable function $f$ on $\R_+ \times \R^d\times \R^d$ with
$f(s, y, y)=0$ for all $y\in \R^d$, any stopping time $T$ (with
respect to the filtration of $X^b$) and any $x\in \R^d$,
\begin{equation}\label{e:levy}
\E_x \left[\sum_{s\le T} f(s,X^b_{s-}, X^b_s) \right]= \E_x \left[
\int_0^T \left( \int_{\R^d} f(s,X^b_s, y) J^b(X^b_s,y) dy \right)
ds \right].
\end{equation}
A  L\'evy system for $X^b$ describes the jumps of
the process $X^b$.
A Markov process on $\R^d$ is said to have strong Feller property if
its transition semigroup maps bounded measurable functions on $\R^d$
into bounded continuous functions on $\R^d$.
Since $q^b(t,x,y)$ is a  continuous function,
one has by  Theorem \ref{T:1.1} and the dominated convergence theorem
that  the Feller process $X^b$ of Theorem \ref{T:1.3}
has strong Feller property.

\medskip

Condition \eqref{e:1.21a} is always satisfied if $b(x, z)$ is nonnegative.
We emphasize the $m_{b^+, \lambda}$ and $M_{b^+, \lambda}$ terms
appeared in the estimates in Theorem \ref{T:1.3}.
Under condition \eqref{e:1.21a} and the assumption that $\| b \|_\infty
\leq A$,
the value of $b(x, z)$ on $\R^d \times \{ z\in \R^d: |z|\leq \lambda\}$
is irrelevant in the estimates of $q^b(t, x, y)$ in \eqref{e:1.14}.
By selecting suitable $\lambda>0$ in \eqref{e:1.14}, one can get optimal
two-sided estimates on $q^b(t, x, y)$.
The following follows immediately from Theorem \ref{T:1.3}
by taking a suitable $\lambda >0$.

\begin{cor}\label{C:1.4}
 Let $A\geq 0$  and $\eps>0$.
There is a positive constant
$C=C(d,\alpha,\beta, A,   \eps)\geq 1$
 so that for any bounded $b$ satisfying \eqref{e:b}
with $\| b\|_\infty\leq A$ and
$$
 j^b(x, z) \geq \eps \left( \frac{1}{|z|^{d+\alpha}}
 + \frac{1}{|z|^{d+\beta}} \right)
 \qquad \hbox{for a.e. } x, z\in \R^d,
$$
we have
$$ C^{-1} p_1 (t, x, y) \leq q^b(t, x, y) \leq C p_1(t, x, y)
\qquad \hbox{for } t\in (0, 1] \hbox{ and } x, y\in \R^d.
$$
\end{cor}

\medskip

Theorem \ref{T:1.3} in particular implies that
if  $b(x,\cdot)$ is a bounded function satisfying
\eqref{e:b} and \eqref{e:1.15}
so that $b(x, z)=0$ for every $x\in \R^d$ and $|z|\geq R$
for some $R>0$;
or, equivalently if  $\L^b=\Delta^{\alpha/2}+\S^b$ is a lower order perturbation of $\Delta^{\alpha/2}$
by finite range non-local operator $\S^b$,
then the upper bound of the kernel $q^b(t,x,y)$ is dominated by $p_0(t,x,y)$ for each
$(t,x,y)\in (0,1]\times\R^d\times\R^d.$
In fact, we have the following more general result.

\bigskip

\begin{thm}\label{T:1.4}
For every $A>0$ and $M \geq 1$,
 there is a constant $C_5=C_5(d,\alpha,\beta, A, M)\geq 1$
such that
for any bounded $b$
satisfying \eqref{e:b}
with $\| b\|_\infty \leq A$ and
\begin{equation}\label{e:1.21}
M^{-1}\,  |z|^{ -(d+\alpha)} \leq j^b(x,z)
\leq M\,  |z|^{ -(d+\alpha)}
 \quad \hbox{for a.e. }  x,z\in\R^d,
\end{equation}
or equivalently,
\begin{equation}\label{e:1.21'}
-(1-M^{-1})\frac{{\cal A}(d, -\alpha)}{{\cal A}(d, -\beta)}\,  |z|^{\beta -\alpha} \leq b(x, z)
\leq (M-1)\frac{{\cal A}(d, -\alpha)}{{\cal A}(d, -\beta)}
 |z|^{\beta -\alpha}
 \quad \hbox{for a.e.} \, x,z\in\R^d,
\end{equation}
we have
\begin{equation}\label{e:1.22}
C_5^{-1} p_0(t, x,  y)\leq q^b(t,x,y)
\leq C_5 {p}_0 (t,x,y)  \quad \hbox{for }
t\in (0,1] \hbox{ and }  x,y\in\R^d.
\end{equation}
\end{thm}

\medskip

We can restate some of results from Theorems \ref{T:1.1}, \ref{T:1.2},
\ref{T:1.3} and \ref{T:1.4} as follows.

\begin{thm}\label{T:1.5}
Let $b(x, z)$ be a bounded  function  on $\R^d\times \R^d$
satisfying \eqref{e:b} and \eqref{e:1.15}.
For each $x\in \R^d$,
 the martingale problem for $(\L^b, {\cal S}(\R^d ))$
with initial value $x$ is well-posed. These martingale
problem solutions $\{ \P_x, x\in \R^d\}$ form
a strong Markov process $X^b$,
which has infinite lifetime and possesses a jointly
continuous transition density function $q^b(t, x, y)$ with
respect to the Lebesgue measure on $\R^d$.
Moreover, the following holds.
\begin{description}
\item{\rm (i)} The transition density function $q^b(t, x, y)$
can be explicitly constructed as follows.
Define $q^b_0(t, x, y):=p_0(t, x, y)$
 and
$$
q^b_n (t, x, y):=\int_0^t \int_{\R^d} q^b_{n-1} (t-s, x, z)
\S^b_z p_0(s, z, y) dz ds \quad \hbox{for } n\geq 1.
$$
There is $\ee >0$ so that $\sum_{n=0}^\infty q^b_n (t, x, y)$
converges absolutely on $(0, \ee]\times \R^d\times \R^d$
and $q^b(t, x, y)= \sum_{n=0}^\infty q^b_n (t, x, y)$ on $(0, \ee]\times \R^d\times \R^d$.

\item{\rm (ii)}
$\displaystyle
q^b(t, x, y)=p_0(t, x, y)+\int_0^t \int_{\R^d} q^b (t-s, x, z)
\S^b_z p_0 (s, z, y)dz ds $ on $(0,\infty)\times\R^d\times\R^d.$

\item{\rm (iii)} For every $A>0$ and $\lambda >0$, there are positive constants
$c_k=c_k(d,\alpha,\beta,A), k=1,2,3$ and
$c_k=c_k(d,\alpha,\beta, A, \lambda)$, $k=4, \cdots, 9$,
 such that for any bounded function $b(x,z)$ on
  $\R^d \times \R^d$ satisfying \eqref{e:b} and \eqref{e:1.15}
with $\| b\|_\infty\leq A$,
$$
c_1 e^{-c_2t}\overline{p}_0(t,c_3x,c_3y)\leq q^b(t,x,y) \leq c_4 e^{c_5 t} \, {p}_{M_{b^+, \lambda}} (t,x,y) \quad \hbox{on } (0, \infty) \times \R^d
\times \R^d
$$
 and  for any non-negative function $b(x, z)$ on
  $\R^d \times \R^d$ satisfying \eqref{e:b}
with $\| b\|_\infty\leq A$,
$$
c_6 e^{-c_7 t}\,  p_{m_{b, \lambda}} (t,x,y)\leq q^b(t,x,y) \leq c_8 e^{c_9 t} \, {p}_{M_{b, \lambda}} (t,x,y) \quad \hbox{on } (0, \infty) \times \R^d
\times \R^d .
$$

\item{\rm (iv)} For every $A>0$  and $M\geq 1$, there are positive constants
$c_k=c_k(d,\alpha,\beta, A, M)$, $k= 10, \cdots, 13$,
 such that for any bounded function $b(x,z)$ on
  $\R^d \times \R^d$ satisfying \eqref{e:b} and \eqref{e:1.21}
with $\| b\|_\infty\leq A$,
$$
c_{10} e^{-c_{11}t} {p}_0(t,x,y)\leq q^b(t,x,y) \leq c_{12} e^{c_{13} t} \, {p}_0 (t,x,y) \quad \hbox{on } (0, \infty) \times \R^d
\times \R^d.
$$
\end{description}
\end{thm}

\begin{remark}\label{R:1.4} \rm
 (i) In general, we can not expect $q^b$ to have
  comparable lower and upper bound estimates.
  The estimates in \eqref{e:1.14'} and \eqref{e:1.14} are
 sharp in the sense that
$q^b(t, x, y)=p_0(t, x, y)$ when $b\equiv 0$,
$q^b(t, x, y)=p_1(t, x, y)$ when $b\equiv 1$,
and $q^b(t, x, y)=\overline p_0(t, x, y)$ when
$b(x, z)=0$ for $|z|\leq 1$
and $b(x, z)= - \frac{{\cal A}(d, -\alpha)}{ {\cal A}(d, -\beta)}\,
|z|^{\beta-\alpha}$ for $|z|\geq 1$.
 Clearly,
by \eqref{e:1.5}-\eqref{e:1.6}, $p_0(t, x, y)$ and
$p_1(t, x, y)$ are not comparable on $(0, 1]\times \R^d \times \R^d$.
We point out that it follows from \eqref{e:1.6} and \eqref{e:1.14}
that  every $A\geq 1 $, there is a constant
$\wt C = \wt C (d,\alpha,\beta, A)\geq 1 $
so  that for any non-negative $b$ on $\R^d\times \R^d$
satisfying \eqref{e:b} with $ 1/A \leq b(x, z)\leq A$ a.e.
\begin{equation}\label{e:1.16}
(1/ \wt C ) \,  p_1 (t,x,y)\leq q^b(t,x,y)
\leq \wt C \, {p}_1 (t,x,y)  \quad \hbox{for }
t\in (0,1] \hbox{ and }  x,y\in\R^d.
\end{equation}

\smallskip

(ii)  Heat kernel estimates for fractional Laplacian $\Delta^{\alpha/2}$
under gradient perturbation
and (possibly non-local) Feynman-Kac perturbation have recently been studied
in \cite{BJ, CKS, CKS2, W}. In both of these cases, under a Kato
class condition on the coefficients, the fundamental solution
of the perturbed operator is always strictly positive and is comparable
to the fundamental solution $p_0(t, x, y)$ of the fractional Laplacian
$\Delta^{\alpha/2}$ on $(0, 1]\times \R^d\times \R^d$. 

The novelty of this paper is on non-local perturbations.
The analysis of non-local perturbations with infinite jumping intensity measure
is much harder and is in fact very challenging. 
While the idea of using Duhamel's method \eqref{e:1.7n} in the study of
operator perturbation is not new, the key is how to implement it to obtain  two-sided
sharp heat kernel estimates where the lower bound is comparable to the upper bound,
and to establish the uniqueness of the fundamental solution.
It requires precise estimates on the non-local derivatives of the heat kernel
for fractional Laplacian, which turns out to be quite delicate and challenging.  
To the best of authors' knowledge, this is the first paper
on the study of heat kernels under non-local perturbation with infinite jump
intensity measure in a systematic way.  
We emphasize that the function $b(x, z)$ in \ref{e:1.1} is only measurable.  
Our Theorems \ref{T:1.2} and \ref{T:1.3} reveal  
some new phenomenon that heat kernels under non-local perturbation $\S^b$ are typically unstable.
This is  is  in stark contrast with $\Delta^{\alpha/2}$ under either gradient
(local) perturbations or (possibly non-local) Feynman-Kac perturbations.
However, Theorem \ref{T:1.4} of this paper in particular indicates
 that the heat kernel estimate for
 $\Delta^{\alpha /2}$
  is stable under finite range lower order perturbation.

 \smallskip

(iii) Kolokoltsov \cite{K} studied heat kernel estimates for symmetric pseudo-differential operators
(or stable-like jump diffusions) 
with smooth symbols. However neither the results nor the approach in \cite{K}
 applies to our case even when
$b(x, z)$ is assumed to be smooth. In addition to the  smooth  symbol requirement,
the operators  (1.8)-(1.9) considered in \cite{K} would require $\alpha=\beta$,
excluding the case where there are two different stable scales as are considered in this paper. 
In particular, it does not apply to SDE \eqref{e:sde}.  
 For  information on the connection between pseudo-differential operators and
 discontinuous Markov processes, we refer the reader to
 \cite{Ja, SW} and the references therein. 

\smallskip

(iv) Martingale problem for non-local operators (with or without
elliptic differential operator component) has been studied by
many authors. See, e.g.,
\cite{BC, BT, Ko1, Ko2, MP, NT, St, Ts} and the references
therein. 
In particular,
Komatsu \cite{Ko2} and Mikulevicious-Pragarauskas \cite{MP}
considered martingale problem for a class of non-local operators
that is directly related to $\L^b$. In fact, the uniqueness
of the martingale problem for $(\L^b,  {\cal S}(\R^d ))$
stated in Theorem \ref{T:1.3} above is a direct consequence
of  \cite[Theorem 3]{Ko2}, while it follows from \cite[Theorem 5]{MP}
that for any bounded $b$ satisfying  \eqref{e:b} and \eqref{e:1.15},
there is a unique solution to the martingale problem
$(\L^b, C^\infty_c(\R^d))$.  
  The main contribution of Theorem \ref{T:1.3} is on the two-sided
transition density function estimates for the martingale problem
solution $X^b_t$. We also mention that the well-posedness of martingale
problem for $(\Delta^{\alpha/2}+ b(x) \cdot \nabla, C^\infty_c(\R^d))$ with $b(x)$ an
$\R^d$-valued Kato class function has recently been established in \cite{CW}.

(v)  There are several directions to extend our results.
For example,  one can replay ${\cal A}(d, -\alpha)/|z|^{d+\alpha}$ in \eqref{e:1.2n}
and $1/|z|^{d+\beta}$ in \eqref{e:1.1} by the L\'evy kernel of pure jump subordinate
Brownian motion.  This is doable by following the ideas and approach of this paper. 
 Another direction is to consider Laplacian under non-local perturbation;
 that is, to replace $\Delta^{\alpha/2}$ in $\L^b$ by Laplacian operator $\Delta$.
 This has recently been carried out  in Wang \cite{Wa}. 
 
(v)  There are several directions to extend our results.
For example,  one can replay ${\cal A}(d, -\alpha)/|z|^{d+\alpha}$ in \eqref{e:1.2n}
and $1/|z|^{d+\beta}$ in \eqref{e:1.1} by the L\'evy kernels of pure jump subordinate
Brownian motions.  This should be doable by following the ideas and approach of this paper. 
 Another direction is to consider Laplacian under non-local perturbation;
 that is, to replace $\Delta^{\alpha/2}$ in $\LL^b$ by Laplacian operator $\Delta$.
 This has recently been carried out  in Wang \cite{Wa}.     \qed
\end{remark}

The rest of the paper is organized as follows. In Section
\ref{S:2}, we derive some estimates on $\overline{\Delta}^{\beta/2}_x p_0(t, x, y)$ and
 $\Delta^{\beta/2}_x p_0(t, x, y)$
that will be used in later.  The existence and uniqueness of the fundamental
solution $q^b(t, x, y)$ of $\L^b$ are given in Section \ref{S:3}.
This is done through a series of lemmas and theorems,
which provide more detailed information on $q^b(t, x, y)$
and $q^b_n(t, x, y)$.
Theorem \ref{T:1.1} then follows from these results.
We show in Section \ref{S:4}
that the semigroup $\{T^b_t; t>0\}$ associated with $q^b(t,x,y)$ is a strongly continuous
semigroup in $C_\infty (\R^d)$.
We then apply Hille-Yosida-Ray theorem and Courr\'ege's first theorem
to establish Theorem \ref{T:1.2}.
When $b$ satisfies \eqref{e:b}, \eqref{e:5.1} and \eqref{e:1.15}, $q^b(t, x, y)$ determines a conservative Feller process $X^b$.
We first derive a L\'evy system of $X^b$ and also prove $(X^b, \P_x)$ is the unique solution to the martingale problem for
$(\L^b, \S(\R^d))$  in Section \ref{S:5}.
We next establish,
for any given $A>0$,  the equi-continuity
of $q^b(t, x, y)$ on each $[1/M, M]\times \R^d \times \R^d$
for any $b$ that satisfies \eqref{e:b} with
$\| b\|_\infty \leq A$. Using this, we can drop the condition \eqref{e:5.1}
and establish the Feller process $X^b$ with transition density $q^b(t,x,y)$
for general bounded $b$ that satisfies \eqref{e:b} and
\eqref{e:1.15}
by approximating it with a sequence of $\{k_n(x, z), n\geq 1\}$
that satisfy \eqref{e:b}, \eqref{e:5.1} and \eqref{e:1.15}.
The upper bound estimate for $q^b(t,x,y)$ in \eqref{e:1.14'} and \eqref{e:1.14} can be obtained from that of $q^{\wh{b}_\lambda}(t,x,y)$
due to the Meyer's construction of $X^{\wh{b}_\lambda}$ from $X^b,$
 where
 $\wh{b}_\lambda(x,z)=b(x,z)1_{\{|z|\leq \lambda\}}(z)+b^+(x,z)1_{\{|z|>\lambda\}}(z).$
The lower bound estimates in \eqref{e:1.14'} and \eqref{e:1.14} are established by the L\'evy system of $X^b$ and some probability estimates.
Finally, we use the estimates in \eqref{e:1.14} for $b$ with support in $\{(x,z)\in\R^d\times\R^d: |z|\leq 1\}$
and the non-local Feynman-Kac perturbation results from
\cite{CKS2} to obtain Theorem \ref{T:1.4}.

Throughout this paper, we use the capital letters $C_1,C_2, \cdots $
to denote constants in the statement of the results, and their
labeling will be fixed. The lowercase constants $c_1, c_2, \cdots$
will denote generic constants used in the proofs, whose exact values
are not important and can change from one appearance to another.
 We will use ``$:=$" to denote a definition. For a differentiable function
 $f$ on $\R^d$, we use  $\partial_i f$ and $\partial^2_{ij}f$ to
 denote the partial derivatives $\frac{\partial f}{\partial x_i}$
 and $\frac{\partial^2 f}{\partial x_i \partial x_j}$.

\section{Preliminaries}\label{S:2}

Suppose that $Y$ is a symmetric $\alpha$-stable process, and $Z$ is a
symmetric $\beta$-stable process on $\R^d$ that is independent of $Z$.
  For any $a \ge 0$, we define $Y^a$ by $Y_t^a:=Y_t+ a^{1/\beta} Z_t$. We will call the process $Y^a$ the independent sum of the
symmetric $\alpha$-stable process  $Y$ and the symmetric
$\beta$-stable process  $Z$ with weight $a^{1/\beta}$. The infinitesimal
generator of $Y^a$ is $\Delta^{\alpha/2}+a \Delta^{\beta/2}$.
Let $p_a (t, x, y)$ denote the transition density of $Y^a$ (or
equivalently the heat kernel of $\Delta^{\alpha/2}+ a
\Delta^{\beta/2}$) with respect to the Lebesgue measure on $\R^d$.
 Recently it is proven in \cite{CK2} that
\begin{eqnarray}
p_1 (t, x, y)\asymp \left(  t^{-d/\alpha} \wedge t^{-d/\beta}
\right)\wedge \left( \frac{t}{|x-y|^{d+\alpha}} +
\frac{t}{|x-y|^{d+\beta}} \right) \quad \mbox{on } (0, \infty)\times
\R^d \times \R^d. \label{e:1.0}
\end{eqnarray}

Unlike the case of  the symmetric $\alpha$-stable process $Y:=Y^0$,
$Y^a$ does not have the stable  scaling for $a>0$. Instead, the
following approximate scaling property holds : for every
$\lambda>0$, $\{\lambda^{-1} Y^{a}_{\lambda^\alpha t}, t\geq 0\}$ has the same distribution as $\{Y^{a \lambda^{(\alpha-\beta)} }_t, t\geq 0\}$.
Consequently, for any $\lambda>0$, we have
\begin{equation}\label{e:scaling}
p_{a\lambda^{(\alpha-\beta) }}  ( t,  x, y) =
\lambda^d p_{a}  (\lambda^{\alpha}t, \lambda x, \lambda y) \qquad
\hbox{for } t>0 \hbox{ and } x, y \in \R^d.
\end{equation}
In particular, letting $a=1$, $\lambda= \gamma^{1/(\alpha -\beta)},$
 we get
$$
p_\gamma(t,x,y)= \gamma^{d/(\alpha-\beta)}p_1(\gamma^{\alpha/
 (\alpha-\beta)} t ,  \gamma^{1/(\alpha-\beta)} x ,
\gamma^{1/(\alpha-\beta)} y) \qquad \hbox{for } t>0
\hbox{ and } x, y\in \R^d.
$$
So we deduce from \eqref{e:1.0} that there exists
a constant $C>1$ depending only on
 $d$, $\alpha$ and $\beta$ such that for every $a>0$
 and $(t,x,y) \in (0,
\infty)\times \R^d \times \R^d$
\begin{equation}\label{e:h}
C^{-1} h_a(t, x, y) \leq p_a (t, x, y)\leq C h_a(t, x, y),
\end{equation}
where
$$
h_a(t, x, y):=  \left(    t^{-d/\alpha} \wedge
(a t )^{-d/\beta} \right) \wedge \left(\frac{t}{|x-y|^{d+\alpha}} +
\frac{a \, t}{|x-y|^{d+\beta}}\right).
$$
In fact, \eqref{e:h} also holds when $a=0$.
Observe  (see \eqref{e:1.7a}) that for every $A>0$, there is a constant
$c=c(d, \alpha, \beta, A)\geq 1$ so that  for every
$(t, x, y)\in (0, 1]\times \R^d\times \R^d$
 and  $0\leq a \leq A$,
\begin{equation}\label{e:ha2}
c^{-1}\,  t^{-d/\alpha} \wedge \left(\frac{t}{|x-y|^{d+\alpha}} +
\frac{a \, t}{|x-y|^{d+\beta}}\right)
\leq h_a(t, x, y) \leq c\, t^{-d/\alpha} \wedge \left(\frac{t}{|x-y|^{d+\alpha}} +
\frac{a \, t}{|x-y|^{d+\beta}}\right)
\end{equation}

\bigskip
Recall that
$p_0(t, x, y)=p_0(t, x-y)$
is the transition density function of
the symmetric $\alpha$-stable process $Y^0$.

\bigskip

\begin{lem}\label{0}
There exists a constant $C_6=C_6(d, \alpha)>0$
such that for every $t>0$, $x\in \R^d$ and
$i, j=1, \dots, d$,
$$
\left| \frac{\partial}{\partial x_i} p_0(t, x )\right| \leq C_6 t^{-(d+1)/\alpha} \left(1\wedge \frac{t^{1/\alpha}}{|x |}\right)^{d+1+\alpha},$$
$$
\left| \frac{\partial^2}{\partial x_i\partial x_j}   p_0 (t,x ) \right|\leq
C_6 t^{-(d+2)/\alpha} \left(1\wedge \frac{t^{1/\alpha}}{|x |}\right)^{d+2+\alpha}.
$$
\end{lem}

\pf
By  \cite[Lemma 5]{BJ}, there is a positive constant $c_1$ so that for all $t>0$ and $x,y\in\R^d$
$$
\left| \nabla_x p_0(t, x )\right| \leq c_1| x|  \left(t^{-(d+2)/\alpha}\wedge \frac{t}{|x |^{d+2+\alpha}}\right)
\leq c_1\left( t^{-(d+1)/\alpha}\wedge \frac{t}{|x |^{d+1+\alpha}}\right).
$$
That is, the first inequality holds.
Let $\eta_t(r)$ be the density function of the $\alpha/2$-stable subordinator at time $t$ and
$g(t,x)=(4\pi t)^{-d/2}e^{-|x|^2/4t}$  be the Gaussian kernel on $\R^d$.
There is a constant $c$ so that
 $\eta_t(r)\leq ctr^{-1-\alpha/2}$ for all $r,t>0,$ see \cite[Lemma 5]{BJ}. Note that
$$\left| \frac{\partial^2}{\partial x_i\partial x_j} g(s,x)\right|
\leq  \left( \frac{|x|^2}{s^2}+\frac{2}{s}\right) g(s,x)
=(4\pi)^2|x|^2 g^{(d+4)}(s,x_1)+8\pi g^{(d+2)}(s,x_2),
$$
where $x_1\in\R^{d+4}$ and $x_2\in\R^{d+2}$ with $|x_1|=|x_2|=|x|,$ $g^{(d+2)}(s,x_2)$ and $g^{(d+4)}(s,x_1)$ are the
Gaussian kernels on $\R^{d+2}$ and $\R^{d+4}$, respectively.
Since $p_0(t, x)=\int_0^\infty g(s, x) \eta_t (s) ds$, we have
 by the dominated convergence theorem that
 there is a positive constant $c_2$ so that for all $t>0$ and $x\in\R^d$
$$\begin{aligned}
\left|\frac{\partial^2}{\partial x_i\partial x_j}p_0(t,x)\right|&\leq \int_0^\infty \left| \frac{\partial^2}{\partial x_i\partial x_j} g(s,x)\right| \eta_t(s)\,ds\\
&\leq (4\pi)^2|x|^2 p^{(d+4)}_0(t,x_1)+8\pi p^{(d+2)}_0(t,x_2)\\
&\leq c_2 \left( t^{-(d+2)/\alpha}\wedge \frac{t}{|x|^{d+2+\alpha}}\right),
\end{aligned}$$
where  $p^{(d+2)}_0(t,x_2)$ and $p^{(d+4)}_0(t,x_1)$ are the
transition density functions of the symmetric $\alpha$-stable processes
in $\R^{d+2}$ and $\R^{d+4}$, respectively.
This establishes the second inequality in Lemma \ref{0}.
 \qed

\bigskip

 Define for $t>0$ and $x,y\in \R^d$,  the function
$$
|\Delta^{\beta/2}_x| p_0(t, x, y)
 \begin{cases}
={\cal A}(d, -\beta) \Big(\int_{|z|\leq t^{1/\alpha}}
\big| p_0 (t, x+z, y)-p_0(t,x, y)- \frac{\partial}{\partial x}p_0(t, x, y)
\cdot z  \big| \, \dfrac{1} {|z|^{d+\beta}}dz \\
 \hskip 0.3truein
  +\int_{|z|> t^{1/\alpha}}|p_0(t,x+z, y)-p_0(t,x, y)|\dfrac{dz}{|z|^{d+\beta}}\Big)
   \hskip 0.8truein \hbox{for } |x-y|^\alpha\leq t,\\
= {\cal A}(d, -\beta) \Big( \int_{|z|\leq |x-y|/2}|p_0(t,x+z, y)-p_0(t,x, y)- \frac{\partial}{\partial x}p_0(t,x, y)\cdot z |\dfrac{1}{|z|^{d+\beta}} dz\\
 \hskip 0.3truein  +\int_{|z|>|x-y|/2}|p_0(t,x+z, y)-p_0(t,x, y)|\dfrac{dz}{|z|^{d+\beta}}
\Big) \hskip 0.6truein \hbox{for } |x-y|^\alpha> t .
\end{cases}
$$
Let
\begin{equation}\label{e:3.9}
f_0(t, x, y):=\left( t^{1/\alpha} \vee |x-y| \right)^{-(d+\beta)}
= t^{-(d+\beta)/\alpha} \left( 1 \wedge \frac{t^{1/\alpha}}{|x-y|} \right)^{d+\beta}.
\end{equation}

\begin{lem}\label{1}
There exists a constant $C_7=C_7(d, \alpha, \beta )>0$ such that
\begin{equation}\label{e:3.2}
 |\Delta^{\beta/2}_x| p_0(t,x,y) \leq C_7 f_0(t, x, y)
\qquad \hbox{on } (0, \infty)\times \R^d \times \R^d.
\end{equation}
\end{lem}

\pf We only need to prove $ |\Delta^{\beta/2}_x| p_0(t,x) \leq C_7 f_0(t, x, 0)$ for all  $t>0$  and $x\in \R^d$.

 (i) We first consider the case $|x|^\alpha\leq t.$ In
this case,
$$\begin{aligned}
|\Delta^{\beta/2}_x| p_0(t,x)&= {\cal A}(d, -\beta)\int_{|z|\leq t^{1/\alpha}}|p_0(t,x+z)-p_0(t,x)- \frac{\partial}{\partial x}p_0(t,x)
\cdot z |\dfrac{dz}{|z|^{d+\beta}}\\
&\qquad +{\cal A}(d, -\beta)\int_{|z|\geq t^{1/\alpha}}|p_0(t,x+z)-p_0(t,x)|\dfrac{dz}{|z|^{d+\beta}} \\
&=I+II.
\end{aligned}$$

Note that by Lemma \ref{0},
$$\sup_{u\in \R^d} \left|\frac{\partial^2}{\partial u_i\partial u_j}p_0(t, u)\right|\leq C_6t^{-(d+2)/\alpha}, $$
and so  by Taylor's formula,
$$
I \leq {\cal A}(d, -\beta)\sup_{u\in \R^d} \left|\frac{\partial^2}{\partial u_i\partial u_j}p_0(t, u)\right|
\int_{|z|\leq t^{1/\alpha}} \frac{|z|^2}{|z|^{d+\beta}}\,dz
\leq c_1t^{-(d+2)/\alpha}t^{(2-\beta)/\alpha}\leq
c_1t^{-(d+\beta)/\alpha}.
$$
On the other hand, by \eqref{e:1.5}
$$
II\leq   {\cal A}(d, -\beta)\int_{|z|\geq t^{1/\alpha}}
\left(p_0(t,x+z)+p_0(t,x)\right) \frac{dz}{|z|^{d+\beta}}
 \leq c_2t^{-d/\alpha}\int_{|z|\geq t^{1/\alpha}}
\frac{1}{|z|^{d+\beta}}dz \leq c_3t^{-(d+\beta)/\alpha}.
$$

(ii) Next, we consider the case $|x|^\alpha\geq t.$ In this case,
$$\begin{aligned}
|\Delta_x^{\beta/2}| p_0(t,x)&= {\cal A}(d, -\beta)\int_{|z|\leq |x|/2}|p_0(t,x+z)-p_0(t,x)- \frac{\partial}{\partial x}p_0(t,x)\cdot z |\dfrac{dz}{|z|^{d+\beta}}\\
&\qquad +{\cal A}(d, -\beta)\int_{|z|\geq |x|/2}|p_0(t,x+z)-p_0(t,x)|\dfrac{dz}{|z|^{d+\beta}}\\
&=:I+II.
\end{aligned}$$
Note that $|x+z|\geq |x|/2$ for $|z|\leq |x|/2$. So by Lemma \ref{0},
$$\sup_{|z|\leq |x|/2}\Big|\frac{\partial^2}{\partial x_i\partial x_j}p_0(t,x+z) \Big|
\leq C_6\sup_{|z|\leq |x|/2}t|x+z|^{-(d+2+\alpha)}\leq 2^{(d+2+\alpha)}C_6t|x|^{-(d+2+\alpha)}.$$ Hence, by
Taylor's formula
\begin{equation}\label{e:w1}
\begin{aligned}
I&\leq {\cal A}(d, -\beta)\sup_{|z|\leq |x|/2} \Big|
\frac{\partial^2}{\partial x_i\partial x_j}p_0(t, x+z) \Big|
\int_{|z|\leq |x|/2} \frac{|z|^2}{|z|^{d+\beta}}\,dz\\
&\leq c_4t|x|^{-(d+2+\alpha)}|x|^{2-\beta}=
c_4t|x|^{-(d+\alpha+\beta)}.
\end{aligned}\end{equation}
Noting that $|x|^\alpha\geq t,$ thus $I\leq c_4|x|^{-(d+\beta)}.$
On the other hand, note that symmetric $\alpha$-stable process is a subordinate Brownian motion, so
$p_0(t,x+z)\leq p_0(t,x)$ if $|x+z|\geq |x|$ and $p_0(t,x)\leq p_0(t,x+z)$ if $|x+z|\leq |x|$.
Hence, by \eqref{e:1.5} and the condition that
$|x|^\alpha\geq t,$ we obtain
\begin{equation}\label{e:w2}\begin{aligned}
II&\leq {\cal A}(d, -\beta)\int_{|z|\geq |x|/2, |x+z|\geq |x|}
2p_0(t,x)\frac{dz}{|z|^{d+\beta}}
+{\cal A}(d, -\beta)\int_{|z|\geq |x|/2, |x+z|\leq |x|} 2p_0(t,x+z)\frac{dz}{|z|^{d+\beta}}\\
&\leq 2{\cal A}(d, -\beta)p_0(t,x)\int_{|z|\geq |x|/2}\frac{dz}{|z|^{d+\beta}}
+2^{d+1+\beta}{\cal A}(d, -\beta)|x|^{-(d+\beta)}\int_{z\in\R^d}p_0(t,x+z)\,dz\\
&\leq c_5t|x|^{-(d+\alpha)}|x|^{-\beta}+2^{d+1+\beta}{\cal A}(d, -\beta)|x|^{-(d+\beta)}\leq
c_6|x|^{-(d+\beta)}.
\end{aligned}\end{equation}
This establishes the lemma. \qed

\bigskip
In order to get the upper bound estimates in \eqref{e:1.7}
in terms of weight $M_{b, \lambda}$ rather than $\| b\|_\infty$,
we define,  for $t>0, \lambda>0$ and $x,y\in \R^d$,  the function
$$
|\Delta^{\beta/2}_{\lambda, x}| \, p_0(t, x, y)
 \begin{cases}
={\cal A}(d, -\beta) \Big(\int_{|z|\leq \lambda\wedge t^{1/\alpha}}
\big| p_0 (t, x+z, y)-p_0(t,x, y)- \frac{\partial}{\partial x}p_0(t, x, y)
\cdot z  \big| \, \dfrac{1} {|z|^{d+\beta}}dz \\
 \hskip 0.3truein
  +\int_{\lambda>|z|> (\lambda\wedge t^{1/\alpha})}|p_0(t,x+z,y)-p_0(t,x,y)|\dfrac{dz}{|z|^{d+\beta}}\Big)
   \hskip 0.6truein \hbox{for } |x-y|^\alpha\leq t,\\
= {\cal A}(d, -\beta) \Big( \int_{|z|\leq \lambda\wedge |x-y|/2}|p_0(t,x+z,y)-p_0(t,x,y)- \frac{\partial}{\partial x}p_0(t,x,y)\cdot z |\dfrac{1}{|z|^{d+\beta}} dz\\
 \hskip 0.3truein  +\int_{\lambda>|z|>(\lambda\wedge |x-y|/2)}|p_0(t,x+z,y)-p_0(t,x,y)|\dfrac{dz}{|z|^{d+\beta}}
\Big) \hskip 0.4truein \hbox{for } |x-y|^\alpha> t .
\end{cases}
$$
Observe that
$$
|\Delta^{\beta/2}_{\lambda, x}|\,  p_0(t, x, y)|
\leq |\Delta^{\beta/2}_x | \, p_0(t, x, y).
$$
Set
$$
f_{0, \lambda}(t,x,y)= \begin{cases}
t^{-(d+\beta)/\alpha} \qquad  &\hbox{when }  |x-y|\leq t^{1/\alpha},\\
|x-y|^{-(d+\beta)} \1_{\{|x-y|\leq \lambda\}}
+|x-y|^{-(d+\alpha)} \1_{\{|x-y|>\lambda\}}
&\hbox{when } |x-y|>t^{1/\alpha}.
\end{cases}
$$
Observe that when $\lambda =\infty$, $f_{0, \infty}$ is just the
function $f_0$ defined in \eqref{e:3.9}.

\begin{lem}\label{Lw1}
For each $\lambda>0$ and $T>0,$ there exists a constant $C_8=C_8(d, \alpha, \beta, \lambda, T)>0$ such that
\begin{equation}\label{e:3.2'}
 |\Delta^{\beta/2}_{\lambda, x}| \, p_0(t,x,y) \leq C_8 f_{0, \lambda}(t, x, y)
\qquad \hbox{on } (0, T]\times \R^d \times \R^d.
\end{equation}
\end{lem}

\pf  (i) We first consider the case $|x-y|^\alpha\leq t.$ Note that
$$\begin{aligned}
|\Delta^{\beta/2}_{\lambda, x}| \, p_0(t,x,y)&\leq |\Delta^{\beta/2}_x| p_0(t,x,y).
\end{aligned}$$
Hence, by the first part $(i)$ in the proof of Lemma \ref{1}, there exists a positive constant $c_1$ so that
$$|\Delta^{\beta/2}_{\lambda, x}| \, p_0(t,x,y)\leq c_1t^{-(d+\beta)/\alpha}.$$

(ii) Next, we consider the case $|x-y|^\alpha> t.$ In this case
$$\begin{aligned}
|\Delta^{\beta/2}_{\lambda, x}| \, p_0(t,x,y)&\leq {\cal A}(d, -\beta)\int_{|z|\leq |x-y|/2}|p_0(t,x+z,y)-p_0(t,x, y)- \frac{\partial}{\partial x}p_0(t,x, y)\cdot z |\dfrac{dz}{|z|^{d+\beta}}\\
&\qquad +{\cal A}(d, -\beta)\int_{\lambda\geq |z|\geq (\lambda\wedge |x-y|/2)}|p_0(t,x+z, y)-p_0(t,x, y)|\dfrac{dz}{|z|^{d+\beta}}\\
&=:I+II.
\end{aligned}$$
 By \eqref{e:w1}, there is a positive constant $c_2$  so that
$$
I\leq c_2t|x-y|^{-(d+\alpha+\beta)}
\leq c_3 \left( |x-y|^{-(d+\beta)} \1_{\{|x-y|\leq 2\lambda\}}+|x-y|^{-(d+\alpha)}\1_{\{|x-y|>2\lambda\}} \right).
$$
Here the last inequality holds since $t|x-y|^{-(d+\alpha+\beta)}\leq T(2\lambda)^{-\beta}|x-y|^{-(d+\alpha)}$ when $|x-y|>2\lambda$ and $t|x-y|^{-(d+\alpha+\beta)}\leq |x-y|^{-(d+\beta)}$ due to $|x-y|^\alpha\geq t.$

It is clear that $II=0$ if $|x-y|>2\lambda.$ On the other hand, if $|x-y|\leq 2\lambda,$ then there exists a positive constant $c_4$ so that
 $II\leq c_4 |x-y|^{-(d+\beta)}$ by \eqref{e:w2}.
 Finally, we note that $|x-y|^{-(d+\beta)}\asymp |x-y|^{-(d+\alpha)}$ for $\lambda<|x-y|\leq 2\lambda.$
This establishes the lemma. \qed

\bigskip

For each $\lambda>0$ and $a\geq 0$, we extend the definition
of $f_{0, \lambda}(t, x, y)$   to define
\begin{equation}\label{e:fa}
f_{a, \lambda}(t, x, y):= \begin{cases}
t^{-(d+\beta)/\alpha} \hskip 2.6truein\hbox{when }
 |x-y|\leq t^{1/\alpha},\\
|x-y|^{-(d+\beta)}\1_{\{|x-y|\leq \lambda\}}
+\left(|x-y|^{-(d+\alpha)}
+a \, |x-y|^{-(d+\beta)}\right) \1_{\{|x-y|>\lambda\}} \\
\hskip 3.2truein \hbox{when } |x-y|>t^{1/\alpha}.
\end{cases}
\end{equation}
Note that $f_{a, \infty}(t, x, y)=f_0 (t, x, y)$.

\begin{lem}\label{n1}
For each $\lambda>0,$ there is a constant $C_9=C_9(d, \alpha, \beta, \lambda)>0$
such that for every $a\in [0, 1]$,
\begin{equation}\label{L1}
\int_0^t\int_{\R^d}f_{a, \lambda}(s, z, y) dz ds  \leq C_9 \, (t^{1-\beta/\alpha}+t), \qquad  t\in(0, \infty),\, y\in\R^d.
\end{equation}
\end{lem}

\pf  By the definition of $f_{a, \lambda},$
 \begin{eqnarray*}
&&\int_0^t\int_{\R^d}f_{a, \lambda}(s, z, y)\,dz\,ds\\
&\leq& \int_0^t\int_{|y-z|\leq s^{1/\alpha}} s^{-(d+\beta)/\alpha}\,dz\,ds
+\int_0^t\int_{\lambda\geq|y-z|> s^{1/\alpha}} \frac{1}{|y-z|^{d+\beta}}\,dz\,ds \\
&& +\int_0^t\int_{|y-z|\geq \lambda}(|y-z|^{-(d+\alpha)}+ |y-z|^{-(d+\beta)})\,dz\,ds\\
&\leq & c_1 \int_0^t (s^{-\beta/\alpha}+1)\,ds \leq c_2(t^{1-\beta/\alpha}+t).
\end{eqnarray*}
\qed

For every $a\geq 0$,  define
\begin{equation}\label{ep}
g_a(t,x,y)= \begin{cases}
t^{-d/\alpha} \qquad  &\hbox{when }  |x-y|\leq t^{1/\alpha},\\
\frac{t}{|x-y|^{d+\alpha}}+\frac{a t}{|x-y|^{d+\beta}} &\hbox{when } |x-y|>t^{1/\alpha}.
\end{cases}
\end{equation}
Observe that
\begin{equation}\label{e:2.6}
\int_{\R^d} g_a(t, x, y) dy \asymp 1+ a t^{1-\beta/\alpha}
\qquad \hbox{on } (0, \infty) \times \R^d.
\end{equation}
Recall that $ {p}_a (t,x,y)$ is the heat kernel of the operator $\Delta^{\alpha/2}+ a \Delta^{\beta/2}.$
 Moreover,
in view of \eqref{e:1.7a},
\begin{equation}\label{e:3.3}
g_a (t, x, y)\asymp   p_a (t, x, y) \qquad \hbox{on }
 (0, 1]\times \R^d \times \R^d.
 \end{equation}

\medskip

\begin{lem}\label{2'}
For each $\lambda>0$ and $T>0,$ there exists $C_{10}=C_{10}(d,\alpha,\beta, \lambda, T)>0$ such that for every $a\in [0, 1]$ and all $t\in (0,T], x,y\in\R^d,$
$$\int_0^t \int_{\R^d}g_a(t-s,x,z)f_{a, \lambda}(s,z,y)\,dz\,ds
\leq C_{10} g_a(t,x,y).$$
\end{lem}

\pf  Denote by $I=\int_0^t\int_{\R^d} g_a (t-s,x,z)
f_{a, \lambda}(s,z,y)\,dz\,ds.$

(i) Suppose that $|x-y|\leq t^{1/\alpha}$. Then
 \begin{eqnarray*}
I &=& \int_0^t\int_{|x-z|\leq 2t^{1/\alpha}} g_a(t-s,x,z)
f_{a, \lambda}(s,z,y)\,dz\,ds\\
&&\quad+\int_0^t\int_{|x-z|> 2t^{1/\alpha}}
g_a(t-s,x,z) f_{a, \lambda}(s,z,y)\,dz\,ds\\
&&=:I_1+I_2.
\end{eqnarray*}

We write $I_1$ as
 \begin{eqnarray*}
 I_1&=&\int_0^{t/2}\int_{|x-z|\leq 2t^{1/\alpha}}
g_a(t-s,x,z)f_{a, \lambda}(s,z,y)\,dz\,ds\\
&&\quad+\int_{t/2}^{t}\int_{|x-z|\leq 2t^{1/\alpha}}
g_a(t-s,x,z) f_{a, \lambda}(s,z,y)\,dz\,ds\\
&=&I_{11}+I_{12}.
\end{eqnarray*}

If $s\in (0,t/2),$ then $t-s\in (t/2,t)$.
In this case,
$g_a(t-s,x,z)\leq c_1t^{-d/\alpha}$ when $|x-z|\leq 2t^{1/\alpha}$  by (\ref{ep}).
Hence, by Lemma \ref{n1},
$$
I_{11}\leq c_1t^{-d/\alpha}\int_0^t\int_{\R^d}f_{a, \lambda}(s,z,y)\,dz\,ds
\leq c_2(T^{1-\beta/\alpha}+T)  \, t^{-d/\alpha}.
$$

 When $s\in [t/2,t]$, since $|x-y|\leq t^{1/\alpha}$ and $|x-z|\leq 2t^{1/\alpha}$,
  $|y-z|\leq 3t^{1/\alpha}\leq 3 (2s)^{1/\alpha}$.  Thus $f_{a, \lambda}(s,z,y)\leq c_3s^{-(d+\beta)/\alpha}\leq 2^{(d+\beta)/\alpha}c_3t^{-(d+\beta)/\alpha}$.
Hence,
$$
I_{12}\leq 2^{(d+\beta)/\alpha}c_3t^{-(d+\beta)/\alpha}\int_0^t\int_{\R^d}g_a(t-s,x,z)\,dz\,ds
\leq c_4T^{1-\beta/\alpha}(1+T^{1-\beta/\alpha}) \, t^{-d/\alpha}.
$$

Next we consider $I_2$. Noting that
$|x-z|> 2t^{1/\alpha}$,
so  we have by \eqref{ep} and Lemma \ref{n1},
\begin{eqnarray*}
I_2
&\leq&  c_5\int_0^t\int_{|x-z|> 2t^{1/\alpha}}
\left(\frac{t-s}{|x-z|^{d+\alpha}}+\frac{t-s}{|x-z|^{d+\beta}}\right)
f_{a, \lambda}(s,z,y)\,dz\,ds\\
& \leq & c_6 t^{-d/\alpha} \left(1+t^{1-\beta /\alpha}\right)
\int_0^t\int_{\R^d} f_{a, \lambda}(s,z,y) \,dz\,ds\\
&\leq & c_7(1+T^{1-\beta/\alpha})(T^{1-\beta/\alpha}+T) \, t^{-d/\alpha} .
\end{eqnarray*}
We thus conclude from the above that there is a $c_8=c_8(d,\alpha,\beta, \lambda, T)>0$ such that $I\leq c_8 \, t^{-d/\alpha}$
for every $t\in (0,T]$ whenever $|x-y|\leq t^{1/\alpha}$.

(ii) Next assume that $|x-y|>t^{1/\alpha}$.  Then
 \begin{eqnarray*}
I &=& \int_0^t\int_{|x-z|\leq |x-y|/2} g_a(t-s,x,z)f_{a, \lambda}(s,z,y)\,dz\,ds\\
&&\quad+\int_0^t\int_{|x-z|>|x-y|/2}
g_a(t-s,x,z)f_{a, \lambda}(s,z,y)\,dz\,ds\\
&=:& I_1+I_2.
\end{eqnarray*}

If $|x-z|\leq |x-y|/2,$ then $|y-z|\geq |x-y|/2>t^{1/\alpha}/2$.
Hence, there is a constant $c_9$ so that
$$
f_{a, \lambda}(s,z,y)\leq
c_9 \left(|x-y|^{-(d+\alpha)}+a|x-y|^{-(d+\beta)}\right)
$$
for $s\in (0,t).$ Therefore,
 \begin{eqnarray*}
 I_1&\leq& c_9(|x-y|^{-(d+\alpha)}+a|x-y|^{-(d+\beta)})
\cdot
\int_0^t\int_{\R^d}
g_a(t-s,x,z)\,dz\,ds\\
&\leq& c_{10}(1+T^{1-\beta/\alpha})\left(\frac{t}{|x-y|^{d+\alpha}}+\frac{a t}{|x-y|^{d+\beta}}\right).
\end{eqnarray*}

If $|x-z|>|x-y|/2,$ then $|x-z|>t^{1/\alpha}/2$. Hence
$g_a(t-s,x,z)\leq
c_{11}\left( \frac{t}{|x-y|^{d+\alpha}}+\frac{a t}{|x-y|^{d+\beta}}\right)$
 by (\ref{ep}). Thus by Lemma \ref{n1},
 we obtain
 \begin{eqnarray*}
I_2
 &\leq& c_{11}\left( \frac{t}{|x-y|^{d+\alpha}}+\frac{a t}{|x-y|^{d+\beta}}\right)
 \int_0^t\int_{\R^d} f_{a, \lambda}(s,z,y)\,dz\,ds\\
&\leq& c_{12}(T^{1-\beta/\alpha}+T) \, \left( \frac{t}{|x-y|^{d+\alpha}}+\frac{a t}{|x-y|^{d+\beta}}\right).
\end{eqnarray*}
This completes the proof of the Lemma. \qed

\section{Fundamental solution}\label{S:3}

Throughout the  rest of this paper, $b(x, z)$ is a bounded function on
$\R^d\times \R^d$ satisfying condition \eqref{e:b}.
Recall the definition of the non-local operator $\S^b$ from \eqref{e:1.1}.
Let $|q^b|_0(t,x,y)=p_0(t,x,y),$
and define for each $n\geq 1$,
$$
|q^b|_n(t,x,y)=\int_0^t\int_{\R^d} |q^b|_{n-1}(t-s,x,z)|\S^b_z p_0
(s,z,y)|\,dzds.
$$
For each $\lambda>0,$ define
$$ b_\lambda(x,z)=b(x,z)1_{\{|z|>\lambda\}}(z).$$

In view of \eqref{e:1.5}, there exists a constant $C_{11}=C_{11}(d, \alpha, \beta)>0$
such that $p_0(t,x,y)\leq C_{11}g_{a}(t,x,y)$
for all $t>0, a\in[0,1]$ and $x,y\in\R^d$, where $g_a$ is the function defined by
\eqref{ep}.
On the other hand, note that
$$\begin{aligned}
|\S^b f(x)|
&=\left|{\cal A}(d, -\beta)
\int_{\R^d} \left( f(x+z)-f(x)-\langle\nabla f(x),
z \> \1_{\{|z|\leq \lambda\}} \right) \frac{b(x, z)}{|z|^{d+\beta}}dz\right|\\
&\leq \left|{\cal A}(d, -\beta)
\int_{|z|\leq \lambda} \left( f(x+z)-f(x)-\langle\nabla f(x),
z \>  \right) \frac{b(x, z)}{|z|^{d+\beta}}dz\right|\\
&\quad +\left|{\cal A}(d, -\beta)
\int_{\R^d} ( f(x+z)-f(x))\frac{b_\lambda(x, z)}{|z|^{d+\beta}}dz\right|\\
&\leq \| b\|_\infty \cdot |\Delta^{\beta/2}_{\lambda,x}| f(x)+\|b_\lambda\|_\infty \cdot |\Delta^{\beta/2}_x|f(x)\\
\end{aligned}$$
where $ |\Delta^{\beta/2}_{\lambda,x}| f(x)$ is defined in the similar way as  $ |\Delta^{\beta/2}_{\lambda,x}| p_0(t,x, y).$
Then by Lemma \ref{1}  and Lemma \ref{Lw1},
for every $A>0, \lambda>0$ and $T>0$ and every bounded function $b$ with $\|b\|_\infty\leq A, $
\begin{equation}\label{L}
\begin{aligned}
|\S^b_z p_0(t,z,y)|
&\leq \| b\|_\infty \cdot |\Delta^{\beta/2}_{\lambda,z}| p_0(t,z,y)+\|b_\lambda\|_\infty \cdot |\Delta^{\beta/2}_z| p_0(t,z,y)\\
&\leq C_8 A \, f_{0, \lambda}(t,z,y)+ C_7M_{b, \lambda} \, f_0(t,z,y)\\
&\leq (C_7+C_8) A f_{M_{b, \lambda}/A, \lambda}(t,z,y), \qquad t\in (0,T].
\end{aligned}
\end{equation}
Here  recall that $M_{b, \lambda} = {\rm esssup}_{x, z\in \R^d, |z|>\lambda} |b(x, z)|,$
$f_{a, \lambda}$ is the function defined in \eqref{e:fa}.
The above estimate is a refinement of Lemma \ref{1}. The latter
corresponds to the case where $\lambda =\infty$.

\begin{lem}\label{L:3.3}
For each $\lambda>0, A>0$ and $T>0$ and every bounded function $b$ on $\R^d\times\R^d$ satisfying condition \eqref{e:b}
with $\|b\|_\infty\leq A,$
\begin{equation}\label{k}
|q^b|_n(t,x,y)
\leq C_{11} \left(A(C_7+C_8) C_{10}  \right)^{n}g_{M_{b, \lambda}/A}(t,x,y)<\infty, \quad t\in (0,T],\, x,y\in\R^d.
\end{equation}
\end{lem}

\pf We prove this lemma by induction.
Since $p_0(t,x,y)\leq C_{11}g_{M_{b,\lambda}/A}(t,x,y)$ and $M_{b,\lambda}/A\leq 1,$
in view of
Lemma \ref{2'} and \eqref{L},   \eqref{k} clearly holds for $n=1$.
Suppose that \eqref{k} holds for $n=j\geq 1$.  Then by
 Lemma \ref{2'} and (\ref{L}),
\begin{eqnarray*}
&& |q^b|_{j+1}(t,x,y)\\
&\leq& C_{11} \left(A(C_7+C_8) C_{10} \right)^j
 \,
\int_0^t\int_{\R^d} g_{M_{b, \lambda}/A}(t-s,x,z)  |\S^{b}_z p_0(s,z,y)|\,dz\,ds\\
&\leq& C_{11} \left(A(C_7+C_8) C_{10} \right)^j
\, (C_7+C_8)A \,
\int_0^t\int_{\R^d} g_{M_{b, \lambda}/A}(t-s,x,z)
f_{M_{b,\lambda}/A,
\lambda} (s,z,y)\,dz\,ds\\
&\leq& C_{11} \left(A(C_7+C_8) C_{10} \right)^{j+1}g_{M_{b, \lambda}/A}(t,x,y)
\end{eqnarray*}
for $t\in (0, T]$ and $x, y\in \R^d$.
This proves that \eqref{k} holds for $n=j+1$ and thus for every $n\geq 1$. \qed

\bigskip

Now we define $q_n^b: (0,\infty)\times\R^d\times\R^d\rightarrow \R$ as follows.
 For $t>0$ and $x,y\in\R^d,$
let $q^b_0(t,x,y)=p_0(t,x,y),$ and for each $n\geq 1$, define
\begin{equation}\label{e:3.6}
q_n^b(t,x,y)=\int_0^t\int_{\R^d} q^b_{n-1}(t-s,x,z)\S^b_z p_0(s,z,y)\,dz\,ds.
\end{equation}
Clearly by Lemma \ref{L:3.3},
each $q_n^b(t,x,y)$ is well defined.

\begin{lem}\label{con}
For every $n\geq 0$, $q_n^b(t,x,y)$ is jointly continuous on  $(0,\infty)\times\R^d\times\R^d.$
\end{lem}

\pf  We prove it by induction.
Clearly $q_0^b(t,x,y)$ is continuous on $(0,\infty)\times\R^d\times\R^d$.
Suppose that $q_n^b(t, x, y)$ is continuous on
$(0,\infty)\times\R^d\times\R^d$.
For every $M\geq 2$ ,
it follows from \eqref{L} , Lemma \ref{L:3.3} and the dominated convergence theorem that for $\ee<1/(2M),$
$$ (t, x, y) \mapsto \int_\ee^{t-\ee}\int_{\R^d} q_n^b(t-s,x,z)\S^b_z p_0(s,z,y) dz ds
$$
is jointly continuous on $[1/M, M]\times \R^d\times \R^d$.
On the other hand, it follows from  \eqref{L} and \eqref{e:2.6}
that
 \begin{eqnarray*}
&&\sup_{t\in [1/M, M]} \sup_{x,y}\int_{t-\ee}^t\int_{\R^d} g_{M_{b, \lambda}}(t-s,x,z)
|\S^b_z p_0(s,z,y)|\,dz\,ds\\
&\leq &c_1 A \left(\sup_{t\in [1/M, M]}[(t-\ee)^{-(d+\beta)/\alpha}+(t-\ee)^{-(d+\alpha)/\alpha}]\right)
\sup_{x\in \R^d}\sup_{t\in [1/M, M]}\int_{t-\ee}^t\int_{\R^d} g_{M_{b, \lambda}}(t-s,x,z) \,dz\,ds \\
&\leq & c_2 A (2M)^{(d+\alpha)/\alpha}  \int_0^\ee  (1+r^{1-\beta/\alpha}) dr \\
&\leq&  c_3 A  (2M)^{(d+\alpha)/\alpha} \ee,
\end{eqnarray*}
which goes to zero as $\ee\rightarrow 0$; while
by \eqref{L} and \eqref{L1}, there exist $c_4$ and $c_5$ such that
\begin{eqnarray}\label{h1}
&& \sup_{t\in [1/M, M]} \sup_{x,y}\int_0^\ee\int_{\R^d} g_{M_{b, \lambda}}(t-s,x,z)|\S^b_z p_0 (s,z,y)|\,dz\,ds \nonumber \\
&\leq &  c_4 \left(\sup_{t\in [1/M, M]} (t-\ee)^{-d/\alpha}\right)
\sup_{y\in \R^d} \int_0^\ee\int_{\R^d} |\S^b_z p_0 (s,z,y)|\,dz\,ds \nonumber \\
&\leq &  c_5 (2M)^{d/\alpha} \, \| b\|_\infty \,
\ee^{1-\beta/\alpha}\rightarrow 0
\end{eqnarray}
as $\ee\rightarrow 0.$
We conclude from Lemma \ref{L:3.3} and the above argument that
$$ q^b_{n+1}(t,x,y)=\int_0^t\int_{\R^d} q^b_{n}(t-s,x,z)
\S^b_z p_0(s,z,y)\,dz\,ds
$$
is jointly continuous in $(t, x, y)\in [1/M, M]\times\R^d\times\R^d$
and so in $(t, x, y)\in (0, \infty)\times\R^d\times\R^d$.
This completes the proof of the lemma. \qed

\bigskip

Recall $f_0(t, x, y)$ is the function defined in \eqref{e:3.9} and
$$ |\Delta^{\beta/2}_x| p_0(t,x,y) \leq C_7 f_0(t, x, y)
\qquad \hbox{on } (0, \infty)\times \R^d \times \R^d.$$

\bigskip

\begin{lem} \label{L:3.5}
  There is a constant $C_{12} =C_{12} (d, \alpha, \beta)>0$
so that for every $A>0$ and every bounded function $b$ on $\R^d\times\R^d$ with $\|b\|_\infty\leq A$ and for every integer $n\geq 0$ and $\eps >0$,
\begin{equation}\label{e:3.6b}
 \left|\int_{\{z\in \R^d: |z|>\eps\}}
\left( q^b_n (t, x+z, y) -q^b_n(t, x, y)\right)
\frac{{\cal A}(d, -\beta) b(x, z)}{|z|^{d+\beta}} dz \right|
\leq  (C_{12}A)^{n+1} f_0(t, x, y)
\end{equation}
for $(t, x, z) \in (0, 1]\times \R^d\times \R^d$,
 and $\S^b_x q^b_n (t, x, y)$ exists pointwise for $(t, x, z)
\in (0, 1]\times \R^d\times \R^d$ in the sense of
\eqref{e:1.5a} with
\begin{equation}\label{e:qbn1}
\S^b_x q^b_{n+1}(t, x, y) = \int_0^t \int_{\R^d}
\S^b_x q^b_n (t-s, x, z) \S^b_z p_0(s, z, y)  dz ds
\end{equation}
and
\begin{equation}\label{e:3.11}
|\S^b_x q^b_n (t, x, y)|\leq  (C_{12}A)^{n+1} f_0(t, x, y)
\qquad \hbox{on } (0, 1]\times \R^d\times \R^d.
\end{equation}
Moreover,
\begin{equation}\label{e:qbn2}
q^b_{n+1}(t, x, y)=\int_0^t\int_{\R^d} p_0(t-s, x, z)\S^b_z q^b_{n}(s, z, y)
  dz ds
 \qquad \hbox{for } (t, x, y)\in (0, 1]\times \R^d\times \R^d.
\end{equation}
\end{lem}

\pf Let $q(t, x, y)$ denote the transition density function
  of the symmetric $\beta$-stable process on $\R^d$.
Then by \eqref{e:1.5} but with $\beta$ in place of $\alpha$,
we have
\begin{equation}\label{e:3.10}
q(t, x, y)\asymp t^{-d/\beta} \left( 1 \wedge \frac{t^{1/\beta}}{|x-y|} \right)^{d+\beta}\qquad \hbox{on } (0, \infty) \times \R^d \times \R^d.
\end{equation}
Observe that \eqref{e:3.9} and \eqref{e:3.10} yield
\begin{equation}
 f_0(t, x, y) \asymp  t^{-\beta/\alpha} q(t^{\beta/\alpha}, x, y)
 \qquad \hbox{on } (0, \infty)\times \R^d \times \R^d.
\end{equation}
Hence on $(0, \infty)\times \R^d\times \R^d$,
\begin{eqnarray*}
&& \int_0^t  \int_{\R^d} f_0(t-s, x, z)   f_0(s, z, y)  ds dz\\
&\asymp&  \int_0^t  (t-s)^{-\beta/\alpha} s^{-\beta/\alpha}
\left(\int_{\R^d} q((t-s)^{\beta/\alpha}, x, z)
  q( s^{\beta/\alpha}, z, y)  dz \right) ds \\
&=& \int_0^t  (t-s)^{-\beta/\alpha} s^{-\beta/\alpha}
q( (t-s)^{\beta/\alpha}+s^{\beta/\alpha}, x, y) ds \\
&\asymp& q( t^{\beta/\alpha}, x, y) \int_0^t  (t-s)^{-\beta/\alpha}
   s^{-\beta/\alpha} ds  \\
&=& q( t^{\beta/\alpha}, x, y)\, t^{1-(2\beta/\alpha)} \int_0^1 (1-u)^{-\beta/\alpha} u^{-\beta/\alpha} du \\
&\asymp &  t^{1-\beta/\alpha} f_0(t, x, y).
\end{eqnarray*}
In the second $\asymp$ above, we used the fact that
$$
(t/2)^{\beta/\alpha}\leq (t-s)^{\beta/\alpha}+s^{\beta/\alpha}
\leq 2t^{\beta/\alpha} \qquad \hbox{for every } s\in (0, t)
$$
and the estimate \eqref{e:3.10}, while in the last equality, we
used a change of variable $s=tu$. So there is a constant
$c_1=c_1(d, \alpha, \beta)>0$ so that
\begin{equation}\label{e:3.13}
\int_0^t  \int_{\R^d} f_0(t-s, x, z)   f_0(s, z, y)  ds dz
\leq c_1 \, f_0(t, x, y) \qquad \hbox{for every } t\in (0, 1]
\hbox{ and } x, y\in \R^d.
\end{equation}
By increasing the value of $c_1$ if necessary, we may and do assume
that $c_1$ is larger than $1$.

We now proceed by induction. Let $C_{12}:=c_1C_7.$
Note that
\begin{equation}\label{e:f0'}
|\S^b_xp_0(t,x,y)|\leq A|\Delta^{\beta/2}_x| p_0(t,x,y) \leq C_7A f_0(t, x, y).
\end{equation}
When $n=0$, \eqref{e:qbn2} holds by definition.
By Lemma \ref{1},
\eqref{e:3.6b} and \eqref{e:3.11} hold for $n=0$.
Suppose that  \eqref{e:3.6b} and \eqref{e:3.11} hold for $n=j$.
Then for every $\eps >0$,
 by the definition of $q^b_{j+1}$, Lemma \ref{L:3.3}, \eqref{e:3.13} and Fubini's theorem,
\begin{eqnarray}\label{e:3.11a}
&&   \int_{\{w\in \R^d: |\omega|>\eps\}}
\left( q^b_{j+1} (t, x+w, y) -q^b_{j+1} (t, x, y)\right)
\frac{{\cal A}(d, -\beta) b(x, w)} {|w|^{d+\beta}} dw  \\
&=& \int_0^t\int_{\R^d} \left(
 \int_{\{w\in \R^d: |w|>\eps\}}
\left( q^b_j (t-s, x+w, z) -q^b_j (t-s, x, z)\right)
 \frac{{\cal A}(d, -\beta) b(x, w)} {|w|^{d+\beta}} dw\right)
  \nonumber \\
&&  \hskip 0.3truein \times \, \S^b_z p_0(s, z, y)\, dz ds \nonumber
\end{eqnarray}
and so
\begin{eqnarray*}
&&  \left| \int_{\{w\in \R^d: |w|>\eps\}}
\left( q^b_{j+1} (t, x+w, y) -q^b_{j+1} (t, x, y)\right)
\frac{{\cal A}(d, -\beta) b(x, w)} {|w|^{d+\beta}} dw \right| \nonumber \\
&\leq& \int_0^t\int_{\R^d} (C_{12}A)^{j+1} f_0(t-s, x, z)\,
  |\S^b_z p_0(s, z, y)| dz ds  \\
&\leq & \int_0^t\int_{\R^d} (C_{12}A)^{j+1} f_0(t-s, x, z) \, C_7A  f_0(s, z, y)   dz ds \\
&\leq & (C_{12}A)^{j+2} f_0(t, x, y).
 \end{eqnarray*}
By \eqref{e:3.11a} and Lebesgue dominated convergence theorem,
we conclude that
\begin{eqnarray*}
&& S^b_x q^b_{j+1}(t, x, y) \\
&:=& \lim_{\eps \to 0}
 \int_{\{w\in \R^d: |w|>\eps\}}
\left( q^b_{j+1} (t, x+w, y) -q^b_{j+1} (t, x, y)\right)
\frac{{\cal A}(d, -\beta) b(x, w)} {|w|^{d+\beta}} dw \\
&=&\int_0^t\int_{\R^d} \left(
 \lim_{\eps \to 0} \int_{\{w\in \R^d: |w|>\eps\}}
\left( q^b_j (t-s, x+w, z) -q^b_j (t-s, x, z)\right)
 \frac{{\cal A}(d, -\beta) b(x, w)} {|w|^{d+\beta}} dw\right)
  \\
&& \hskip 0.3truein \times \, \S^b_z p_0(s, z, y)\, dz ds \\
&=& \int_0^t\int_{\R^d} \S^b_x q^b_j (t-s, x, z)\,
     \S^b_z p_0(s, z, y)\, dz ds
\end{eqnarray*}
exists and \eqref{e:qbn1} as well as \eqref{e:3.11} holds for $n=j+1$.
(The same proof verifies \eqref{e:qbn1} when $n=0$.)
On the other hand, in view of \eqref{e:3.11} and \eqref{e:qbn2}
for $n=j$, we have by the Fubini theorem,
\begin{eqnarray*}
&&  q^b_{j+1} (t, x, y) \\
&=& \int_0^t \int_{\R^d}
q^b_j (s, x, z) \S^b_z p_0(t-s, z, y)  dz ds\\
&=& \int_0^t \int_{\R^d}
\left( \int_0^{s} \int_{\R^d} p_0(r, x, w) \S^b_w q^b_{j-1}
(s-r, w, z) dr dw \right) \S^b_z p_0(t-s, z, y)  dz ds \\
&=& \int_0^t \int_{\R^d} p_0(r, x, w)
\left( \int_r^t \int_{\R^d} \S^b_w q^b_{j-1} (s-r, w, z)
\S^b_z p_0(t-s, z, y) ds dz \right)  dw dr \\
&=& \int_0^t \int_{\R^d} p_0(r, x, w)
\S^b_w q^b_j(t-r, w, y)   dw dr . \\
\end{eqnarray*}
This verifies that \eqref{e:qbn2} also holds for $n=j+1$.
The lemma is now established by induction.
\qed

\bigskip
Recall that $M_{b, \lambda} = {\rm esssup}_{x, z\in \R^d, |z|>\lambda} |b(x, z)| =\|b_\lambda(x, z)\|_\infty.$

\begin{lem}\label{u1}
For each $\lambda>0,$ there are  positive constants
$A_0=A_0(d, \alpha, \beta, \lambda)$
and $C_{13}=C_{13}(d,\alpha,\beta, \lambda)$ so that if $\| b\|_\infty
\leq A_0$, then for every integer $n\geq 0$,
\begin{equation}\label{e:3.8}
|q_{n+1}^b(t,x,y)|
\leq C_{13} 2^{-n} \,  p_{M_{b, \lambda}}(t,x,y)
\quad \hbox{for } t\in(0,1] \hbox{ and }  x,y\in\R^d,
\end{equation}
\eqref{e:3.6b} holds and so
$\S^b_x q_n^b(t,x,y)$ exists pointwise in the sense of \eqref{e:1.5a}
 with
\begin{equation}\label{e:3.15}
  |\S^b_x q_n^b(t,x,y)|
  \leq   2^{-n} \, f_0(t,x,y) \quad \hbox{for } t\in(0,1] \hbox{ and }  x,y\in\R^d,
\end{equation}
and
\begin{equation}\label{e:3.17}
  \sum_{n=0}^\infty q_n^b(t,x,y)\geq \frac{1}{2}p_0(t,x,y) \quad \hbox{for } t\in(0,1]
\hbox{ and } |x-y|\leq 3 t^{1/\alpha} .
\end{equation}
\end{lem}

\pf  We take a positive constant $A_0$ so that $A_0\leq 1\wedge[2(C_7+C_8)C_{10}+2C_{12}]^{-1}.$
  We have by Lemma \ref{L:3.3} and Lemma \ref{L:3.5} that
for every $b$ with $\| b\|_\infty \leq A_0$,
$$
|q_{n+1}^b(t,x,y)|\leq
C_{11} 2^{-n} g_{M_{b, \lambda}/A_0}(t, x, y)\leq C_{11}A_0^{-1}2^{-n}g_{M_{b, \lambda}}(t, x, y)
\quad \hbox{ and } \quad |\S^b_x q^b_n(t, x, y)|\leq 2^{-n} f_0(t, x, y)
$$
for every $t\in (0,1]$ and  $x,y\in\R^d$.
This together with \eqref{e:3.3} establishes \eqref{e:3.8}
 and \eqref{e:3.15}.

On the other hand, by  \eqref{ep}, there exists $c=c(d, \alpha, \beta)\geq 1$  so that $g_{a}(t,x,y)\leq cp_0(t,x,y)$
for $a\in [0,1]$ and $|x-y|\leq 3t^{1/\alpha}$ and $t\in (0,1]$.
Take $A_0$ small enough so that $A_0\leq 1\wedge[2(C_7+C_8)C_{10}+2C_{12}]^{-1} $ and $\sum_{n=1}^\infty (A_0(C_7+C_8)C_{10})^n\leq \frac{1}{2cC_{11}}.$
 Then for every $b$ with $\|b\|_\infty \leq A_0$,
 we have by Lemma \ref{L:3.3} for $|x-y|\leq 3t^{1/\alpha}$ and $t\in (0,1]$
 that
$$
\sum_{n=1}^\infty |q^b|_n(t,x,y)\leq    cC_{11}\sum_{n=1}^\infty (A_0(C_7+C_8)C_{10})^n
 p_0(t,x,y)\leq \frac{1}{2} p_0(t,x,y).
$$
Consequently, for $|x-y|\leq 3t^{1/\alpha}$ and $t\in (0,1]$,
$$
\sum_{n=0}^\infty q_n^b(t,x,y)\geq p_0(t,x,y)-\sum_{n=1}^\infty |q_n^b(t,x,y)|\geq \frac{1}{2}p_0(t,x,y).
$$
\qed

\medskip

We now extend the results in Lemma \ref{u1}
to any bounded $b$ that satisfies condition \eqref{e:b}.
For $\lambda >0$, define
\begin{equation}\label{e:4.1}
b^{(\lambda)} (x, z)= \lambda^{\beta/\alpha -1}
b(\lambda^{-1/\alpha} x, \lambda^{-1/\alpha} z).
\end{equation}
For a function $f$ on $\R^d$, set
$$ f^{(\lambda)} (x):= f(\lambda^{-1/\alpha} x).
$$
By a change of variable, one has from \eqref{e:1.2n} and \eqref{e:1.1}
that
$$ \Delta^{\alpha/2} f^{(\lambda)} (x)= \lambda^{-1} (\Delta^{\alpha/2}  f) (\lambda^{-1/\alpha}x)
$$
and
\begin{equation}\label{e:4.2}
 \S^{b^{(\lambda)}} f^{(\lambda)} (x)= \lambda^{-1} ( \S^b f) (\lambda^{-1/\alpha}x).
\end{equation}
We remark here that condition \eqref{e:b} used in establishing \eqref{e:4.2}.
Note that the transition density function $p_0(t, x, y)$ of
the symmetric $\alpha$-stable process has the following scaling property:
\begin{equation}\label{e:4.3}
p_0(t, x, y)= \lambda^{-d/\alpha} p_0(\lambda^{-1}t, \lambda^{-1/\alpha} x, \lambda^{-1/\alpha} y)
\end{equation}
Recall  $q^b_n(t, x, y)$ is the function defined inductively by \eqref{e:3.6} with $q^b_0(t, x, y):=p_0(t, x, y)$.

\begin{lem}\label{L:3.5b}
Suppose that $b$ is a bounded function on $\R^d\times \R^d$ satisfying \eqref{e:b}. For every $\lambda >0$ and for every integer $n\geq 0$
and $t>0,$
\begin{equation}\label{e:4.4}
 q^{b^{(\lambda)}}_n (t, x, y)=
\lambda^{-d/\alpha} \, q^{b}_n  (\lambda^{-1}t, \lambda^{-1/\alpha} x, \lambda^{-1/\alpha} y), \qquad x, y\in \R^d;
\end{equation}
or, equivalently,
\begin{equation}\label{e:4.5}
q^b_n (t, x, y)=
\lambda^{d/\alpha} \, q^{b^{(\lambda)}}_n  (\lambda t, \lambda^{1/\alpha} x, \lambda^{1/\alpha} y), \qquad x, y\in \R^d.
\end{equation}
\end{lem}

\pf We prove it by induction. Clearly in view of \eqref{e:4.3},
\eqref{e:4.4} holds when $n=0$. Suppose that \eqref{e:4.4} holds
for $n=j\geq 0$. Then by the definition \eqref{e:3.6}, \eqref{e:4.2}
and \eqref{e:4.3},
\begin{eqnarray*}
&& q^{b^{(\lambda)}}_{j+1} (t, x, y) = \int_0^t \int_{\R^d}
q^{b^{(\lambda)}}_j (t-s, x, z) \S^{b^{(\lambda)}}_z p_0(s, z , y)\, dz ds \\
&=&  \int_0^t \int_{\R^d} \lambda^{- d/\alpha}
 q^{b}_j  (\lambda^{-1}(t-s), \lambda^{-1/\alpha} x, \lambda^{-1/\alpha} z)
\lambda^{-d/\alpha-1} \left( \S^b_z p_0 (\lambda^{-1}s, \cdot , \lambda^{-1/\alpha}y)\right) (\lambda^{-1/\alpha}z) \, dz ds \\
&=& \lambda^{-d/\alpha} \int_0^{\lambda^{-1}t} \int_{\R^d}
 q^{b}_j  (\lambda^{-1}t-r, \lambda^{-1/\alpha} x, w)
 \left( \S^b_w p_0 (r, \cdot , \lambda^{-1/\alpha} y)\right)(w) \,  dw dr\\
&=& \lambda^{-d/\alpha}  q^{b}_{j+1} (\lambda^{-1}t, \lambda^{-1/\alpha} x, \lambda^{-1/\alpha} y).
\end{eqnarray*}
This proves that \eqref{e:4.4} holds for $n=j+1$ and so, by induction,
it holds for every $n\geq 0$. \qed

\bigskip

Recall that $A_0$ is the positive constant in Lemma \ref{u1}.

\medskip

\begin{thm}\label{T:3.7}
 For every $\lambda>0$ and $A>0$, there is a positive constant $C_{14}=C_{14}(d, \alpha, \beta, A, \lambda)>0$
 so that for every bounded function $b$ with $\| b\|_\infty \leq A$,
  that satisfies condition
\eqref{e:b} and $n\geq 0$,
\begin{equation}\label{e:3.16n}
  |q_n^b(t,x,y)|
  \leq C_{14} 2^{-n} \, \left( t^{-d/\alpha}\wedge
  \left( \frac{t}{|x-y|^{d+\alpha}} + \frac{ M_{b, \lambda} \, t}{|x-y|^{d+\beta}} \right)\right)
\end{equation}
for every $0<t \leq 1\wedge ( A_0/\|b\|_\infty)^{\alpha/(\alpha -\beta)}$
and $x,y\in\R^d$, and
\begin{equation}\label{e:3.17n}
  \sum_{n=0}^\infty q_n^b(t,x,y)\geq \frac{1}{2}p_0(t,x,y) \quad \hbox{for }
  0<t \leq 1\wedge( A_0/\|b\|_\infty)^{\alpha/(\alpha -\beta)}
\hbox{ and } |x-y|\leq 3 t^{1/\alpha} .
\end{equation}
Moreover, for every $n\geq 0$, \eqref{e:3.6b} holds and so
$\S^b_x q_n^b(t,x,y)$ exists
pointwise in the sense of \eqref{e:1.5a} with
\begin{equation}\label{e:3.24}
  |\S^b_x q_n^b(t,x,y)|
  \leq  2^{-n}(\|b\|_\infty/A_0)  f_0(t, x, y)
\end{equation}
for every $0<t \leq 1\wedge ( A_0/\|b\|_\infty)^{\alpha/(\alpha -\beta)}$
and $x,y\in\R^d$. Moreover, \eqref{e:qbn1} and \eqref{e:qbn2} hold.
\end{thm}

\pf In view of Lemma \ref{u1}, it suffices to prove the theorem
 for $A_0<\| b\|_\infty \leq A$.
Set $r=( \|b\|_\infty/A_0)^{\alpha/(\alpha -\beta)}$.
The function $b^{(r)}$ defined by \eqref{e:4.1} has
the property $\|b^{(r)} \|_\infty =A_0$.
Thus by Lemma \ref{u1},
  there is a constant
$C_{14}=C_{14}(d, \alpha, \beta, A, \lambda):=C_{13} (d, \alpha, \beta,
 r^{1/\alpha} \lambda)>0$ so that
  for every integer $n\geq 0$,
\begin{equation}\label{e:3.8c}
|q_{n}^{b^{(r)}}(t,x,y)|
\leq C_{14} \, 2^{-n} \,  p_{M_{b^{(r)}, r^{1/\alpha} \lambda }}(t,x,y)
\quad \hbox{for } t\in(0,1] \hbox{ and }  x,y\in\R^d,
\end{equation}
\eqref{e:3.6b} holds and so
$\S^{b^{(r)}}_x q_n^{b^{(r)}}(t,x,y)$ exists pointwise in the sense of \eqref{e:1.5a}
 with
\begin{equation}\label{e:3.15c}
  |\S^{b^{(r)}}_x q_n^{b^{(r)}}(t,x,y)|
  \leq   2^{-n} \, f_0(t,x,y) \quad \hbox{for } t\in(0,1] \hbox{ and }  x,y\in\R^d,
\end{equation}
and
\begin{equation}\label{e:3.17c}
  \sum_{n=0}^\infty q_n^{b^{(r)}}(t,x,y)\geq \frac{1}{2}p_0(t,x,y) \quad \hbox{for } t\in(0,1]
\hbox{ and } |x-y|\leq 3 t^{1/\alpha} .
\end{equation}
Noting  $r^{1-\beta/\alpha}
  M_{b^{(r)},  r^{1/\alpha}\lambda}=M_{b, \lambda}$, we have
by   \eqref{e:4.5}, \eqref{e:3.8c} and \eqref{e:h} that
for every $0<t\leq 1/r= (A_0/\|b\|_\infty)^{\alpha/(\alpha -\beta)}$
and $x, y\in \R^d$,
\begin{eqnarray*}
|q_n^b(t,x,y)|&=&
r^{d/\alpha} \,  |q^{b^{(r)}}_n  (r t, r^{1/\alpha} x, r^{1/\alpha} y) |  \\
&\leq& C_{14} 2^{-n} r^{d/\alpha} \, p_{M_{b^{(r)}, r^{1/\alpha} \lambda}}
 (r t, r^{1/\alpha} x, r^{1/\alpha} y)  \\
 &\leq & 2C \, C_{14}\, 2^{-n} \left( t^{-d/\alpha}\wedge
  \left( \frac{t}{|x-y|^{d+\alpha}} + \frac{ r^{1-\beta/\alpha}
  M_{b^{(r)}, r^{1/\alpha}\lambda} \, t}{|x-y|^{d+\beta}} \right)\right)\\
&\leq& 2C \, C_{14}\, 2^{-n} \left( t^{-d/\alpha}\wedge
  \left( \frac{t}{|x-y|^{d+\alpha}} + \frac{
  M_{b,   \lambda} \, t}{|x-y|^{d+\beta}} \right)\right),
\end{eqnarray*}
which establishes \eqref{e:3.16n}.
Similarly,  \eqref{e:3.17n}  follows
from \eqref{e:4.3}, and \eqref{e:3.17c}, while
the conclusion of \eqref{e:3.24} is a direct consequence of
\eqref{e:4.2}, \eqref{e:4.5} and \eqref{e:3.15c}.
That \eqref{e:qbn1} and \eqref{e:qbn2} hold follows directly
from Lemma \ref{L:3.5} and Lemma \ref{L:3.5b}.
\qed

\bigskip

Recall that  $q^b(t,x,y):=\sum_{n=0}^\infty q_n^b(t,x,y)$, whenever
it is convergent.
The following theorem follows immediately from Lemmas \ref{con},
 \ref{u1} and Theorem \ref{T:3.7}.

\begin{thm}\label{c1}
For every $\lambda>0$ and $A>0$,
let $C_{14}=C_{14}(d, \alpha, \beta, A , \lambda)$
 be the constant in Theorem
\ref{T:3.7}. Then for every bounded function
$b$ with $\| b\|_\infty \leq A$ that satisfies condition
\eqref{e:b}, $q^b(t,x,y)$ is well defined and is
jointly continuous in $(0, 1\wedge ( A_0/\|b\|_\infty)^{\alpha/(\alpha -\beta)} ]\times\R^d\times\R^d$.
Moreover,
$$
  |q^b(t, x, y)| \leq 2C_{14}\,  \left( t^{-d/\alpha}\wedge
  \left( \frac{t}{|x-y|^{d+\alpha}} +  \frac{  M_{b, \lambda} \, t}{|x-y|^{d+\beta}} \right)\right)
$$
and $\S^b_x q^b(t, x, y)$ exists pointwise in the sense of
\eqref{e:1.5a} with
$$ | \S^b_x q^b(t, x, y) | \leq  2(\|b\|_\infty/A_0) f_0(t, x, y)
$$
  for every $0<t \leq 1\wedge ( A_0/\|b\|_\infty)^{\alpha/(\alpha -\beta)}$
and $x,y\in\R^d$, and
\begin{equation}\label{e:ql1}
q^b(t,x,y)\geq \frac{1}{2}p_0(t,x,y) \quad \hbox{for }
  0<t \leq 1\wedge ( A_0/\|b\|_\infty)^{\alpha/(\alpha -\beta)}
\hbox{ and } |x-y|\leq 3 t^{1/\alpha} .
\end{equation}
Moreover, for every
$0<t \leq 1\wedge ( A_0/\|b\|_\infty)^{\alpha/(\alpha -\beta)}$ and $x,y\in\R^d$,
\begin{eqnarray}\label{e:3.21}
q^b(t, x, y)&=&p_0(t, x, y)+\int_0^t \int_{\R^d} q^b(t-s, x, z)
\S^b_z p_0(s, z, y)  dz ds\\
&=& p_0(t, x, y)+\int_0^t \int_{\R^d} p_0 (t-s, x, z)
\S^b_z q^b (s, z, y)dz ds \label{e:Du2}.
\end{eqnarray}
\end{thm}

\bigskip

\begin{thm}\label{T:3.9}
Suppose that $b$ is a bounded function on $\R^d\times \R^d$ satisfying \eqref{e:b}. Let $A_0$ be the constant in Lemma \ref{u1}.
 Then for every  $t, s>0$
with $t+s\leq 1\wedge ( A_0/\|b\|_\infty)^{\alpha/(\alpha -\beta)}$ and $x, y\in \R^d$,
\begin{equation}\label{e:3.26}
\int_{\R^d} q^b(t, x, z) q^b(s, z, y) dz = q^b(t+s, x, y).
\end{equation}
\end{thm}

\pf In view of Theorem \ref{T:3.7}, we have
$$ \int_{\R^d} q^b(t, x, z) q^b(s, z, y) dz
=\sum_{j=0}^\infty \sum_{n=0}^j \int_{\R^d} q^b_n(t, x, z)
q^b_{j-n}(s, z, y) dz.
$$
So it suffices to show that for every $j\geq 0$,
\begin{equation}\label{e:3.12}
\sum_{n=0}^j \int_{\R^d} q^b_n(t, x, z)
q^b_{j-n}(s, z, y) dz=q^b_j (t+s, x, y)
\end{equation}
Clearly, \eqref{e:3.12} holds for $j=0$.
Suppose that \eqref{e:3.12} holds for $j=l\geq 1$.
Then we have by Fubini's theorem and the estimates
in \eqref{L} and Theorem \ref{T:3.7},
\begin{eqnarray*}
&& \sum_{n=0}^{l+1} \int_{\R^d} q^b_n(t, x, z)
   q^b_{l+1-n}(s, z, y) dz \\
&=&  \int_{\R^d} q^b_{l+1}(t, x, z)
   p_0 (s, z, y) dz + \sum_{n=0}^{l} \int_{\R^d} q^b_n(t, x, z)
   q^b_{l+1-n}(s, z, y) dz  \\
&=& \int_{\R^d} \left(   \int_0^t \int_{\R^d}
     q^b_{l}(t-r, x, w) \S^b_w p_0(r, w, z)  dw dr \right)
     p_0 (s, z, y) dz \\
&& + \sum_{n=0}^{l} \int_{\R^d} q^b_n(t, x, z)
   \left( \int_0^s\int_{\R^d} q^b_{l-n}(s-r, z, w)\S^b_w
   p_0(r, w, y)  dw dr\right)  dz \\
&=& \int_0^t \int_{\R^d}  q^b_{l}(t-r, x, w) \S^b_w
    p_0(r+s, w, y)  dw dr\\
&& +  \int_0^s  \int_{\R^d} q^b_l(t+s-r, x, w) \S^b_w
   p_0(r, w, y) dw dr   \\
&=& q^b_{l+1} (t+s, x, y).
\end{eqnarray*}
This proves that \eqref{e:3.12} holds for $j=l+1$.
So by induction, we conclude that \eqref{e:3.12} holds for
every $j\geq 0$.  \qed

\bigskip

For notational simplicity,
denote $1\wedge ( A_0/\|b\|_\infty)^{\alpha/(\alpha -\beta)}$
 by $\delta_0$.
 In view of Theorem \ref{T:3.9}, we can uniquely extend the definition of
$q^b(t, x, y)$ to $t> \delta_0$
by using the Chapman-Kolmogorov equation recursively as follows.

Suppose that $q^b(t, x, y)$ has been defined and satisfies the
Chapman-Kolmogorov equation \eqref{e:3.26} on $(0, k\delta_0]\times \R^d
\times \R^d$. Then for $t\in (k\delta_0, (k+1)\delta_0]$, define
\begin{equation}\label{pk}
q^b(t,x,y) =\int_{\R^d}q^b(s,x,z)q^b(r,z,y)\,dz, \quad  x, y\in \R^d
\end{equation}
for any $s, r\in (0, k\delta_0]$ so that $s+r=t$.
Such $q^b(t, x, y)$ is well defined on $(0, \infty)\times \R^d\times \R^d$
and satisfies \eqref{e:3.26} for every $s, t>0$.
Moreover, since Chapman-Kolmogorov equation holds for $q^b(t, x, y)$ for all $t, s>0$,
we have by Theorem \ref{c1} and \eqref{e:h}-\eqref{e:ha2}
that for every $A\geq A_0$,
there are constants
$c_i=c_i(d, \alpha, \beta, A, \lambda)$, $i=1, 2$, so that for every
$b(x, z)$ satisfying \eqref{e:b} with $\| b\|_\infty \leq A$,
\begin{equation}\label{e:3.26n}
|q^b(t, x, y)| \leq c_1 \, e^{c_2t} \,   p_{M_{b, \lambda}} (t, x, y)
\qquad \hbox{for every } t>0 \hbox{ and } x, y\in \R^d .
\end{equation}

\begin{thm}\label{T:3.10a}
$q^b(t, x, y)$ satisfies \eqref{e:3.21} and \eqref{e:Du2} for every $t>0$ and $x,y\in\R^d.$
\end{thm}

\pf Let $\delta_0:=1\wedge ( A_0/\|b\|_\infty)^{\alpha/(\alpha -\beta)}.$
It suffices to prove that for every $n\geq 1,$
\eqref{e:3.21} and \eqref{e:Du2}  hold for all $t\in (0,n\delta_0]$ and $x,y\in\R^d.$

Clearly, \eqref{e:3.21} holds for $t\in (0,n\delta_0]$ with $n=1.$
Suppose that  \eqref{e:3.21} holds for $t\in (0,n\delta_0]$ with $n=k.$
For $t\in (k\delta_0, (k+1)\delta_0]$, take $l,s\in (0, k\delta_0]$ so that $l+s=t.$
Then we have by Fubini's theorem, Chapman-Kolmogorov equation of $q^b$, Lemma \ref{2'}, \eqref{L} and (\ref{e:3.26n}),
 \begin{eqnarray*}
 q^b(l+s, x,y)
&=& \int_{\R^d} q^b(l,x,z)q^b(s,z,y)\,dz\\
&=& \int_{\R^d} q^b(l,x,z) \left(p_0(s,z,y)+\int_0^s\int_{\R^d} q^b(s-r, z, \omega) \S^b_\omega p_0(r, \omega, y)\,d\omega\,dr\right)\,dz\\
&=& \int_{\R^d} p_0(l,x,z) p_0(s,z,y)\,dz\\
&& +\int_{\R^d} \left(\int_0^l\int_{\R^d} q^b(l-u, x, \eta) \S^b_\omega p_0(u, \eta, z)\,d\eta\,du\right) p_0(s,z,y)\,dz\\
&& +\int_0^s\int_{\R^d} q^b(l+s-r, x, \omega) \S^b_\omega p_0(r, \omega, y)\,d\omega\,dr\\
&=& p_0(l+s,x,y)+\int_0^l\int_{\R^d} q^b(l-u, x, \eta) \S^b_\omega p_0(u+s, \eta, y)\,d\eta\,du\\
&& +\int_0^s\int_{\R^d} q^b(l+s-r, x, \omega) \S^b_\omega p_0(r, \omega, y)\,d\omega\,dr\\
&=& p_0(l+s,x,y)+\int_0^{l+s} \int_{\R^d} q^b(l+s-r, x, z)\S^b_z p_0(r,z,y)\,dz\,dr.
\end{eqnarray*}
By the similar procedure as above, we can also prove that \eqref{e:Du2} holds for every $t>0$ and $x,y\in\R^d.$
\qed
\bigskip

\begin{thm}\label{T:3.10}
Suppose that $b$ is a bounded function on $\R^d\times \R^d$ satisfying \eqref{e:b}.
Then $q^b(t, x, y)$ is the unique continuous kernel that satisfies
the Chapman-Kolmogorov equation \eqref{e:3.26} on $(0, \infty)\times \R^d
\times \R^d$ and that
 for some $\ee>0$,
\begin{equation}\label{e:3.22}
|q^b(t, x, y)| \leq c \, p_1 (t, x, y)
\end{equation}
 and \eqref{e:3.21} hold
 for $(t, x, y)\in (0, \ee]\times \R^d \times \R^d$.
 Moreover, \eqref{e:3.26n} holds for $q^b(t, x, y)$.
\end{thm}

\pf Suppose that $\overline  q$ is any continuous kernel that satisfies,
for some $\ee>0$,  \eqref{e:3.21} and \eqref{e:3.22} hold for
$(t, x, y)\in (0, \ee]\times \R^d \times \R^d$.
Without loss of generality, we may and do assume that
$\ee <1\wedge (A_0/\| b\|_\infty)^{\alpha/(\alpha -\beta)}$.
Using \eqref{e:3.21} recursively, one gets
\begin{equation} \label{e:3.29}
\overline  q (t, x, y)=\sum_{j=0}^n q^b_j (t, x, y)
 + \int_0^t \int_{\R^d}  \overline  q (t-s, x, z)
 (\S^b p_0)^{*, n+1}_z (s,z,y)
ds dz.
 \end{equation}
Here $(\S^b p_0)^{*, n}_z (s,z,y) $  denotes the $n$th
convolution operation of the function $\S^b_z p_0 (s, z, y)$;
that is, $(\S^b p_0)^{*, 1}_z (s,z,y)= \S^b_z p_0 (s,z,y)$
and
\begin{equation}\label{e:s}
(\S^b p_0)^{*, n}_z (s,z,y)= \int_0^s \int_{\R^d}
\S^b_z p_0(r,z,w) \, (\S^b p_0)^{*, n-1}_w (s-r,w,y) dw dr
\quad \hbox{for } n\geq 2.
\end{equation}
It follows from \eqref{e:qbn1} that
$$
\S^b_zq^b_n(s,z,y)=(\S^b p_0)^{*, n+1}_z (s,z,y).
$$
Thus, by \eqref{e:3.29} we have
\begin{equation}\label{e:f1}
\overline  q (t, x, y)=\sum_{j=0}^n q^b_j (t, x, y)
 + \int_0^t \int_{\R^d}  \overline  q (t-s, x, z)
 \S^b_z q^b_n(s,z,y)
ds dz.
\end{equation}
By the condition
\eqref{e:3.22} and \eqref{e:3.24}, there is a constant $c_1>0$ so that for every $n\geq 1$,
$$ \left| \int_0^t \int_{\R^d}  \overline q (t-s, x, z)
(\S^b p_0)^{*, n}_z (s,z,y) ds dz \right| \leq c_1 2^{-n} \int_0^t \int_{\R^d}  p_1 (t-s, x, z)
f_0(s, z, y) ds dz.
$$
Noting that $p_1(t,x,y)\asymp g_1(t,x,y)$ on $(0,1]\times\R^d\times\R^d$ and
\begin{eqnarray}
 \int_0^t\int_{\R^d}f_0(s, z, y)\,dz\,ds
&\leq& \int_0^t\int_{|y-z|\leq s^{1/\alpha}} s^{-(d+\beta)/\alpha}\,dz\,ds
+\int_0^t\int_{|y-z|> s^{1/\alpha}} \frac{1}{|y-z|^{d+\beta}}\,dz\,ds \
\nonumber \\
&=&  c_2t^{1-\beta/\alpha}. \label{e:f0}
\end{eqnarray}
Then by the similar proof in Lemma \ref{2'}, we can get
$$\int_0^t \int_{\R^d}  p_1 (t-s, x, z)
f_0(s, z, y) ds dz \leq c_3p_1(t,x,y).$$
It follows that
$$  \overline q (t, x, y)=\sum_{n=0}^\infty q_n^b (t,x, y)
=q^b(t, x, y)
$$
for every $t\in (0, \ee]$ and $x, y\in \R^d$.
Since both $ \overline q$ and $q^b$ satisfy the Chapman-Kolmogorov
equation \eqref{e:3.26},
$\overline q=q^b$ on $(0, \infty) \times \R^d \times \R^d$.
\qed

\bigskip

\begin{remark}\rm
It follows from the definition of $q^b_n(t, x, y)$ and
Lemma \ref{L:3.5} that
$ (\S^b p_0)^{*, n+1} (s,z,y) = \S^b_z q^b_n(s, z, y)$.
\qed
\end{remark}

\bigskip

\bigskip

In view of Lemma \ref{L:3.5b} and Chapman-Kolmogorov equation, we have

\begin{thm}\label{T:3.11}
Suppose that $b$ is a bounded function on $\R^d\times \R^d$ satisfying \eqref{e:b}.
 \ $q^b(t, x, y)=\lambda^{d/\alpha} q^{b^{(\lambda)}}(\lambda t, \lambda^{1/\alpha} x, \lambda^{1/\alpha} y)$ on $(0, \infty)\times
\R^d\times \R^d$, where $b^{(\lambda)}(x, z):=\lambda^{\beta/\alpha-1}
 b(\lambda^{-1/\alpha} x, \lambda^{-1/\alpha} z)$.
\end{thm}

\bigskip

For a bounded function $f$ on $\R^d$, $t>0$ and $x\in\R^d$,
 we define
$$
T^b_tf(x)=\int_{\R^d} q^b(t,x,y)f(y)\,dy
\quad \hbox{and} \quad
P_t f(x)=\int_{\R^d} p_0(t, x, y) f(y) dy.
$$
The following lemma follows immediately from \eqref{e:3.26}
and \eqref{pk}.

\medskip

\begin{lem}\label{lq^b}
Suppose that $b$ is a bounded function on $\R^d\times \R^d$ satisfying \eqref{e:b}.
For all $s,t>0,$ we have $T^b_{t+s}=T^b_tT^b_s.$
\end{lem}

\bigskip

\begin{thm}\label{T:3.13}
Let $b$ be a bounded function on $\R^d\times \R^d$ satisfying \eqref{e:b}.
Then for every $f\in C^2_b (\R^d)$,
$$
T^b_t f(x) -f(x)=\int_0^t T^b_s \L^b f(x) ds
\qquad \hbox{for every } t>0, \, x\in \R^d .
$$
\end{thm}

\pf  Note that by Theorem \ref{T:3.10a}, for each bounded Borel function $f$ in $\R^d,$
\begin{equation}\label{e:3.32}
T^b_t f(x)=P_t f(x)+\int_0^t T^b_{t-s} \S^b P_s f(x) ds
=P_t f(x)+\int_0^t T^b_s \S^b P_{t-s} f(x) ds.
\end{equation}
Hence,  for $f\in C^2_b (\R^d)$,
\begin{eqnarray*}
&& T^b_t f(x)-f(x) \\
&=&P_t f(x)-f(x) +\int_0^t T^b_s \S^b f(x) ds +
\int_0^t T^b_s \S^b (P_{t-s} f-f) (x) ds \\
&=& \int_0^t P_s \Delta^{\alpha/2} f(x) ds
+ \int_0^t T^b_s \S^b f(x) ds +
\int_0^t T^b_s \S^b (P_{t-s} f-f) (x) ds      \\
&=& \int_0^t T^b_s \Delta^{\alpha/2} f(x) ds
 - \int_0^t \left(\int_0^{s} T^b_r \S^b P_{s-r}(\Delta^{\alpha/2}f)(x) dr \right) ds \\
&& + \int_0^t T^b_s \S^b f(x) ds +
\int_0^t T^b_s \S^b (P_{t-s} f-f) (x) ds  \\
&=& \int_0^t T^b_s \left( \Delta^{\alpha/2}+\S^b\right) f(x) ds
 - \int_0^t \left(\int_r^{t} T^b_r \S^b P_{s-r}(\Delta^{\alpha/2}f)(x) ds
  \right) dr \\
&& +
\int_0^t T^b_s \S^b (P_{t-s} f-f) (x) ds  \\
&=& \int_0^t T^b_s  \L^b  f(x) ds
 - \int_0^t   T^b_r \S^b ( P_{t-r}f -f)(x)   dr  +
\int_0^t T^b_s \S^b (P_{t-s} f-f) (x) ds  \\
&=& \int_0^t T^b_s  \L^b  f(x) ds.
\end{eqnarray*}
Here in the third inequality, we used \eqref{e:3.32};
while in the fifth inequality we used Lemma \ref{1}
and \eqref{e:3.26n}, which
allow the interchange of the integral sign $\int_r^t$
with $T^b_r \S^b$, and
the fact that
$$  \int_r^t  P_{s-r}(\Delta^{\alpha/2}f )(x)ds
= \int_r^t \left(\frac{d}{ds} P_{s-r}f (x)\right) ds
= P_{t-r}f(x)-f(x).
$$
\qed

\bigskip

\begin{thm}\label{lambda1}
Let $b$ be a bounded function on $\R^d\times \R^d$ satisfying \eqref{e:b}.
 Then $q^b(t,x,y)$ is jointly continuous in $(0,\infty)\times\R^d\times\R^d$  and $\int_{\R^d}q^b(t,x,y)\,dy=1$ for
 every $x\in \R^d$ and  $t>0.$
\end{thm}

\pf  By Lemma \ref{lq^b}, we have
\begin{equation}\label{eq^b}
q^b(t+s,x,y)=\int_{\R^d}q^b(t,x,z)q^b(s,z,y)\,dz, \quad x,y\in\R^d, s,t>0.
\end{equation}
Continuity of $q^b(t,x,y)$ in $(t, x, y)\in (0, \infty) \times \R^d
\times \R^d$ follows from Theorem \ref{c1}, (\ref{eq^b}) and the dominated convergence theorem.
For $n\geq 1$ and $t\in (0, T]$,
it follows from \eqref{L}, Lemma \ref{2'}, Theorem \ref{T:3.7} and  Fubini's Theorem
that for every $t\in (0, 1\wedge( A_0/\|b\|_\infty)^{\alpha/(\alpha -\beta)}]$,
$$\begin{aligned}
\int_{\R^d}q_n^b(t,x,y)\,dy
&=\int_{\R^d}\int_{\R^d}\int_0^t q^b_{n-1}(t-s,x,z) \S_z^b p_0(s,z,y)\,ds\,dz\,dy\\
&=\int_{\R^d}\int_0^t q^b_{n-1}(t-s,x,z)\S^b_z \left(\int_{\R^d} p_0(s,z,y)\,dy\right)\,ds\,dz =0.\\
\end{aligned}$$
Hence we have by Lemma \ref{u1},
$$\int_{\R^d}q^b(t,x,y)\,dy=\int_{\R^d} p_0(t,x,y)\,dy=1$$ for
$t\in (0,1\wedge( A_0/\|b\|_\infty)^{\alpha/(\alpha -\beta)}]$. This conservativeness property extends to all $t>0$ by (\ref{eq^b}). \qed

\bigskip

Theorem \ref{T:1.1} now follows from
\eqref{e:h}-\eqref{e:ha2},
Theorems \ref{c1}, \ref{T:3.10a}, \ref{T:3.10},
 \ref{T:3.13} and \ref{lambda1}.

\section{$C_\infty $-Semigroups   and Positivity }\label{S:4}

Recall that $A_0$ is the positive constant in Lemma \ref{u1}.

\begin{lem}\label{lambda2}
Suppose that $b$ is a bounded function on $\R^d\times \R^d$
satisfying condition \eqref{e:b}.
Then $\{T^b_t, t>0\}$ is a strongly continuous
  semigroup in $C_\infty(\R^d).$
\end{lem}

\pf  The  following proof is a minor modification of that for \cite[Proposition 2.3]{CKS}. For reader's convenience, we spell out
the details.
Since $q^b(t,x,y)$ is continuous by Theorem \ref{lambda1},
it follows  that $T^b_t$ maps bounded continuous functions
to continuous function for every $t>0.$
Moreover, by
\eqref{e:3.26n}
 and the semigroup property of $q^b(t,x,y)$,
 there are constants $c_1$ and $c_2$ so that
\begin{equation}\label{e:4.2a}
|q^b(t,x,y)|\leq c_1e^{c_2t}   p_{M_{b, \lambda}} (t, x, y)
\qquad \hbox{for every } t>0 \hbox{ and } x, y \in \R^d.
\end{equation}
Thus, for every $f\in C_\infty(\R^d)$ and $t>0,$
$$\lim_{x\rightarrow\infty}|T^b_tf(x)|\leq \lim_{x\rightarrow\infty}c_1e^{c_2t}
\int_{\R^d}   p_{M_{b, \lambda}}(t, x, y) | f(y)|\,dy=0
$$
and so $T^b_tf\in C_\infty(\R^d)$.
On the other hand, given $f\in C_\infty (\R^d)$,
for every $\ee >0$, there is $\delta>0$ so that
$|f(x)-f(y)|\leq \ee$ whenever $|x-y|\leq \delta$.
Since
$$
\lim_{t_0\to 0} \sup_{t\leq t_0}\sup_{x\in\R^d}\int_{|y-x|\geq \delta} |q^b(t,x,y)|\,dy
 \leq \lim_{t_0\to 0} ce^{c_2t_0}\sup_{t\leq t_0}\sup_{x\in\R^d}\int_{|y-x|\geq \delta}
   p_{M_{b, \lambda}}(t, x, y)\,dy =0,
$$
we have
\begin{eqnarray*}
&&\lim_{t\rightarrow 0}\sup_{x\in\R^d}|T^b_tf(x)-f(x)|\\
&=&\lim_{t\rightarrow 0}\sup_{x\in\R^d} \left|\int_{\R^d}q^b(t,x,y)(f(y)-f(x))\,dy \right|\\
&\leq&  \lim_{t\rightarrow 0}\sup_{x\in\R^d} \int_{|y-x|<\delta} |q^b(t,x,y)| \, |f(y)-f(x)|\,dy
+2\|f\|_\infty\lim_{t\rightarrow 0}\sup_{x\in\R^d}\int_{|y-x|\geq \delta} |q^b(t,x,y)| \,dy\\
&\leq & \ee \lim_{t\rightarrow 0}\sup_{x\in\R^d} \int_{\R^d} c_1e^{c_2t}
  p_1(t, x, y) dy =c_1 \ee.
\end{eqnarray*}
Thus $\lim_{t\rightarrow 0}\sup_{x\in\R^d}|T^b_tf(x)-f(x)|=0$.
This proves that $T^b_t$ is a strongly continuous semigroup in $C_\infty(\R^d)$.
\qed

\bigskip

\begin{lem}\label{L:4.2}
Let  $b$ be a bounded function on $\R^d\times \R^d$ satisfying
\eqref{e:5.1}.
For each $f\in C^2_\infty (\R^d)$, $\L^b f(x)$ exists
pointwise and is in $C_\infty (\R^d)$.
\end{lem}

\pf Suppose that $\gamma \in (0, 2)$ and $f\in C^2_\infty (\R^d)$.
Denote $\sum_{i, j=1}^d | \partial^2_{ij} f(x)|$ by
$|D^2 f(x)|$.
Let $R>1$ to be chosen later.
Then for each $x\in \R^d$, we have by Taylor expansion,
\begin{eqnarray*}
\Phi_f (x)&:=& \int_{\R^d} \left| f(x+z)-f(x)-\nabla f (x)
\cdot z \1_{\{|z|\leq 1\}} \right|
\frac{1}{|z|^{d+\gamma}} dz \non
&\leq & \int_{|z|\leq 1} \left| f(x+z)-f(x)-\nabla f (x) \cdot z
\1_{\{|z|\leq 1\}} \right|
\frac{1}{|z|^{d+\gamma}} dz \non
&& + \int_{1<|z|\leq R} \left| f(x+z)-f(x) \right|
\frac{1}{|z|^{d+\gamma}} dz  + \int_{|z|>R} \left| f(x+z)-f(x) \right|
\frac{1}{|z|^{d+\gamma}} dz \non
&\leq &   \sup_{|y|\leq 1} | D^2 f(x+y)| \int_{|z|\leq 1} |z|^{2-d-\gamma}dz
+ \int_{1<|z|\leq R} \left| f(x+z)-f(x) \right|
\frac{1}{|z|^{d+\gamma}} dz \non
&&+ 2 \|f \|_\infty \int_{|z|>R} |z|^{-d-\gamma} dz \non
&=& c \sup_{|y|\leq 1} | D^2 f(x+y)|+ \int_{1<|z|\leq R} \left| f(x+z)-f(x) \right|
\frac{1}{|z|^{d+\gamma}} dz
 + c R^{-\gamma}  \|f \|_\infty .
\end{eqnarray*}
For any given $\ee >0$, we can take $R$ large so that
$c R^{-\gamma}  \|f \|_\infty <\ee/2$ to conclude
that
\begin{equation}\label{e:4.3b}
\lim_{|x|\to \infty}   \int_{\R^d} \left| f(x+z)-f(x)-\nabla f (x) \cdot z \1_{\{|z|\leq 1\}} \right|
\frac{1}{|z|^{d+\gamma}} dz =0.
\end{equation}
By the same reason, applying the above argument to function
$x\mapsto f(x+y)-f(x)$ in place of $f$ yields that for every
$\ee>0$ and $x_0\in \R^d$, there is $\delta>0$ so that
\begin{equation}\label{e:4.4b}
\Phi_{f(\cdot +y)-f}(x_0)<\ee \quad \hbox{for every } |y|<\delta.
\end{equation}
It follows from the last two displays, the definition of $\L^b$ and
\eqref{e:1.4a}
that $\L^bf(x)$ exists for every $x\in \R^d$ and $\L^bf\in C_\infty (\R^d)$.
\qed

\bigskip

\noindent{\bf Proof of Theorem \ref{T:1.2}.}
 Since $b$ satisfies condition \eqref{e:5.1}, then $\L^b f\in C_\infty (\R^d)$ for every $f\in C^2_c(\R^d)$ by Lemma \ref{L:4.2}.
Let $\wh \L^b$ denote the infinitesimal generator of the strongly continuous  semigroup $\{T^b_t; t\geq 0\}$ in $C_\infty (\R^d)$, which is
a closed linear operator.
It follows from Theorem \ref{T:3.13},
 Lemmas \ref{lambda2} and \ref{L:4.2}
that for every
$f\in C^2_\infty (\R^d)$,
$(T^b_t f(x)-f(x))/t$  converges uniformly to $\L^b f(x)$ as $t\to 0$.
So
\begin{equation}\label{e:4.5b}
C^2_\infty(\R^d)\subset
D(\wh \L^b) \quad \hbox{ and } \quad
\wh \L^b f=\L^b f \quad \hbox{for } f\in C^2_\infty (\R^d).
\end{equation}
In view of Theorem \ref{c1},
there are constants $c_1, c_2>0$ so that
\eqref{e:4.2a} holds. This implies that
$$
\sup_{x\in \R^d}\int_0^\infty e^{-\lambda t} |T^b_t f|(x) dt \leq c_\lambda \|f \|_\infty,
\quad f\in C_\infty (\R^d),
$$
for every $ \lambda>c_2$. Observe that $e^{-c_2 t}T^b_t$ is a strongly continuous
semigroup in $C_\infty (\R^d)$ whose infinitesimal generator is $\wh \L^b-c_2$.
The above display
implies that $(0, \infty)$ is contained in the residual set $\rho (\wh \L^b -c_2)$ of
$\wh \L^b-c_2$.  Therefore by Theorem \ref{lambda1} and the Hille-Yosida-Ray theorem \cite[p165]{EK}, $\{e^{-c_2 t}T^b_t; t\geq 0\}$
is a positive preserving semigroup on $C_\infty (\R^d)$
if and only if $\wh \L^b-c_2$ satisfies the positive maximum principle.
On the other hand,   Courr\'ege's first theorem (see \cite[p158]{D})
tells us that $\wh \L^b-c_2$ satisfies the positive maximum principle
if and only if for each $x\in \R^d$,
$$ \frac{{\cal A}(d, -\alpha)}{|z|^{d+\alpha}}
+ \frac{{\cal A}(d, -\beta) b(x, z)}{|z|^{d+\beta}}\geq 0
\quad \hbox{for a.e. } z\in \R^d.
$$
Since $ e^{-c_2 t}T^b_t$ has a continuous integral kernel $e^{-c_2t}q^b(t, x, y)$, it follows that
$q^b(t, x, y)\geq 0$ on $(0, \infty)\times \R^d \times \R^d$
if and only if for each $x\in \R^d$,  \eqref{e:1.15} holds.
If $b(x, z)=b(x)$ is a function of $x$ only, then
by taking $|z|\to \infty$, one concludes that
 \eqref{e:1.15} holds
if and only if $b(x) \geq 0$ on $\R^d$.
  \qed

\section{Feller process and heat kernel estimates}\label{S:5}

Throughout this section, $b$ is a bounded  function
satisfying condition \eqref{e:b} and \eqref{e:1.15}.
We will show that $q^b(t, x, y) > 0$
and so it generates a Feller process $X^b$ that has strong Feller property.
We further derive the upper and lower bound estimates on $q^b(t, x, y)$.
We will first establish the Feller process $X^b$ and its connection to the martingale problem
for $(\L^b, \S(\R^d))$ under an additional assumption
\eqref{e:5.1}.
We will then remove this additional assumption using an approximation method
and the uniqueness result on $q^b(t, x, y)$ from Theorem \ref{T:3.10}.

 \bigskip
Suppose that $b$ is a bounded  function
satisfying conditions \eqref{e:b}, \eqref{e:5.1} and \eqref{e:1.15}.
Then it follows from Theorem \ref{T:1.2}, Theorem \ref{lambda1},
Lemma \ref{lambda2} and Theorem \ref{T:3.9}, $T^b$ is a Feller semigroup.
So it uniquely determines a conservative
Feller process $X^b=\{X^b_t, t\geq 0, \P_x, x\in \R^d\}$ having
$q^b(t, x, y)$ as its transition density function.
Since, by Theorem \ref{T:3.10},
$q^b(t, x, y)$ is continuous and $q^b(t, x, y)\leq c_1 e^{c_2t} p_{M_{b, \lambda}} (t, x, y)$
for some positive constants $c_1$ and $c_2$,
$X^b$ enjoys the strong Feller property.

\bigskip

\begin{prp}\label{p2}
Suppose that $b$ is a bounded  function
satisfying conditions \eqref{e:b}, \eqref{e:5.1} and \eqref{e:1.15}.
For each $x\in\R^d$ and $f\in C_b^2 (\R^d),$
$$M_t^f:=f(X^b_t)-f(X^b_0)-\int_0^t \L^b f(X^b_s)\,ds$$
is a martingale under $\P_x$.  So in particular,
the Feller process $(X^b, \P_x, x\in \R^d)$ solves the martingale
problem for $(\L^b, C^2_\infty (\R^d))$.
\end{prp}

\pf This follows immediately from
Theorem \ref{T:3.13}
and the Markov property of $X^b$. \qed

We next determine the L\'evy system of $X^b$. Recall that
\begin{equation}\label{e:5.2}
J^b(x, y)= \frac{{\cal A}(d, -\alpha)}{|x-y|^{d+\alpha}}
+ \frac{{\cal A}(d, -\beta) \, b(x, y-x)}{|x-y|^{d+\beta}}.
\end{equation}

\begin{prp}\label{p3}
Suppose that $b$ is a bounded  function
satisfying conditions \eqref{e:b}, \eqref{e:5.1} and \eqref{e:1.15}.
Assume that $A$ and $B$ are disjoint compact sets in $\R^d$.
Then
$$
\sum_{s\leq t} \1_{\{X^b_{s-}\in A, X^b_s\in B\}}-
\int_0^t \1_A(X^b_s) \int_B J^b(X^b_s, y) dy\,ds
$$
is a $\P_x$-martingale for each $x\in \R^d$.
\end{prp}

\pf  The proof is similar to that for \cite[Theorem 2.6]{CKS}.
For reader's convenience, we give the details here.
 Let $f\in
C^\infty(\R^d)$ with $f=0$ in an open neighborhood of $A$
and $f=1$ in an open neighborhood of  $B$.
Define
$$
M_t^f:=f(X^b_t)-f(X^b_0)-\int_0^t \L^bf(X^b_s)\,ds .
$$
Then $M_t^f$ is a martingale under $\P_x$ by Proposition \ref{p2},
and so is $N_t^f:=\int_0^t \1_A(X^b_{s-})dM_s^f$.
 Proposition \ref{p2} in particular implies that $X^b_t=(X^{b,
 1}_t, \dots, X^{b, d}_t) $ is a semi-martingale. So by Ito's formula,
we have that,
\begin{equation}\label{e:ito}
f(X^b_t)-f(X^b_0)=\sum^d_{i=1}\int^t_0{\partial}_if(X^b_{s-})dX^{b,
i}_s +\sum_{s\le t}\eta_s(f) +\frac12 A_t(f),
\end{equation}
where
\begin{equation}\label{e:ito2}
\eta_s(f)=f(X^b_s)-f(X^b_{s-})-\sum^d_{i=1}{\partial}_if(X^b_{s-})(X^{b,
i}_s-X^{b, i}_{s-})
\end{equation}
and
\begin{equation}\label{e:ito3}
A_t(f)=\sum^d_{i, j=1}\int^t_0{\partial}_i{\partial}_jf(X^b_{s-})d
\langle (X^{b, i})^c, (X^{b, j})^c\rangle_s.
\end{equation}
Since $f$  vanishes in an open neighborhood of $A$,
we have by
  \eqref{e:ito}--\eqref{e:ito3}, \eqref{e:1.2n} and \eqref{e:1.1} that
\begin{eqnarray*}
N^f_t&=&\sum_{s\le t}{\bf 1}_A(X^b_{s-})f(X^b_s) -
\int^t_0{\bf 1}_A(X^b_s)\left(
\L^b f(X^b_s)\right)ds\\
&=&\sum_{s\le
t}{\bf 1}_A(X^b_{s-})f(X^b_s)-\int^t_0{\bf 1}_A(X^b_s)\int_{\R^d}f(y)J^b(X^b_s,
y)dyds.
\end{eqnarray*}
By taking a sequence of functions $f_n\in C^\infty_c(\R^d)$ with
$f_n=0$ on $A$,  $f_n=1$ on $B$ and $f_n\downarrow {\bf 1}_B$, we get
that, for any $x\in \R^d$,
$$
\sum_{s\le t}{\bf 1}_A(X^b_{s-}){\bf 1}_B(X^b_s) -\int^t_0{\bf
1}_A(X^b_s)\int_B J^b(X^b_s, y)dyds
$$
is a $\P_x$-martingale   for every $x\in \R^d$. \qed

Proposition \ref{p3} implies that
$$
\E_x\left[ \sum_{s\le t}{\bf 1}_A(X^b_{s-}){\bf 1}_B(X^b_s)
\right]=
\E_x\left[\int^t_0\int_{\R^d} {\bf 1}_A(X^b_s){\bf 1}_B(y)J^b(X^b_s, y)dyds\right].
$$
Using this and a routine measure theoretic argument, we get
$$\E_x\left[ \sum_{s\le
t}f(s, X^b_{s-}, X^b_s) \right]
=\E_x\left[\int^t_0\int_{\R^d}f(s, X^b_s, y)J^b(X^b_s, y)dyds\right]
$$
for any non-negative measurable function $f$ on $(0, \infty)
\times \R^d\times \R^d$
vanishing on $\{(x, y)\in \R^d\times \R^d: x=y\}$. Finally,
following the same arguments as in \cite[Lemma 4.7]{CK} and
\cite[Appendix A]{CK2}, we get

\begin{prp}\label{p4}
Suppose that $b$ is a bounded  function
satisfying conditions \eqref{e:b}, \eqref{e:5.1} and \eqref{e:1.15}.
Let $f$ be a nonnegative function on $\R_+\times
\R^d\times\R^d$ vanishing on the diagonal. Then for stopping time $T$
with respect to the minimal admissible filtration generated by $X^b$,
$$
\E_x \left[ \sum_{s\leq T}f(s, X^b_{s-}, X^b_s) \right]
=\E_x \left[ \int_0^T\int_{\R^d}
f(s,X^b_s,u) J^b(X^b_s,u) \,du\,ds \right].
$$
\end{prp}

 \bigskip

 To remove the assumption \eqref{e:5.1} on $b$, we approximate
 a general measurable function $b(x, z)$ by continuous $k_n(x, z)$.
 To show that $q^{k_n}(t, x, y)$ converges to $q^b (t, x, y)$,
 we establish equi-continuity of $q^b(t, x, y)$ and apply
 the uniqueness result, Theorem \ref{T:3.10}.

 \bigskip

\begin{prp}\label{H1}
 For each  $0<t_0<T<\infty$ and $A>0$,
 the function $q^b(t,x,y)$ is uniformly continuous in $(t,x)\in (t_0,T)\times\R^d$ for
every $b$ with $\|b\|_\infty\leq A$ that satisfies \eqref{e:b} and for
all $y\in \R^d.$
\end{prp}

\pf In view of Theorem \ref{T:3.11},
it suffices to prove the theorem
for $A=A_0$, where $A_0$ is the constant in Lemma \ref{u1}
(or in Theorem \ref{T:1.1}).
Using the Chapman-Kolmogorov equation for $q^b(t, x, y)$ (see Lemma \ref{lq^b})
and \eqref{e:3.26n}, it suffices to prove the Proposition for $T=1$.

Noting that by \eqref{e:qbn2}
$$q^b_n(t,x,y)=\int_0^t\int_{\R^d} p_0(t-r, x, z) \S^b_z q^b_{n-1} (r, z, y)\,dz\,dr.$$
Hence, for $T>t>s>t_0, x_1,x_2\in\R^d$ and $y\in\R^d,$ we have
$$\begin{aligned}
&|q^b_n(s,x_1,y)-q^b_n(t,x_2,y)|\\
\leq & \int_0^s\int_{\R^d}|p_0(s-r, x_1, z)-p_0(t-r, x_2, z)||\S^b_z q^b_{n-1} (r, z, y)|\,dz\,dr\\
&+\int_s^t\int_{\R^d} p_0(t-r, x_2, z)|\S^b_z q^b_{n-1} (r, z, y)|\,dz\,dr\\
=:& I+II.
\end{aligned}$$
It is known (see \cite{CK}) that there are positive constants $c_1$ and $\theta$ so that
 for any $ t,s\in [t_0,T]$
and $x_i\in \R^d$ with $i=1,2,$
$$
|p_0(s,x_1,y)-p_0(t, x_2, y)|\leq
c_1\,  t_0^{-(d+\theta)/\alpha}
\left(|t-s|^{1/\alpha}+|x_1-x_2|\right)^\theta,
\quad y\in\R^d,
$$
we have by \eqref{e:3.9}, \eqref{e:3.15} and \eqref{e:f0}, for $\rho\in (0,s/2),$
\begin{equation}\label{e:h1}\begin{aligned}
I=&\int_0^{s-\rho}\int_{\R^d}|p_0(s-r, x_1, z)-p_0(t-r, x_2, z)||\S^b_z q^b_{n-1} (r, z, y)|\,dz\,dr\\
&+\int_{s-\rho}^s\int_{\R^d}|p_0(s-r, x_1, z)-p_0(t-r, x_2, z)||\S^b_z q^b_{n-1} (r, z, y)|\,dz\,dr\\
\leq &c_2 2^{-(n-1)}\rho^{-(d+\theta)/\alpha}\left(|t-s|^{1/\alpha}+|x_1-x_2|\right)^\theta\int_0^{s-\rho}\int_{\R^d} f_0(r,z,y)\,dz\,dr\\
&+c_2 2^{-(n-1)}\int_{s-\rho}^s\int_{\R^d}(p_0(s-r, x_1, z)+p_0(t-r, x_2, z))f_0(r,z,y)\,dz\,dr\\
\leq &c_3 2^{-(n-1)}\rho^{-(d+\theta)/\alpha}\left(|t-s|^{1/\alpha}+|x_1-x_2|\right)^\theta s^{1-\beta/\alpha}+c_3 2^{-(n-1)}(s-\rho)^{-(d+\beta)/\alpha}\rho.
\end{aligned}\end{equation}
Moreover, by \eqref{e:3.9} and \eqref{e:3.15},
\begin{equation}\label{e:h2}
II \leq 2^{-(n-1)}\int_s^t\int_{\R^d} p_0(t-r, x_2, z)f_0(r,z,y)\,dz\,dr
 \leq 2^{-(n-1)}s^{-(d+\beta)/\alpha}|t-s|.
 \end{equation}
Therefore, noting that
$$
|q^b(s,x_1,y)-q^b(t,x_2,y)|
\leq |p_0(s,x_1,y)-p_0(t,x_2,y)|+\sum_{n=1}^\infty|q^b_n(s, x_1, y)-q^b_n(t, x_2, y)|,
$$
then  first taking $|t-s|$ and $|x_1-x_2|$ small, and then making $\rho$ small in \eqref{e:h1} and \eqref{e:h2}
 yields the conclusion of this Proposition.
\qed

\begin{prp}\label{H2}
For each  $0<t_0<T<\infty$ and $A>0$,
 the function $q^b(t,x,y)$ is uniformly continuous in $y$ for
every $b$ with $\|b\|_\infty\leq A$ that satisfies \eqref{e:b}
 and for
all $(t,x)\in (t_0, T)\times \R^d.$
\end{prp}

\pf  In view of Theorem \ref{T:3.11},
it suffices to prove the theorem
for $A=A_0$, where $A_0$ is the constant in Lemma \ref{u1}
(or in Theorem \ref{T:1.1}).
Using the Chapman-Kolmogorov equation for $q^b(t, x, y)$ (see Lemma \ref{lq^b})
and \eqref{e:3.26n}, it suffices to prove the Proposition for $T=1$.

Define $P(s, x, y)=p_0(s,x)-p_0(s,y).$ For $s>0,$ we have
\begin{equation}\label{Lp}\begin{aligned}
&|\S^b p_0(s, y_1)-\S^b p_0(s, y_2)|\\
\leq &c_1\int_{\R^d}|P(s,y_1+h, y_2+h)-P(s, y_1, y_2)-\langle\nabla_{(y_1,y_2)} P(s, y_1, y_2), h
\1_{|h|\leq 1}\rangle|\frac{dh}{|h|^{d+\beta}}\\
\leq &c_1\int_{|h|\leq 1}|h|^2\sup_{\theta\in(0,1)}|\frac{\partial^2}{\partial y_1^2}p_0(s,y_1+\theta h)
-\frac{\partial^2}{\partial y_2^2}p_0(s,y_2+\theta h)|\frac{dh}{|h|^{d+\beta}}\\
&+c_1\int_{|h|>1} |p_0(s,y_1+h)-p_0(s, y_2+h)-p_0(s, y_1)+p_0(s,y_2)|\frac{dh}{|h|^{d+\beta}}\\
\leq &c_2 \sup_{y}|\frac{\partial^3}{\partial y^3}p_0(s,y)||y_1-y_2|\int_{|h|\leq 1}|h|^2\frac{dh}{|h|^{d+\beta}}
+c_2\sup_{y}|\frac{\partial}{\partial y}p_0(s,y)||y_1-y_2|\int_{|h|> 1}\frac{dh}{|h|^{d+\beta}}\\
\leq &c_3|y_1-y_2|[s^{-(d+3)/\alpha}+s^{-(d+1)/\alpha}],\\
\end{aligned}\end{equation}
where in the fourth inequality, $|\frac{\partial^3}{\partial y^3}p_0(s,y)|\leq c_3 s^{-(d+3)/\alpha}$ can be proved similarly by the argument in Lemma \ref{0}.
Take $\rho \in (0, t_0/2)$.
 Then for each $n\geq 1,$
 we have  by
  \eqref{e:1.6}, \eqref{e:f0}, Lemma \ref{n1}, Lemma \ref{u1} and (\ref{Lp}),
  that  for $(t, x, y)\in (t_0, 1)\times \R^d\times \R^d$,
\begin{eqnarray*}
&&|q_n^b(t,x,y_1)-q_n^b(t,x,y_2)|\\
&\leq &\int_0^\rho\int_{\R^d}q_{n-1}^b(t-s,x,z)|\S_z^b p_0(s,z,y_1)-
\S_z^b p_0(s,z,y_2)|\,dz\,ds\\
&&+\int_\rho^t\int_{\R^d}q_{n-1}^b(t-s,x,z)|\S_z^b p_0(s,z,y_1)-\S_z^b p_0(s,z,y_2)|\,dz\,ds\\
&\leq &c_4 2^{-(n-1)}\int_0^\rho\int_{\R^d}  p_1
(t-s,x,z)|\S_z^b p_0(s,z,y_1)-\S_z^b p_0(s,z,y_2)|\,dz\,ds\\
&& +c_4 2^{-(n-1)} \int_\rho^t\int_{\R^d}
  p_1 (t-s,x,z)
\left| \S_z^b p_0(s,z-y_1)-\S_z^b p_0(s,z-y_2)\right|\,dz\,ds\\
&\leq &c_5 2^{-(n-1)}t_0^{-d/\alpha}\int_0^\rho\int_{\R^d} \left(|\S_z^b p_0(s,z,y_1)|+|\S_z^b p_0(s,z,y_2)|\right)\,dz\,ds\\
&&+ c_5 2^{-(n-1)}\rho^{-(d+3)/\alpha}|y_1-y_2|\int_\rho^t\int_{\R^d}
  p_1 (t-s,x,z)\,dz\,ds\\
&\leq & c_6 \, 2^{-(n-1)}\, t_0^{-d/\alpha}\rho^{1-\beta/\alpha}
+c_6 2^{-(n-1)}\rho^{-(d+3)/\alpha}|y_1-y_2| .
\end{eqnarray*}
Therefore we have
$$\begin{aligned} &|q^b(t, x, y_1)-q^b(t, x, y_2)|\\
 \leq  &|p_0(t, x, y_1)-p_0(t, x, y_2)|+
 \sum_{n=1}^\infty c_62^{-(n-1)} \,   t_0^{-d/\alpha}\rho^{1-\beta/\alpha}
+\sum_{n=1}^\infty c_62^{-(n-1)}  \rho^{-(d+3)/\alpha}|y_1-y_2|.
\end{aligned}$$
 By first taking $|y_1-y_2|$ small and then making $\rho$ small
 yields the desired uniform continuity of $q^b(t, x, y)$.
 \qed

\begin{thm}\label{T:m1}
Suppose $b$ is a bounded function on $\R^d\times\R^d$
satisfying \eqref{e:b} and \eqref{e:1.15}.
The kernel $q^b(t, x, y)$ uniquely determines a Feller process $X^b
=(X^b_t, t\geq 0, \P_x, x\in \R^d)$
on the canonical Skorokhod space
${\mathbb D}([0, \infty), \R^d)$
such that
$$ \E_x \left[ f(X^b_t)\right] =\int_{\R^d} q^b (t, x, y) f(y) dy
$$
for every bounded continuous function $f$ on $\R^d$.
The Feller process $X^b$ is conservative
 and has a L\'evy system $(J^b(x, y)dy, t)$,
where
$$ J^b(x, y)= \frac{{\cal A}(d, -\alpha)}{|x-y|^{d+\alpha}}
+ \frac{{\cal A}(d, -\beta) \, b(x, y-x)}{|x-y|^{d+\beta}}.
$$
Moreover, for each $x\in \R^d$, $(X^b, \P_x)$ is the unique
solution to the martingale problem $(\L^b, {\cal S}(\R^d ))$
with initial value $x$. Here ${\cal S} (\R^d)$ denotes the
space of tempered functions on $\R^d$.
\end{thm}
\pf When $b$ is a bounded function satisfying \eqref{e:b}, \eqref{e:5.1} and \eqref{e:1.15},
 the theorem has already been established via
 Propositions \ref{p2}-\ref{p4}.
We now remove the assumption \eqref{e:5.1}.
Suppose that $b(x, z)$ is a bounded
function that satisfies \eqref{e:b} and \eqref{e:1.15}.
Let  $\varphi$ be a non-negative smooth function with compact
support in $\R^d$ so that
$\int_{\R^d}\varphi (x) dx =1$. For each $n\geq 1$,
define $\varphi_n(x)=n^d  \varphi(nx)$ and
$$ k_n (x, z):= \int_{\R^d} \varphi_n (x-y) b(y, z) dy.
$$
Then $k_n$ is a function that satisfies \eqref{e:b}, \eqref{e:5.1} and \eqref{e:1.15}
with $\| k_n\|_\infty \leq \| b \|_\infty$.
By Theorem \ref{T:1.1}, Proposition \ref{H1} and Proposition \ref{H2}, $q^{k_n}(t,x,y)$ is uniformly bounded
and equi-continuous on $[1/M, M]\times\R^d\times\R^d$ for each $M\geq 1$, then
there is a subsequence $\{n_j\}$ of $\{n\}$ so that $q^{k_{n_j}}(t,x,y)$ converges  boundedly and uniformly on compacts of $(0, \infty)\times \R^d\times \R^d$,
to some continuous function $\overline q (t,x,y)$, which again satisfies \eqref{e:1.7}. Obviously, $\overline q(t,x,y)$ also satisfies the Chapman-Kolmogorov equation and $\int_{\R^d}\overline q(t,x,y)\,dy=1.$
By \eqref{e:3.21} and Theorem \ref{c1},
$$
q^{k_{n_j}}(t, x, y)
=p_0(t, x, y)+\int_0^t \int_{\R^d} q^{k_{n_j}} (t-s, x, z)
\S^{k_{n_j}}_z p_0(s, z, y)dz ds
$$
and
$$  q^{k_{n_j}}(t, x, y)\leq c\, g_{M_{b, \lambda}}(t, x, y)
$$
for every $0<t\leq 1\wedge(A_0/\|b\|_\infty)^{\alpha/(\alpha -\beta )}$
and $x, y\in \R^d$, where $c$ is a positive constant that depends only
on $d$, $\alpha, \beta$ and $\|b\|_\infty$. Letting $j\to \infty$, we have
by \eqref{L}, Lemma \ref{2'} and the dominated convergence theorem
that
$$
\overline q(t, x, y)=p_0 (t, x, y)+\int_0^t \int_{\R^d} \overline q (t-s, x, z)
\S^b_z p_0 (s, z, y)dy ds
$$
and $  \overline q(t, x, y)\leq c\, g_{M_{b, \lambda}}(t, x, y)$
for every $0<t\leq 1\wedge(A_0/\|b\|_\infty)^{\alpha/(\alpha -\beta )}$
and $x, y\in \R^d$. Hence we conclude from Theorem \ref{T:3.10}
that $\overline q(t, x, y)= q^b(t, x, y)$. This in particular implies
that $q^b(t, x, y)\geq 0$. So there is a Feller process $X^b$ having $q^b(t, x,y)$
as its transition density function.
The proof of Propositions \ref{p2}-\ref{p4}
only uses the condition (\ref{e:5.1}) through its implication
 that $q^b(t, x, y)\geq 0$. So in view of what we just established,
Propositions \ref{p2}-\ref{p4}  continue to
hold for $X^b$ under the current setting without the additional assumption \eqref{e:5.1}.
The non-local operator $\L^b$ satisfies the assumptions $[A_1]$ and $[A_2]$
of \cite{Ko2}. So by \cite[Theorem 3]{Ko2},  solution to
the martingale problem $(\L^b, {\cal S}(\R^d))$ is unique.
Since ${\cal S}(\R^d) \subset C^2_\infty (\R^d)$,
  the proof of the theorem is now complete.
\qed

\medskip
For each $\lambda>0,$ define
\begin{equation}\label{e:hatb}
\wh{b}_\lambda(x,z)=b(x,z)1_{\{|z|\leq \lambda\}}(z)+b^{+}(x,z)1_{\{|z|>\lambda\}}(z).
\end{equation}
In the following,
we use a method of Meyer \cite{Mey} to construct from $X^b$, by adding
suitable jumps, a strong Markov process $Y$
 corresponding to the jumping kernel $J^{\wh{b}_\lambda}$
defined by \eqref{e:Jb} but with $\wh{b}_\lambda$ in place of $b.$
Define $$\mathcal{J}(x)=\int_{\R^d}(J^{\wh{b}_\lambda}(x,y)-J^b(x,y))\,dy.$$
Then there exists a positive constant $c_1$ so that $0\leq\mathcal{J}(x)\leq c_1$ for all $x\in\R^d.$
Let $$q(x,y)=\dfrac{J^{\wh{b}_\lambda}(x,y)-J^{b}(x,y)}{\mathcal{J}(x)}.$$
Let $S_1$ be an exponential random variable of parameter $1$ independent of $X^b.$
Set
\begin{equation}\label{e:u1}
C_t=\int_0^t \mathcal{J} (X^b_s)\,ds, \quad U_1=\inf\{t\geq 0: C_t\geq S_1\}.
\end{equation}
We let $Y_t=X^b_t$ for $0\leq t<U_1$ and  define $Y_{U_1}$ with law $q(Y_{U_{1-}}, \cdot)=q(X^b_{U_{1-}}, \cdot),$
and then repeat using an independent exponential random variable $S_2$ to define $U_2,$ etc.
So the construction proceeds now in the same way from the new starting point $(U_1, Y_{U_1}).$
Since $\mathcal{J}(x)$ is bounded, only finitely many new jumps are introduced in any bounded time interval.
In \cite{Mey}, it is proved that the resulting process $Y$ is a strong Markov process.
By slightly abusing the notation, we still use $\P_x$ and $\E_x$ to denote the  above constructed probability law and expectation
induced on such enlarged probability space under which $Y_0=x.$

\medskip

\begin{lem}\label{L:eta}
For each $x\in\R^d$ and $f\in C_b^2 (\R^d),$
$$ \E_x \left[ f(Y_t); t<U_1 \right] =f(x)+
\E_x \left[ \int_0^t \left( \L^b - {\cal J}(Y_s)\right) f(Y_s) \1_{\{s<U_1\}} ds \right],
$$
\end{lem}
\pf
By the definition of  $U_1$ and Ito's formula, for each function $f\in C_b^2(\R^d),$
 \begin{eqnarray*}
 \E_x \left[ f(Y_t); t<U_1 \right]
&=&\E_x \left[ f(X^b_t) \1_{\{U_1>t\}}\right]
=\E_x \left[ f(X^b_t)e^{-C_t} \right]\\
&=&f(x)+\E_x\left[\int_0^t(\L^b-\J(X^b_s))f(X^b_s)e^{-C_s}\,ds\right]\\
&=&f(x)+\E_x \left[ \int_0^t \left( \L^b - {\cal J}(Y_s)\right) f(Y_s) \1_{\{s<U_1\}} ds \right].
\end{eqnarray*}
\qed

\begin{prp}\label{p2Y}
For each $x\in\R^d$ and $f\in C_b^2 (\R^d),$
$$M_t^f:=f(Y_t)-f(Y_0)-\int_0^t \L^{\wh{b}_\lambda} f(Y_s)\,ds$$
is a martingale under $\P_x$.  So in particular,
the strongly Markov process $(Y, \P_x, x\in \R^d)$ solves the martingale
problem for $(\L^{\wh{b}_\lambda}, C^2_\infty (\R^d))$.
\end{prp}

\pf
Note that $M^f_t$ is an additive function of $Y$.
So by the Markov property of $Y$, it suffices to show
that $\E_x \left[ M^f_t \right]=0$ for every $x\in \R^d$
and $t>0$.

Recall that $U_1$ is defined in \eqref{e:u1}, and denote by
$\{U_n, n\geq 2\}$  the subsequent jump adding times inductively defined according to the construction of Meyer \cite{Mey}.
For every $\alpha >0$, set
 $u_\alpha(x)=\E_x \left[ \int_0^{U_1} e^{-\alpha t} f(Y_t) dt\right]$.
Since by Lemma \ref{L:eta},
$$ \E_x \left[ f(Y_t); t<U_1 \right] =f(x)+
\E_x \left[ \int_0^t \left( \L^b - {\cal J}(Y_s)\right) f(Y_s) \1_{\{s<U_1\}} ds \right],
$$
we have by Fubini theorem that
$$
 u_\alpha (x)=  \frac{f(x)}\alpha
  + \frac1{\alpha} \E_x \left[
\int_0^{U_1} e^{-\alpha s } (\L^b -{\cal J}(Y_s)) f (Y_s) ds \right].
$$
Observe that in view of \cite[p.286]{Sh}
(see, for example, the proof of \cite[Proposition 2.2]{CZ}),
for any non-negative function
$\varphi$ on $\R^d$ and $x\in \R^d$,
$$ \E_x \left[ e^{-\alpha U_1} \varphi (Y_{U_1-})\right]
=\E_x \left[ \int_0^{U_1} e^{-\alpha s}
{\cal J}(Y_s) \varphi (Y_s) ds \right].
$$
Set $U_0=0$ and let $\theta_t$ to denote the time shift operator
for the Markov process $Y$.
Then we have from above and the strong Markov
property of $Y$ that
\begin{eqnarray*}
&& \E_x \left[ \int_0^\infty e^{-\alpha t} f(Y_t) dt \right]
= \sum_{j=0}^\infty
\E_x \left[ \int_{U_j}^{U_{j+1}} e^{-\alpha t} f(Y_t) dt \right]
 =  \sum_{j=0}^\infty
\E_x \left[ e^{-\alpha U_j}  u_\alpha ({Y_{U_j}}) \right] \\
&=&  \frac{f(x)}{\alpha} + \frac1{\alpha}\sum_{j=1}^\infty
\E_x \left[ e^{-\alpha U_j} f({Y_{U_j}})\right]
+ \frac1{\alpha} \sum_{j=0}^\infty
\E_x \left[ \int_{U_j}^{U_{j+1} } e^{-\alpha s } (\L^b -{\cal J}(Y_s)) f (Y_s) ds \right] \\
&=&\frac{f(x)}{\alpha} + \frac1{\alpha}\sum_{j=1}^\infty
\E_x \left[ e^{-\alpha U_j}  \int_{\R^d}
f( y) q ( {Y_{U_j -}}, y) dy \right] + \frac1{\alpha}
 \E_x \left[ \int_0^\infty e^{-\alpha s } (\L^b -{\cal J}(Y_s))
  f (Y_s) ds \right]
\\
&=& \frac{f(x)}{\alpha} + \frac1{\alpha}\sum_{j=1}^\infty
\E_x \left[ e^{-\alpha U_{j-1}}
\int_{\R^d} f( y) \left(e^{-\alpha U_1} q ( {Y_{U_1 -}}, y)\right) \circ \theta_{U_{j-1}} dy \right] \\
&&  + \frac1{\alpha}
 \E_x \left[ \int_0^\infty e^{-\alpha s } (\L^b -{\cal J}(Y_s))
  f (Y_s) ds \right]
\\
&=& \frac{f(x)}{\alpha} + \frac1{\alpha}\sum_{j=1}^\infty
\E_x \left[ e^{-\alpha U_{j-1}}
\int_{\R^d} f( y) \left( \int_0^{U_1} e^{-\alpha s} {\cal J}(Y_s)
 q (Y_s, y)ds \right) \circ \theta_{U_{j-1}} dy \right] \\
&&  + \frac1{\alpha}
 \E_x \left[ \int_0^\infty e^{-\alpha s } (\L^b -{\cal J}(Y_s))
  f (Y_s) ds \right] \\
&=& \frac{f(x)}{\alpha} + \frac1{\alpha}
\E_x \left[
\int_{\R^d} f( y) \left( \int_0^\infty e^{-\alpha s} {\cal J}(Y_s)
 q (Y_s, y)ds \right)   dy \right]\\
&& + \frac1{\alpha}
 \E_x \left[ \int_0^\infty e^{-\alpha s } (\L^b -{\cal J}(Y_s))
  f (Y_s) ds \right] \\
&=& \frac{f(x)}{\alpha} + \frac1{\alpha}
 \E_x \left[ \int_0^\infty e^{-\alpha s }
 \left( \L^b f(Y_s) + \int_{\R^d} {\cal J}(Y_s)  q (Y_s, y)
 ( f(y)-f (Y_s)) dy\right) ds \right] \\
&=& \frac{f(x)}{\alpha} + \frac1{\alpha}
 \E_x \left[ \int_0^\infty e^{-\alpha s }
  \L^{\wh b_\lambda}  f (Y_s)  ds \right].
\end{eqnarray*}
By the uniqueness of the Laplace transform,
we conclude from above that $\E_x \left[ M^f_t \right]=0$
for all $t\geq 0$ and $x\in \R^d$. \qed

\medskip

Note that $\wh{b}_\lambda$ defined by \eqref{e:hatb}
is a bounded function on $\R^d\times\R^d$
satisfying \eqref{e:b} and \eqref{e:1.15}. By Theorem \ref{T:m1}, the kernel $q^{\wh{b}_\lambda}(t, x, y)$ uniquely determines a Feller process
$X^{\wh{b}_\lambda}=(X^{\wh{b}_\lambda}_t, t\geq 0, \P_x, x\in \R^d)$
on the canonical Skorokhod space
${\mathbb D}([0, \infty), \R^d),$ and $(X^{\wh{b}_\lambda}, \P_x)$ is the unique solution to
the martingale problem for $(\L^{\wh{b}_\lambda}, {\cal S}(\R^d ))$
with initial value $x$. This, together with
Proposition \ref{p2Y} implies that the process $Y$ coincides with $X^{\wh{b}_\lambda}$ in the sense of distribution.

\begin{thm}\label{T:u1}
\underline{}For every $\lambda>0$ and $A>0$, there is a positive constant
$C_{15}=C_{15}(d,\alpha,\beta, A, \lambda)$
such that
for any bounded $b$ satisfying \eqref{e:b} and \eqref{e:1.15}
with $\| b\|_\infty\leq A$,
$$
 q^b(t,x,y)
\leq C_{15} {p}_{M_{b^+, \lambda}} (t,x,y)  \quad \hbox{for }
t\in (0,1] \hbox{ and }  x,y\in\R^d.
$$
\end{thm}

\pf Noting that  $\wh{b}_\lambda$ is a bounded function on $\R^d\times\R^d$ with $\|\wh{b}_\lambda\|_\infty\leq \|b\|_\infty$
satisfying \eqref{e:b} and \eqref{e:1.15}, then  by Theorem \ref{T:1.1}, there is a positive constant $C=C(d,\alpha, \beta, A, \lambda)$
so that
\begin{equation}\label{e:b1}
 q^{\wh{b}_\lambda}(t,x,y)
\leq C {p}_{M_{b^+, \lambda}} (t,x,y)  \quad \hbox{for }
t\in (0,1] \hbox{ and }  x,y\in\R^d.
\end{equation}
Let $\{\M_t\}_{t\geq 0}$ be the filtration generated by $X^b$.
Note that  $X^{\wh{b}_\lambda}$ has the same distribution as $Y.$
Then  by Lemma 3.6 in \cite{BBCK}, for any $A\in\M_t,$
$$\P^x(X^{\wh{b}_\lambda}_t\in A)=\P^x(Y_t\in A)\geq \P^x(\{Y_s=X^b_s \, \mbox{for all}\, 0\leq s\leq t\}\cap A)\geq e^{-t\|\J\|_\infty}\P^x(X^b_t\in A).$$
Hence, by \eqref{e:b1}
$$q^b(t,x,y)\leq e^{\|\J\|_\infty} q^{\wh{b}_\lambda}(t,x,y)\leq Ce^{\|\J\|_\infty} {p}_{M_{b^+, \lambda}} (t,x,y)  \quad \hbox{for }
t\in (0,1] \hbox{ and }  x,y\in\R^d.$$
\qed

\medskip

For a Borel set $B\subset \R^d$, we define
 $\tau^b_B=\inf\{t>0: X^b_t\notin B\}$ and
  $\sigma^b_B:=\inf\{t\geq 0: X^b_t\in B\}.$

\begin{prp}\label{y1}
For each $A>0$ and $R_0>0$, there exists
a positive constant
$$
\kappa =\kappa (d, \alpha, \beta, A, R_0)< 2^{\alpha}
\left( 1-(1/3)^\alpha\right)
$$
 so that
for every $b$ satisfying \eqref{e:b} and \eqref{e:1.15}
with $\|b\|_\infty \leq A$, $r\in (0, R_0]$
and $x\in \R^d$,
$$
\P_x \left(\tau^b_{B(x,r)}\leq \kappa r^{\alpha} \right)
\leq \dfrac{1}{2} .
$$
\end{prp}

\pf  Let $f$ be a $C^2$ function taking values in $[0,
1]$ such that $f(0)=0$ and $f(u)=1$ if $|u|\geq 1.$ Set
$f_{x, r}(y)=f(\frac{y-x}{r})$.  Note that $f_{x, r}$ is a $C^2$
function taking values in $[0, 1]$ such that $f_{x,r}(x)=0$ and
$f_{x, r}(y)=1$ if $y\notin B(x, r)$. Moreover,
$$
\sup_{y\in\R^d}\left|\frac{\partial^{2}f_{x,r}(y)}{\partial y_i\partial y_j}\right| \leq r^{-2}\,
\sup_{y\in\R^d}\left|\frac{\partial^{2}f (y)}{\partial
y_i\partial y_j}\right|.
$$
Denote $\sum_{i, j=1}^d | \partial^2_{ij} f(x)|$ by
$|D^2 f(x)|$.
By Taylor's formula,
it follows that \begin{equation}\label{10}\begin{aligned} |\L^b
f_{x,r}(u)| &\leq c_1\int \left|f_{x,r}(u+h)-f_{x,r}(u)-\langle\nabla
f_{x,r}(u),h\rangle\1_{\{|h|\leq r\}}\right|
\left(\dfrac{1}{|h|^{d+\alpha}}+\frac{1}{|h|^{d+\beta}}\right)\,dh\\
&= c_1\int_{\{|h|\leq
r\}}\left|f_{x,r}(u+h)-f_{x,r}(u)-\langle\nabla
f_{x,r}(u),h\rangle\right|
\left(\dfrac{1}{|h|^{d+\alpha}}+\frac{1}{|h|^{d+\beta}}\right)\,dh\\
&\quad+c_1\int_{\{|h|>
r\}}\left|f_{x,r}(u+h)-f_{x,r}(u)\right| \left(\dfrac{1}{|h|^{d+\alpha}}+\frac{1}{|h|^{d+\beta}}\right)\,dh\\
&\leq c_2\|D^2f\|_\infty r^{-2}\int_{|h|\leq r} |h|^2\left(\dfrac{1}{|h|^{d+\alpha}}+\frac{1}{|h|^{d+\beta}}\right)\,dh\\
&\quad +c_2\|f\|_\infty \int_{\{|h|>
r\}}\left(\dfrac{1}{|h|^{d+\alpha}}+\frac{1}{|h|^{d+\beta}}\right)\,dh\\
&\leq  c_3 (r^{-\alpha}+r^{-\beta})
\leq c_3(1+R_0^{\alpha-\beta})r^{-\alpha},
\end{aligned}\end{equation}
where $c_i=c_i(d,\alpha,\beta, A), i=1,2,3$ are positive constants.
Therefore, for each $t>0,$
$$\begin{aligned}
\P_x(\tau^b_{B(x,r)}\leq t)
&\leq \E_x \left[f_{x,r}(X^b_{\tau^b_{B(x,r)}\wedge t})\right]-f_{x,r}(x)\\
& =  \E_x \left[ \int_0^{\tau^b_{B(x,r)}\wedge t} \L^b f_{x,r}(X^b_s)\,ds \right]
\leq   c_3 \big( 1+R_0^{\alpha-\beta}\big)\dfrac{t}{r^{\alpha}}.
\end{aligned}$$
Set $\kappa=(2^{\alpha}[1-(1/3)^\alpha])\wedge (2c_3(1+R_0^{\alpha-\beta}))^{-1},$ then
$$\P_x(\tau^b_{B(x,r)}\leq \kappa r^{\alpha})\leq \frac{1}{2}.$$
\qed

Recall that $m_{b,\lambda}={\rm essinf}_{x,z\in\R^d, |z|>\lambda} b(x,z).$

\begin{prp}\label{x1}
For every $A>0, \lambda>0, 0<\ee<1$ and $R_0>0$,
there exists a constant $C_{16}=C_{16}(d,\alpha,\beta, A, \lambda, \ee, R_0)>0$
so that for every  $b$ satisfying \eqref{e:b} and \eqref{e:1.21a}
with $\| b\|_\infty \leq A$,
$r\in (0, R_0]$ and $x, y\in \R^d$ with $|x-y|\geq 3r$,
$$
\P_x \big(\sigma^b_{B(y, r)}<\kappa r^{\alpha}\big)
\geq C_{16}\, r^{d+\alpha} \, \left(\frac{1}{|x-y|^{d+\alpha}}+
\frac{m_{b^+, \lambda} }{  |x-y|^{d+\beta}}\right).
$$
\end{prp}

\pf  By Proposition \ref{y1} ,
$$
\E_x \left[ \kappa r^{\alpha}\wedge \tau^b_{B(x,r)} \right]
\geq \kappa r^{\alpha}\, \P_x \big(\tau^b_{B(x,r)}\geq
\kappa r^{\alpha} \big)\geq \frac{1}{2}\kappa r^{\alpha}.
$$
Since $J^b(x, y)\geq m_{b^+, \lambda} \,
{\cal A}(d, -\beta)  |x-y|^{-d-\beta}\1_{\{|x-y|>\lambda\}}$,
\eqref{e:1.21a} implies that
$$
J^b(x,y)\geq \frac12 \left( \ee|x-y|^{-(d+\alpha)}+
m_{b^+, \lambda}\, {\cal A}(d, -\beta)  |x-y|^{-(d+\beta)}
\1_{\{|x-y|>\lambda\}} \right)
$$
Thus by Proposition \ref{p4}, there are positive constants
$c_1= c_1(d,\alpha,\beta )$ and
\hfill \break
$c_2= c_2(d,\alpha,\beta, A, \lambda, \ee, R_0)$
so that
\begin{eqnarray*}
\P_x(\sigma^b_{B(y, r)}<\kappa r^{\alpha})
  &\geq&
\P_x(X^b_{\kappa r^{\alpha}\wedge \tau^b_{B(x,r)}}\in B(y,r))\\
&=& \E_x\int_0^{\kappa r^{\alpha}\wedge\tau^b_{B(x,r)}}\int_{B(y,r)}
J^b(X^b_s, u)\,du\,ds\\
&\geq & c_1   \E_x \big[\kappa r^{\alpha}\wedge
\tau^b_{B(x,r)}\big]\int_{B(y,r)} \left(\frac{\eps}{|x-y|^{d+\alpha}}
 +\frac{m_{b^+, \lambda} }{  |x-y|^{d+\beta}}  \1_{\{|x-y|>\lambda\}}\right) du\\
&\geq&  c_2 \ee\kappa r^{d+\alpha} \, \left(\frac{1}{|x-y|^{d+\alpha}}+
\frac{m_{b^+, \lambda} }{  |x-y|^{d+\beta}}\right).
\end{eqnarray*}
Here in the last inequality, we used the fact that $|x-y|^{-(d+\alpha)}\geq
(1+\lambda^{\alpha-\beta} A)^{-1}[|x-y|^{-(d+\alpha)}+m_{b, \lambda} \cdot|x-y|^{-(d+\beta)}]$
for $|x-y|\leq \lambda.$
\qed

\begin{prp}\label{x1'}
For every $A>0$,
there exists a constant $C_{17}=C_{17}(d,\alpha,\beta, A)>0$
so that for every bounded $b$ that  satisfies \eqref{e:b} and \eqref{e:1.15}
with $\| b\|_\infty \leq A$,
 and  $3r\leq |x-y|\leq R_\ast
 :=\frac{1}{3} \left( 2A\frac{{\cal A} (d, -\beta)}{{\cal A} (d, -\alpha)}\right)^{1/(\beta-\alpha)}$,
$$
\P_x \big(\sigma^b_{B(y, r)}<\kappa r^{\alpha}\big)
\geq C_{17}\, \frac{r^{d+\alpha}}{|x-y|^{d+\alpha}}.
$$
\end{prp}

\pf Note that when $|x-u|\leq 3R_\ast$,
$\frac{1}{2}\frac{{\cal A} (d, -\alpha)}{{\cal A} (d, -\beta)}|x-u|^{\beta-\alpha}\geq A$ and so

\begin{equation}\label{a0}\begin{aligned}
J^b(x,u)&= \, \frac{{\cal A}(d, -\alpha)}{|x-u|^{d+\alpha}}
+ \frac{{\cal A}(d, -\beta) \, b(x, u-x)}{|x-u|^{d+\beta}}\\
&\geq \, \frac{{\cal A}(d, -\alpha)}{|x-u|^{d+\alpha}}
-A\frac{{\cal A}(d, -\beta)}{|x-u|^{d+\beta}}
\, \geq \, \frac{1}{2}\frac{{\cal A}(d,-\alpha)}{|x-u|^{d+\alpha}}.
\end{aligned}\end{equation}
By Propositions \ref{p4} and \ref{y1},
we have
$$\begin{aligned}
\P_x \left(\sigma^b_{B(y, r)}<\kappa r^{\alpha} \right)
&\geq \P_x \left( X^b_{\kappa r^{\alpha}\wedge \tau^b_{B(x,r)}}
 \in B(y,r) \right)\\
&=\E_x\int_0^{\kappa r^{\alpha}\wedge\tau^b_{B(x,r)}}\int_{B(y,r)}
J^b(X^b_s, u)\,du\,ds\\
&\geq c_1\E_x \left[ \kappa r^{\alpha}\wedge
\tau^b_{B(x,r)} \right] \int_{B(y,r)}
\frac{1}{|x-y|^{d+\alpha}} du\\
&\geq c_2\kappa r^{d+\alpha} \, \frac{1}{|x-y|^{d+\alpha}},
\end{aligned}$$
where the second inequality holds due to $|X^b_s-u|\leq 3 |x-y|\leq 3R_\ast$ and (\ref{a0}).
\qed

\begin{thm}\label{t2'}
For every $\lambda>0, \ee\in (0,1)$ and $A>0$, there are positive constants
$C_{18}=C_{18} (d,\alpha,\beta, A, \lambda, \ee)$ and $C_{19}=C_{19}
(d,\alpha,\beta, A, \lambda)$
 such that
for any  $b$ with $\| b\|_\infty\leq A$
that satisfies \eqref{e:b} and \eqref{e:1.21a},
\begin{equation}\label{e2'}
C_{18} \, p_{ \, m_{b^+, \lambda}} (t,x,y)\leq q^b(t,x,y)
\leq C_{19} \,  {p}_{M_{b^+, \lambda}} (t,x,y), \quad t\in
(0,1],\, x,y\in\R^d.
\end{equation}
\end{thm}

\pf Noting that the condition \eqref{e:1.21a} in particular
implies \eqref{e:1.15}, so the upper bound estimate follows immediately
from Theorem \ref{T:u1}. We only need to prove the lower bound.
Let $\delta_0:=1\wedge (A_0/A)^{\alpha/(\alpha-\beta)}.$
\eqref{e:ql1} together with (\ref{e:1.5})  also yields that for
any  $\|b\|_\infty\leq A$,
\begin{equation}\label{e:lq^b}
q^b(t,x,y)\geq c_0 t^{-d/\alpha}  \quad \hbox{for }
t\in(0,\delta_0] \hbox{ and } \,|x-y|\leq 3t^{1/\alpha}.
\end{equation}
Here $c_0=c_0(d,\alpha, \beta)$ is a positive constant.
 For every $t\in (0, \delta_0]$,
by Proposition \ref{y1} and Proposition
\ref{x1} with  $R_0=1$,
$r=t^{1/\alpha}/2$ and the strong Markov property of the process $X^b,$ we get for $|x-y|>3t^{1/\alpha},$
\begin{eqnarray}\label{P}
&& \P_x(X^b_{2^{-\alpha}\kappa t}\in B(y, t^{1/\alpha})) \nonumber \\
&\geq& \P_x\left( X^b \hbox{ hits} B(y, t^{1/\alpha}/2) \hbox{ before }
2^{-\alpha} \kappa t \hbox{ and stays there for at least }
 2^{-\alpha} \kappa t \hbox{ units of time} \right) \nonumber \\
&\geq & \P_x \left(\sigma^b_{B(y, t^{1/\alpha}/2)}<2^{-\alpha}\kappa t \right) \inf_{z\in B(y, t^{1/\alpha}/2)}
\P_z \left( \tau^b_{B(y, t^{1/\alpha})}\geq 2^{-\alpha}\kappa t\right) \nonumber \\
&\geq & \P_x \left( \sigma^b_{B(y, t^{1/\alpha}/2)}<2^{-\alpha}\kappa t \right) \inf_{z\in B(y, t^{1/\alpha}/2)}
\P_z \left( \tau^b_{B(z, t^{1/\alpha}/2)}\geq 2^{-\alpha}\kappa t \right) \nonumber \\
&\geq&  c_1\, t^{(d+\alpha)/\alpha}
\left( \frac{1}{|x-y|^{d+\alpha}} +
\frac{m_{b^+, \lambda}}{|x-y|^{d+\beta}}\right).
\end{eqnarray}
Here $c_1=c_1(d,\alpha,\beta, A, \lambda, \ee)$ is a positive constant.
 Hence,  by (\ref{e:lq^b}) and (\ref{P}), for $|x-y|>3t^{1/\alpha}$ and $t\in (0,\delta_0],$
\begin{equation}\label{e:5.ql}\begin{aligned}
q^b(t,x,y)&\geq \int_{B(y,t^{1/\alpha})} q^b(2^{-\alpha}\kappa t,x,z)q^b((1-2^{-\alpha}\kappa) t, z, y)\,dz\\
&\geq \inf_{z\in B(y,t^{1/\alpha})}q^b((1-2^{-\alpha}\kappa)t,z,y)\P_x(X^b_{2^{-\alpha}\kappa t}\in B(y,t^{1/\alpha}))\\
&\geq c_2t^{-d/\alpha} \, t^{(d+\alpha)/\alpha}
\left( \frac{1}{|x-y|^{d+\alpha}} +
 \frac{m_{b^+, \lambda}}{|x-y|^{d+\beta}}\right)\\
&\geq c_2
\left( \frac{t}{|x-y|^{d+\alpha}} +
\frac{t \, m_{b^+, \lambda}}{ |x-y|^{d+\beta}}\right) ,
\end{aligned}\end{equation}
where $c_2=c_2(d,\alpha,\beta, A, \lambda, \ee)>0,$ the third inequality holds due to  $|z-y|\leq t^{1/\alpha}\leq 3((1-2^{-\alpha}\kappa)t)^{1/\alpha}$
when $\kappa\leq 2^{\alpha}(1-3^{-\alpha})$ and (\ref{e:lq^b})-(\ref{P}).
Finally, \eqref{e:lq^b}, \eqref{e:5.ql} together with \eqref{e:1.7a} and the Chapman-Kolmogorov equation yields the desired lower
bound estimate. \qed

\begin{thm}\label{t2''}
For every $\lambda>0$ and $A>0$, there are positive constants $C_k=C_k(d,\alpha,\beta, A),
k=20, 21$ and $C_{22}=C_{22}(d,\alpha,\beta, A, \lambda)$
such that for any  bounded $b$
satisfying \eqref{e:b} and \eqref{e:1.15} with $\| b\|_\infty\leq A$,
\begin{equation}\label{e:2''}
C_{20}\overline{p}_0(t,C_{21} x, C_{21} y)\leq q^b(t,x,y)
\leq C_{22} {p}_{M_{b^+, \lambda}} (t,x,y)  \quad \hbox{for }
t\in (0,1] \hbox{ and }  x,y\in\R^d.
\end{equation}
\end{thm}

\pf
By Theorem \ref{T:u1}, it suffices to prove the lower bound of $q^b$.
Let $\delta_0:=1\wedge (A_0/A)^{\alpha/(\alpha-\beta)}.$
By Chapman-Kolmogorov equation, we only need to consider \eqref{e:2''} for $t\in (0,\delta_0].$
By \eqref{e:trun1}, \eqref{e:trun2} and \eqref{e:ql1}, it suffices to prove \eqref{e:2''} when $|x-y|>3t^{1/\alpha}$ and $t\in (0,\delta_0].$
 Let $R_\ast$ be the constant defined in Proposition \ref{x1'}.

(i) First, we consider the case $R_\ast\geq |x-y|>3t^{1/\alpha}.$
For every $t\in (0, \delta_0]$,
by Proposition \ref{y1} and Proposition
\ref{x1'} with
$r=t^{1/\alpha}/2$ and the strong Markov property of the process $X^b,$
we get, by the similar procedure in (\ref{P}),
for $R_\ast\geq |x-y|>3t^{1/\alpha},$
\begin{equation}\label{P'}
\P_x \left( X^b_{2^{-\alpha}\kappa t}\in B(y, t^{1/\alpha}) \right)
\geq c_1\, t^{(d+\alpha)/\alpha}\frac{1}{|x-y|^{d+\alpha}} .
\end{equation}
Here $c_1=c_1(d,\alpha,\beta, A)$ is a positive constant.
 Hence, for $R_\ast\geq |x-y|>3t^{1/\alpha},$  by (\ref{e:lq^b}) and (\ref{P'}), we have
\begin{equation}\label{l'b1}\begin{aligned}
q^b(t,x,y)&\geq \int_{B(y,t^{1/\alpha})} q^b(2^{-\alpha}\kappa t,x,z)q^b((1-2^{-\alpha}\kappa) t, z, y)\,dz\\
&\geq \inf_{z\in B(y,t^{1/\alpha})}q^b((1-2^{-\alpha}\kappa)t,z,y)
\P_x \left( X^b_{2^{-\alpha}\kappa t}\in B(y,t^{1/\alpha}) \right)\\
&\geq c_2t^{-d/\alpha} \, t^{(d+\alpha)/\alpha}
\frac{1}{|x-y|^{d+\alpha}}\\
&\geq c_2 \frac{t}{|x-y|^{d+\alpha}}
\end{aligned}\end{equation}
where $c_2=c_2(d,\alpha,\beta, A)>0.$

(ii) Next, we consider the case $ |x-y|>R_\ast>3t^{1/\alpha}.$
 Take $C_\ast=R_\ast^{-1}.$ Then $|x-y|>R_\ast=C_\ast^{-1}\geq t/C_\ast$ for $t\in (0,\delta_0].$
The following proof is similar to \cite[Theorem 3.6]{CKK}. For the reader's  convenience, we
spell out the details here.

Let $R:=|x-y|$ and $c_+=R_\ast^{-1}\vee 1.$ Let $l\geq 2$ be a positive integer such that
$c_+R\leq l\leq c_{+}R+1$ and let $x=x_0, x_1, \cdots, x_l=y$ be such that $|x_i-x_{i-1}|\asymp R/l \asymp 1/c_+$
for $i=1,\cdots, l-1.$ Since $t/l\leq C_\ast R/l\leq C_\ast/c_+\leq 1$ and $R/l\leq 1/c_+\leq R_\ast,$ we have by  \eqref{e:lq^b} and \eqref{l'b1},
\begin{equation}\label{l'b}
q^b(t/l, x_i, x_{i+1})\geq c_2\left((t/l)^{-d/\alpha}\wedge\frac{t/l}{(R/l)^{d+\alpha}}\right)
\geq c_2\left((t/l)^{-d/\alpha}\wedge (t/l)\right)\geq c_3t/l.
\end{equation}
Let $B_i=B(x_i,R_\ast),$ by (\ref{l'b}),
\begin{equation}\label{l'b2}\begin{aligned}
q^b(t,x,y)&\geq \int_{B_1}\cdots\int_{B_{l-1}} q^b(t/l,x,x_1)\cdots q^b(t/l,x_{l-1}, y)\,dx_1\cdots dx_{l-1}\\
&\geq (c_4t/l)^l\geq (c_5t/R)^{c_+R+1}\geq c_6(t/R)^{c_7R}\\
&\geq c_6\left(\frac{t}{|x-y|}\right)^{c_7|x-y|}.
\end{aligned}\end{equation}
By \eqref{l'b1}, \eqref{l'b2} and together with the estimates of $\overline{p}_0$ in \eqref{e:trun1}-\eqref{e:trun2},
we get the desired conclusion.
\qed

 \bigskip
 \noindent{\bf Proof of Theorem \ref{T:1.3}.}
Theorem \ref{T:1.3} now follows from Theorems \ref{T:m1}, \ref{t2'} and \ref{t2''}.
\qed

\bigskip
To prove theorem \ref{T:1.4}, we use the main result  in \cite{CKS2} of the heat kernel estimates for non-local operators
under the non-local Feynman-Kac perturbation. For each Borel function $q(x)$ on $\R^d$ and Borel function $F(x,y)$ on $\R^d\times\R^d$
that vanishes along the diagonal, we define a non-local Feynman-Kac transform for the process $X^b$  as follows:
\begin{equation}\label{e:sg}
T^{b,F}_tf(x)=\E_x\left[\exp\left(\int_0^t q(X^b_s) \,ds + \sum_{s\leq t} F(X^b_{s-}, X^b_s)\right)f(X^b_t)\right].
\end{equation}

\begin{prp}\label{P:1}
Suppose $b$ is a bounded function on $\R^d\times\R^d$
satisfying \eqref{e:b} and \eqref{e:1.15},
$q$ is a bounded function on $\R^d$ and $|F(x,y)|\leq c(|x-y|^2\wedge 1)$ for some constant $c.$
Then for each $f$ in $C_b^2(\R^d),$
$$T^{b,F}_tf(x)=f(x)+\int_0^t T^{b,F}_s\L^{b,F}f(x)\,ds,$$
where $$\L^{b,F}f(x)=\L^b f(x)+\int_{\R^d} (e^{F(x,y)}-1) f(y) J^b(x,y)\,dy +q(x)f(x).$$
\end{prp}

\pf First note that since $X^b$ is a semimartingale and $|F(x,y)|\leq c(|x-y|^2\wedge 1)$,
$$
\sum_{s\leq t}|F(X^b_{s-}, X^b_s)|\leq c\sum_{s\leq t} | X^b_s -X^b_{s-}|^2=c \big[ X^b, X^b \big]_t<\infty.
$$
Let $F_1 =e^F-1$ and  define
$$K_t=\int_0^t q(X^b_s) \,ds +\sum_{s\leq t}F_1(X^b_{s-}, X^b_s).$$
Then by \cite[A4.17]{Sh}, the Stieljes exponential
\begin{equation}\label{e:K}
A_t:= \hbox{Exp}(K)_t=e^{K_t^c}\prod_{0<s\leq t}(1+  K_s-K_{s-})=\exp\left(\int_0^t q(X^b_s) \,ds + \sum_{s\leq t} F(X^b_{s-}, X^b_s)\right)
\end{equation}
is the unique solution to
\begin{equation}\label{e:A}
A_t=1+\int_0^t A_{s-} \,dK_s.
\end{equation}

For each function $f$ in $C_b^2(\R^d),$ by Ito's formula, Proposition \ref{p2} and \eqref{e:A}, we have
\begin{eqnarray}\label{e:I1}
 A_tf(X^b_t)
&=&f(X^b_0) + \int_0^t f(X^b_{s-})A_{s-}\,dK_s+\int_0^t A_{s-} \,df(X^b_s)+\sum_{s\leq t} (A_s-A_{s-}) ( f(X^b_s)-f(X^b_{s-}))\nonumber \\
&=& f(X^b_0) + \int_0^t A_{s}f(X^b_{s})q(X^b_s)\,ds+ \sum_{s\leq t} f(X^b_{s-})A_{s-} F_1(X^b_{s-}, X^b_s)\\
&&\quad +\int_0^t A_{s}\L^bf(X^b_s) \,ds+ \int_0^t A_{s-} \,dM^f_s+\sum_{s\leq t} A_{s-}F_1(X^b_{s-}, X^b_s)( f(X^b_s)-f(X^b_{s-})) . \nonumber
\end{eqnarray}
By taking expectation on both sides and using the L\'evy system formula in Proposition \ref{p4}, we get
\begin{eqnarray*}
T^{b,F}_t f(x)&=&\E_x \big[A_tf(X^b_t)\big]\\
&=&f(X^b_0)+\E_x\left[\int_0^t A_{s}f(X^b_{s})q(X^b_s)\,ds+\int_0^t A_{s}\L^bf(X^b_s) \,ds\right]\\
&& + \, \E_x\left[\int_0^t\int_{\R^d} A_s \left(f(X^b_s)+ \big( f(y)-f(X^b_s) \big)\right)
F_1(X^b_{s}, y)J^b(X^b_s, y)\,dy\,ds\right]\\
&=&f(X^b_0)+\E_x\left[\int_0^t A_s \L^{b,F}f(X^b_s) ds \right]\\
&=&f(x)+\int_0^t T^{b,F}_s \L^{b,F} f(x) ds.
\end{eqnarray*}
That completes the proof.
\qed

\bigskip

\noindent{\bf Proof of Theorem \ref{T:1.4}.}
Let $b_0(x,z) =b(x,z) \1_{|z|\leq 1}(z)$, which
 is a bounded function on $\R^d\times\R^d$
satisfying \eqref{e:b} and \eqref{e:1.15}.
By Theorem \ref{T:1.3}, $q^{b_0}(t,x,y)$ is continuous on $(0,\infty)\times\R^d\times\R^d$ and
\begin{equation}\label{e:b01}
C_4p_0(t,x,y)\leq q^{b_0}(t,x,y)\leq C_3p_0(t,x,y)
\end{equation}
for all
$(t,x,y)\in (0,1]\times\R^d\times\R^d.$
In addition, by Proposition \ref{p4} and \eqref{e:Jb}, for each non-negative function $f$ on $\R^d\times\R^d$
that vanishes along the diagonal,
\begin{equation}\label{e:b02}
\E_x \Big[ \sum_{s\leq T}f( X^{b_0}_{s-}, X^{b_0}_s) \Big]
=\E_x \left[ \int_0^T\int_{\R^d}
f(X^{b_0}_s,u) J^{b_0}(X^{b_0}_s,u) \,du\,ds \right].
\end{equation}
and there exist two positive constants $c_1$ and $c_2$ so that
\begin{equation}\label{e:b03}
c_1|x-y|^{-(d+\alpha)}\leq J^{b_0}(x,y)\leq c_2|x-y|^{-(d+\alpha)}.
\end{equation}

Set $F(x,y)=\ln\frac{J^{b}(x,y)}{J^{b_0}(x,y)}$ and $q(x)=\int_{\R^d} (J^{b_0}(x,y)-J^b(x,y))\,dy.$
It is easy to see that $q$ is a bounded function on $\R^d$ and $J^b(x,y)=J^{b_0}(x,y)$ for $|x-y|\leq 1.$
Moreover, by the \eqref{e:1.21} and \eqref{e:b03},  there exist two positive constants $c_3$ and $c_4$ so that
$c_3\leq \frac{J^{b}(x,y)}{J^{b_0}(x,y)}\leq c_4$
for all $|x-y|>1$ and any bounded $b$ with $\|b\|_\infty\leq A.$
Hence, there is a positive constant $c_5$ so that $|F(x,y)|\leq c_5(|x-y|^2\wedge 1).$
Let $T^{b_0, F}_t$ be the semigroup $T^{b,F}_t$ defined by \eqref{e:sg}
but with $b_0$ in place of $b.$
By \eqref{e:b01}-\eqref{e:b03} above and \cite[Theorem 1.3]{CKS2}, the non-local Feynman-Kac semigroup $(T^{b_0,F}_t, t\geq 0)$ has a continuous
density $\wt{q}(t,x,y)$ and there is a positive constant $c_6$ so that for all $(t,x,y)\in (0,1]\times\R^d\times\R^d,$
\begin{equation}\label{e:b04}
c_6^{-1}p_0(t,x,y)\leq \wt{q}(t,x,y)\leq c_6p_0(t,x,y).
\end{equation}
On the other hand, for each $f$ in $C_b^2(\R^d),$
$$\begin{aligned}\L^{b_0,F}f(x)&=\L^{b_0} f(x)+\int_{\R^d} (e^{F(x,y)}-1) f(y) J^{b_0}(x,y)\,dy +q(x)f(x)\\
&=\L^{b_0}f(x)+\int_{\R^d} \left(J^{b}(x,y)-J^{b_0}(x,y)\right) (f(y)-f(x))\,dy \\
&=\L^b f(x).\end{aligned}$$
By taking $f=1$ in Proposition \ref{P:1}, we get  $T^{b_0, F}_t1=1$.
Hence $\wt{q}(t,x,y)$  uniquely determines a conservative
Feller process $\wt{Y}$ with   $\{T^{b_0, F}_t; t\geq 0\}$
as its transition semigroup. Proposition \ref{P:1}
implies that  the distribution  of $\wt{Y}$ on the canonical Skorokhod space ${\mathbb D}([0,\infty), \R^d)$ is a solution to
 the martingale problem $(\L^b, C^2_b(\R^d))$ and in particular
 to the martingale problem $(\L^b, {\cal S}(\R^d))$.
 However by Theorem \ref{T:1.3}, martingale solution to the operator $(\L^b, \S(\R^d))$ is unique. This yields that
$\wt{q}=q^b$ and so we  get the desired conclusion from \eqref{e:b04}.
\qed

\bigskip

  {\bf Acknowledgements.} Part of the main results
of this paper has been
presented at the workshop on ``Nonlocal operators: Analysis, Probability and Geometry and Applications", held at ZiF, Bielefeld, Germany from
July 9 to July 14, 2012 and at the ``Eighth Workshop on Markov Processes
and Related Topics" held at Beijing Normal University and Wuyi Shanzhuang
from July 16 to July 21, 2012. Helpful comments from the audience, in particular those from Mufa Chen, Mateusz Kwasnicki,  and
Ting Yang, are gratefully acknowledged.

\vskip 0.3truein

{\bf Zhen-Qing Chen}

Department of Mathematics, University of Washington, Seattle,
WA 98195, USA

E-mail: \texttt{zqchen@uw.edu}

\bigskip

{\bf Jie-Ming Wang}

School of Mathematics and Statistics,  Beijing Institute of Technology,
Beijing 100081, P. R. China.

E-mail: \texttt{wangjm@bit.edu.cn}

\begin{thebibliography}{99}


\bibitem{D} D. Applebaum,
 \emph{L\'evy Processes and Stochastic Calculus}.
 Cambridge University Press, 2004.

\bibitem{BBCK} M. T. Barlow, R. F. Bass, Z.-Q. Chen and M. Kassmann,
 Non-local Dirichlet forms and symmetric jump processes.
{\it Trans. Amer. Math. Soc. \bf 361} (2009), 1963-1999.


\bibitem{B} R. F. Bass,
Local times for a class of purely discontinuous
martingales. {\it Z. Wahrsch. Verw. Gebiete \bf 67} (1984),
433-459.



\bibitem{BC} R. F. Bass and Z.-Q. Chen,
System of equations driven by stable processes.
{\it Probab. Theory Relat. Fields \bf 134} (2006),
175-214.

\bibitem{BT} R. F. Bass and H. Tang,
 The martingale problem for a class of stable-like processes.
 {\it Stochastic Process. Appl. \bf 119} (2009), 1144-1167.


\bibitem{BJ} K. Bogdan and T. Jakubowski,
Estimates of heat kernel of fractional
Laplacian perturbed by gradient operators.
{\it Comm. Math. Phys. \bf 271} (2007),
179-198.

\bibitem{C} Z.-Q. Chen,
Symmetric jump processes and their heat kernel estimates.
{\it Sci. China Ser. A. \bf 52} (2009), 1423-1445.

 

\bibitem{CKK} Z.-Q.  Chen, P. Kim and T. Kumagai, Weighted Poincar\'e
 inequality and heat kernel estimates for finite range jump processes.
{\it Math. Ann. \bf 342} (2008), 833-883.

\bibitem{CKS} Z.-Q.  Chen, P. Kim and R. Song,
Dirichlet heat kernel estimates for fractional Laplacian
under gradient perturbation.  {\it Ann. Probab. \bf 40} (2012),
2483-2538.

\bibitem{CKS2} Z.-Q.  Chen, P. Kim and R. Song,
Stability of Dirichlet heat kernel estimates for non-local operators under Feynman-Kac perturbation. 
{\it Trans. Amer. Math. Soc. \bf 367} (2015), 5237-5270.

\bibitem{CK} Z.-Q.  Chen and T. Kumagai,
  Heat kernel estimates for stable-like processes on $d$-sets.  {\it Stoch. Process Appl., \bf 108} (2003), 27-62.

\bibitem{CK2} Z.-Q.  Chen and  T. Kumagai,
Heat kernel estimates for jump processes of mixed types
on metric measure spaces. {\it Probab. Theory Relat. Fields, \bf 140}
(2008), 270-317.

\bibitem{CW} Z.-Q. Chen and L. Wang,
Uniqueness of stable processes with drift. 
{\it Proc. Amer. Math. Soc. \bf 144} (2016), 2661-2675.


\bibitem{CZ} Z.-Q. Chen and Z. Zhao,
Potential theory for elliptic systems.
{\it Ann. Probab. \bf 24} (1996), 293-319.

\bibitem{Chung} K. L. Chung,
{\it Lectures from Markov Processes to Brownian Motion}.
Springer-Verlag, Berlin and Heidelberg, 1982.

\bibitem{EK} S. N. Ethier and T. G. Kurtz,
\emph{Markov Processes: Characterization and Convergence}.
  Wiley, New York 1986.

\bibitem{Ja} N. Jacob, {\it Pseudo Differential Operators
and Markov Processes. I, II and III.} Imperial College Press,
London, 2001, 2002 and 2005.

\bibitem{K} V. Kolokoltsov, Symmetric stable laws and stable-like jump diffusions. 
{\it Proc. London
Math. Soc. \bf 80} (2000), 725-768.


\bibitem{Ko1} T. Komatsu, Markov processes associated with
certain integro-differential operators.
{\it Osaka J. Math. \bf 10} (1973), 271-303.

\bibitem{Ko2} T. Komatsu, On the martingale problem
for generators of stable processes with perturbations.
{\it Osaka J. Math. \bf 21} (1984), 113-132.

\bibitem{Mey} P.-A. Meyer, Renaissance, recollements, m\'{e}langes, raletissement
de processus de Markov. {\it Ann. Inst. Fourier \bf 25} (1975), 464-497.


\bibitem{MP} R. Mikulevicious and  G. Pragarauskas, On the martingale problem associated with
nondegenerate L\'evy operators. {\it Liet. Mat. Journal \bf 32} (1992), 297-311.


\bibitem{NT}
A. Negoro,A. and M. Tsuchiya,
 Stochastic processes and semigroups
  associate with degenerate \mbox{L}\'{e}vy generating operators.
  \emph{Stochastic and Stochastic Reports} \textbf{ {26}} (1989),  29--61.

\bibitem{SW} R. Schilling and J. Wang: Some theorems on Feller processes:
transience, local times and ultracontractivity. 
{\it Trans. Amer. Math. Soc. \bf 365} (2013),  3255-3286

\bibitem{Sh} M. Sharpe. {\it General theory of Markov processes}. Academic Press, Boston, 1988.

\bibitem{St} D. Stroock, Diffusion processes associated with
L\'evy generators. {\it Z. Wahrsch. Verw. Gebiete \bf 32} (1975), 209-244.

\bibitem{TTW} H. Tanaka, M. Tsuchiya and S. Watanabe,
Perturbation of drift-type for L\'evy processes.
{\it J. Math. Kyoto Univ. \bf 14} (1974), 73-92.

\bibitem{Ts} M. Tsuchiya, On some perturbation of stable processes.
{\it Proc. 2nd Japan-USSR Symposium on Probab. Theory}.
 Lect. Notes Math. {\bf 330} (1973), 490-497.

\bibitem{W} C. Wang. On estimates of the density of
Feynman-Kac semigroups of $\alpha$-stable-like processes. {\it J.
Math. Anal. Appl.} {\bf 348} (2008), 938--970.

\bibitem{Wa} J. -M. Wang, 
Laplacian perturbed by non-local operators. 
{\it Math.  Z.  \bf 279} (2015),  521-556.


\end{thebibliography}
\end{document}